%% file: cccsimply20200128.tex
\documentclass[12pt,a4paper]{amsart}

\usepackage{amssymb,color}
\usepackage[english]{babel}
\usepackage{verbatim,here}
\usepackage[T1]{fontenc}
\usepackage{floatflt,graphicx,graphics}
\usepackage{color,xcolor}
\usepackage{enumerate}
\usepackage{animate}
\usepackage{a4wide}
\usepackage{amsmath, amssymb}
\usepackage{amsthm}

\usepackage[autolinebreaks]{mcode}

\newcounter{minutes}
\divide\time by 60
\newcounter{hours}
\multiply\time by 60 \addtocounter{minutes}{-\time}
{\small
\curraddr{}
\email{mms.nasser@qu.edu.qa}
\curraddr{}
\email{vuorinen@utu.fi}
}

\dedicatory{}
\commby{}
\theoremstyle{plain}

\theoremstyle{definition}

\theoremstyle{remark}

\newtheorem{nonsec}[equation]{}

\numberwithin{equation}{section}

\newcommand{\beq}{\begin{equation}}
\newcommand{\eeq}{\end{equation}}
\newcommand{\ben}{\begin{enumerate}}
\newcommand{\een}{\end{enumerate}}
\newcommand{\bequu}{\begin{eqnarray*}}
\newcommand{\eequu}{\end{eqnarray*}}
\newcommand{\bequ}{\begin{eqnarray}}
\newcommand{\eequ}{\end{eqnarray}}

\renewcommand{\Im}{{ \rm Im}\,}
\renewcommand{\Re}{{ \rm Re}\,}

\newcommand{\D}{\mathbb{D}}

\usepackage{url}
\usepackage{color}

\newcommand{\CC}{{\mathbb C}}

\renewcommand{\i}{\mathrm{i}}
\newcommand{\bI}{{\bf I}}
\newcommand{\bM}{{\bf M}}
\newcommand{\bN}{{\bf N}}

\renewcommand{\thefootnote}{\number_style{footnote}}
\begin{document}

\def\thefootnote{}

\title[]{Conformal Invariants in Simply Connected Domains}

\author[M.M.S. Nasser]{Mohamed M S Nasser}
\address{Department of Mathematics, Statistics and Physics, Qatar University, P.O. Box 2713, Doha, Qatar.}
\author[M. Vuorinen]{Matti Vuorinen}
\address{Department of Mathematics and Statistics, University of Turku, Turku, Finland.}
\date{}

\begin{abstract}
We study numerical computation of several conformal invariants of simply connected domains in the complex plane including, the hyperbolic distance, the reduced modulus,  the harmonic measure, and the modulus of a quadrilateral. 
The method we use is based on the boundary integral equation with the generalized Neumann kernel. 
Several numerical examples are presented. We validate the performance and accuracy of our method by considering several model problems with known analytic solutions.
\end{abstract}

\maketitle

\footnotetext{\texttt{{\tiny File:~\jobname .tex, printed: \number\year-%
\number\month-\number\day, \thehours.\ifnum\theminutes<10{0}\fi\theminutes}}}
\makeatletter

\makeatother




\input{sec108.tex}

\input{sec208.tex}

\input{sec308.tex}

\input{sec408.tex}

\input{sec508.tex}

\input{sec608.tex}


%

\input{biblio108.tex}

\end{document}

%% file: sec108.tex
\section{Introduction} \label{section1}
\setcounter{equation}{0}

Classical function theory studies analytic functions and conformal maps defined in subdomains of the complex plane $\mathbb{C}\,.$ Most commonly, the domain of  definition of the functions is the unit disk
$\mathbb{D}=\{z \in  \mathbb{C} \,|\,  |z|<1\}\,.$ 
The powerful Riemann mapping theorem says that a given simply-connected domain $G$ with non-degenerate boundary can be mapped conformally onto the unit disk and it enables us to extend results originally proven for functions defined in the unit disk to the case when the domain of definition is a simply-connected domain. Therefore for the convenient analysis of distances and other metric characteristics
of sets it is natural to use conformally invariant distances and set functions.
This works fine in the cases when the Riemann mapping function is known explicitly, such as in the cases described in \cite{KoppStall}. Polygons form a large class of plane domains for which
the Riemann mapping function can be given in terms of the Schwarz-Christoffel formula which is
semi-explicit, it involves unknown accessory parameters. A numerical implementation of conformal
mapping methods based on the Schwarz-Christoffel formula is documented in \cite{dt} and
the Schwarz-Christoffel Toolbox~\cite{dr} is widely used to solve mapping problems in concrete applications. The  so called crowding phenomenon,
an inherent computational challenge in these mapping problems, is described in \cite[pp.20-21]{dt}, \cite[pp.75-77]{ps}. This phenomenon can be observed already in numerical
conformal mapping of rectangles onto half-plane when the quotient $m>1$ of the side lengths is large enough. The critical value of $m$ depends on the computer 
floating point arithmetic and for double precision arithmetics the critical value
lies in the range 
$ [10,20]\,$  \cite[pp.20-21]{dt}, \cite[pp.75-77]{ps}. 

We apply here the boundary integral equation method as developed in \cite{LSN17,Nas-Siam1,Nas-ETNA,Nas-cvee,Weg-Nas} to compute numerically conformal invariants such as the hyperbolic metric, harmonic measure, reduced modulus, and the modulus of a quadrilateral~\cite{du,garmar,ps,Vas02}. The cases considered here include, in particular, classes of domains to which the earlier methods do not seem to apply. 
Our methods are described in the respective sections of the papers, they are implemented in MATLAB, and the results are summarized graphics. We also give experimental error estimate in some simple cases.
We include some code snippets within the text to indicate implementation details. All the computer codes of our computations are available in the internet link \url{https://github.com/mmsnasser/ci-simply}. 

Section 2 reviews the boundary integral equation method for computing the conformal mapping from bounded and unbounded simply connected domains onto circular domains. This method will be applied in the remaining sections, sometimes together with auxiliary procedures. 
In Section 3 we use our method to compute the hyperbolic distance for several examples. Section 4 deals with the reduced modulus for bounded and unbounded simply connected domains. 
Section 5 deals with the harmonic measure for a simply connected polygonal domains.
In Section 6, we present an iterative method for numerical computation of the conformal mapping from a quadrilateral onto a rectangle. We also present an analytic example to illustrate the effect of crowding phenomenon on the accuracy of such mapping.
 

%% file: sec208.tex
\section{Conformal mappings of simply connected domains} \label{sec:cm}

In this section, we review a numerical method for computing the conformal mapping from bounded and unbounded simply connected domains onto the unit disk and the exterior unit disk, respectively. The method is based on the boundary integral equation with the generalized Neumann kernel~\cite{Nas-Siam1,Nas-ETNA,Nas-cmft15,Weg-Nas}.

\nonsec{\bf The generalized Neumann kernel.}\label{sc:gnk}
Let $G$ be a bounded or an unbounded simply connected domain bordered by a closed smooth Jordan curve $\Gamma=\partial G$. The orientation of the boundary $\Gamma$ is  counterclockwise when $G$ is bounded and clockwise when $G$ is unbounded.
We parametrize $\Gamma$ by a $2\pi$-periodic complex function $\eta(t)$, $t\in[0,2\pi]$. We assume that $\eta(t)$ is twice continuously differentiable with $\eta'(t)\ne0$ (the presented method can be applied also if the curve $\Gamma$ has a finite number of corner points but no cusps~\cite{Nas-cvee}).
We denote by $H$ the space of all H\"older continuous real functions on the boundary $\Gamma$.

Let $A$ be the complex function~\cite{Nas-ETNA}
\begin{equation}\label{eq:A}
A(t)= \left\{ \begin{array}{l@{\hspace{0.5cm}}l}
\eta(t)-\alpha,&{\rm if\;} G{\rm \; is\; bounded},\\
1 ,&{\rm if\;} G {\rm\; is\; unbounded}.\\
\end{array}
\right.
\end{equation}
The generalized Neumann kernel $N(s,t)$ is defined for $(s,t)\in [0,2\pi]\times [0,2\pi]$ by
\begin{equation}\label{eq:N}
N(s,t) :=
\frac{1}{\pi}\Im\left(\frac{A(s)}{A(t)}\frac{\dot\eta(t)}{\eta(t)-\eta(s)}\right).
\end{equation}
The kernel $N(s,t)$ is continuous~\cite{Weg-Nas}. Hence, the integral operator $\bN$ defined on $H$ by
\[
\bN\rho(s) := \int_J N(s,t) \rho(t) dt, \quad s\in [0,2\pi],
\]
is compact. The integral equation with the generalized Neumann kernel involves also the following kernel 
\begin{equation}\label{eq:M}
M(s,t) :=
\frac{1}{\pi}\Re\left(\frac{A(s)}{A(t)}\frac{\dot\eta(t)}{\eta(t)-\eta(s)}\right), \quad (s,t)\in [0,2\pi]\times [0,2\pi].
\end{equation}
which is singular and its singular part involves the cotangent function~\cite{Weg-Nas}. The integral operator $\bM$ defined on $H$ by
\[
\bM\rho(s) := \int_J  M(s,t) \rho(t) dt, \quad s\in [0,2\pi],
\]
is singular, but is bounded on $H$~\cite{Weg-Nas}.

\nonsec{\bf Bounded simply connected domain.}
Let $w=\Phi(z)$ be the unique conformal map of the bounded simply connected domain $G$ onto the unit disk $|w|<1$ such that
\begin{equation}\label{eq:cm-cond-b2}
\Phi(\alpha)=0\quad{\rm and}\quad\Phi'(\alpha)>0,
\end{equation}
Then, the mapping function $\Phi$ with normalization~\eqref{eq:cm-cond-b2} can be written for $z\in G\cup\Gamma$ as~\cite[\S~3]{Nas-cmft15}
\begin{equation}\label{eq:Phi-b}
\Phi(z)=c(z-\alpha)e^{(z-\alpha)f(z)}
\end{equation}
where the function $f(z)$ is analytic in $G$ with the boundary values $A(t)f(\eta(t))=\gamma(t)+h+\i\rho(t)$, $A(t)=\eta(t)-\alpha$, the function $\gamma$ is defined by $\gamma(t)=-\log|\eta(t)-\alpha|$, $\rho$ is the unique solution of the integral equation 
\begin{equation}\label{eq:ie}
(\bI-\bN)\rho=-\bM\gamma,
\end{equation}
and $c=e^{-h}$ where the constant $h$ is given by 
\begin{equation}\label{eq:h}
h=[\bM\rho-(\bI-\bN)\gamma]/2.
\end{equation}
Notice that $\Phi'(\alpha)=c=e^{-h}>0$.

Instead of the normalization~\eqref{eq:cm-cond-b2}, we can assume that the mapping function $\Phi$ satisfies the normalization
\begin{equation}\label{eq:cm-cond-b2-1}
\Phi(\alpha)=0\quad{\rm and}\quad\Phi'(\alpha)=1.
\end{equation}
In this case, the function $\Phi$ maps the domain $G$ onto the disk $|w|<R$ with a positive constant $R$. The constant $R$ is uniquely determined by $G$ and the point $\alpha$ and called the conformal radius of $G$ with respect to $\alpha$ and is denoted by $R(G,\alpha)$. For this case, in view of~\eqref{eq:Phi-b}, we can write the mapping function $\Phi$ for $z\in G\cup\Gamma$ as
\begin{equation}\label{eq:Phi-b-1}
\Phi(z)=(z-\alpha)e^{(z-\alpha)f(z)}
\end{equation}
where the function $f$ is as in~\eqref{eq:Phi-b}, i.e., we divide the right-hand side of~\eqref{eq:Phi-b} by $c=e^{-h}$. Hence, the boundary values of the mapping function $\Phi$ satisfy $|\Phi(\eta(t))|=1/c=e^h$ which implies
\begin{equation}\label{eq:cr-b}
R(G,\alpha)=e^h,
\end{equation}
where the constant $h$ is as in~\eqref{eq:h}.

\nonsec{\bf Unbounded simply connected domain.}
For an unbounded simply connected domain $G\subset\overline{\CC}$ with $\infty\in G$, there exists a unique conformal map  $w=\Phi(z)$ from $G$ onto the exterior of the unit disk $|w|>1$ such that
\begin{equation}\label{eq:cm-cond-u2}
\Phi(\infty)=0\quad{\rm and}\quad\Phi'(\infty)>0.
\end{equation}
Then, the mapping function $\Phi$ with normalization~\eqref{eq:cm-cond-b2} can be written for $z\in G\cup\Gamma$ as~\cite[\S~3]{Nas-cmft15}
\begin{equation}\label{eq:Phi-u}
\Phi(z)=c(z-\beta)e^{f(z)}
\end{equation}
where $\beta$ is an auxiliary point in $G^c=\overline{\CC}\backslash\overline{G}$ and $f(z)$ is an analytic function in $G$ with $f(\infty)=0$. The boundary values of the function $f$ are given by $A(t)f(\eta(t))=\gamma(t)+h+\i\rho(t)$ where $A(t)=1$, the function $\gamma$ is defined by $\gamma(t)=-\log|\eta(t)-\beta|$, $\rho$ is the unique solution of the integral equation~\eqref{eq:ie}, and $c=e^{-h}$ where the constant $h$ is given by~\eqref{eq:h}. Notice that $\Phi'(\infty)=c=e^{-h}>0$.

Instead of the normalization~\eqref{eq:cm-cond-u2}, we can assume that the mapping function $\Phi$ satisfies the normalization
\begin{equation}\label{eq:cm-cond-u2-1}
\Phi(\infty)=\infty\quad{\rm and}\quad\Phi'(\infty)=1.
\end{equation}
In this case, the function $\Phi$ maps the domain $G$ onto the disk $|w|>R$ with a positive constant $R$. The constant $R$ is uniquely determined by $G$ and called the conformal radius of $G$ with respect to $\infty$ and is denoted by $R(G,\infty)$. For this case, in view of~\eqref{eq:Phi-u}, we can write the mapping function $\Phi$ for $z\in G\cup\Gamma$ as
\begin{equation}\label{eq:Phi-u-1}
\Phi(z)=(z-\beta)e^{f(z)}
\end{equation}
where the function $f$ is as in~\eqref{eq:Phi-b}. Hence, the boundary values of the mapping function $\Phi$ satisfy $|\Phi(\eta(t))|=1/c=e^h$ which implies
\begin{equation}\label{eq:cr-u}
R(G,\infty)=e^h.
\end{equation}

\nonsec{\bf Numerical solution of the boundary integral equation.}
The integral equation~\eqref{eq:ie} is solved using the MATLAB function \verb|fbie| in~\cite{Nas-ETNA}. The function \verb|fbie| is based on using using the MATLAB function $\mathtt{zfmm2dpart}$ in the toolbox $\mathtt{FMMLIB2D}$~\cite{Gre-Gim12}. The computational cost for the overall method is $O(n\log n)$ operations where $n$ is the number of nodes in the interval $[0,2\pi]$. Let \verb|et|, \verb|etp|, \verb|A|, and \verb|gam| be the discretization vectors of the functions $\eta(t)$, $\eta'(t)$, $A(t)$, and $\gamma(t)$, respectively. Then discretization vectors $\verb|rho|$ and $\verb|h|$ of the solution $\rho(t)$ of the integral equation~\eqref{eq:ie} and the constant $h$ in~\eqref{eq:h}, respectively, can be computed by  
\[
[\verb|rho|,\verb|h|] = \verb|fbie|(\verb|et|,\verb|etp|,\verb|A|,\verb|gam|,\verb|n|,\verb|iprec|,\verb|restart|,\verb|gmrestol|,\verb|maxit|).
\]
In the numerical experiments in this paper, the parameters in $\verb|fbie|$ are chosen as following: $\mathtt{iprec}=5$ (the tolerances of the FMM is $0.5\times 10^{-15}$), $\mathtt{gmrestol}=0.5\times 10^{-14}$ (the tolerances of the GMRES), $\mathtt{restart}=[\,]$ (the GMRES is used without restart), and $\mathtt{maxit}=100$ (the maximum number of iterations for GMRES). 


Finally, the values of the auxiliary points $\alpha$ in~\eqref{eq:Phi-b},~\eqref{eq:Phi-b-1} and $\beta$ in~\eqref{eq:Phi-u},~\eqref{eq:Phi-u-1} have no effects on the values of the conformal mapping $\Phi$ as long as these points are sufficiently far away from the boundary $\Gamma$.

%% file: sec308.tex

\section{Hyperbolic distance}
\label{sc:hyp}

\nonsec{\bf Hyperbolic geometry.} \label{hypgeo} \cite{garmar,be, kela}
 For $x,y \in \mathbb{D}$ the {hyperbolic distance} $\rho_{\mathbb{D}}(x,y)$ is defined by
 \[
 \sinh \frac{\rho_{\mathbb{D}}(x,y)}{2} = \frac{|x-y|}{\sqrt{(1-|x|^2)(1-|y|^2)}}\,.
 \]
The main property of the hyperbolic distance is the invariance under  M\"obius transformations of $\D$ onto $\D$ defined by
\[
z \mapsto \frac{z-a}{1- \overline{a}z}\,
\] 
where $a\in\D$ is fixed.
In the metric space $(\mathbb{D}, \rho_{\mathbb{D}})$ one can build a non-euclidean geometry, where the parallel axiom does not hold. In this geometry, usually called the hyperbolic geometry of the Poincare disk, lines are circular arcs perpendicular to the boundary $\partial \mathbb{D} \,.$ This geometry is fundamentally different from the Euclidean geometry, but many results of Euclidean geometry have counterparts in the hyperbolic geometry \cite{be}.
 
 Let $G$ be a Jordan domain in the plane. One can define the hyperbolic metric on $G$ in terms of the  conformal Riemann mapping function $\Phi: G \to \mathbb{D}= \Phi(G)$ as follows:
 \[
 \rho_G(x,y) = \rho_{\mathbb{D}}(\Phi(x),\Phi(y))\,.
 \]
 This definition yields a well-defined metric, independent of the conformal mapping $\Phi\,$
 \cite{be,garmar,kela}. In hyperbolic geometry the boundary $\partial G$ has the same role as the point of $\{\infty\}$ in Euclidean geometry: both are like ``horizon's'', we cannot see beyond the
 horizon.
It turns out that the hyperbolic geometry is more useful than the Euclidean geometry when studying the inner geometry of domains in geometric function theory.

\nonsec{\bf Computing the hyperbolic distance for simply connected domains.}
Let $G\subset\overline{\CC}$ be a bounded simply connected domain, let $\alpha\in G$, and let $w=\Phi(z)$ be the unique conformal map of $G$ onto the unit disk $|w|<1$ with the normalization~\eqref{eq:cm-cond-b2}. Then for any two points $z_1$ and $z_2$ in $G$, we can define the hyperbolic metric $g$ in the usual way,
\begin{equation}\label{eq:hyp-dis}
\rho_G(z_1,z_2)=\rho_\D(\Phi(z_1),\Phi(z_2)) = 2\sinh^{-1}\left(\frac{|\Phi(z_1)-\Phi(z_2)|}{\sqrt{\left(1-|\Phi(z_1)|^2\right)\left(1-|\Phi(z_2)|^2\right)}}\right).
\end{equation}

A MATLAB function for computing the hyperbolic metric $\rho_G(z_1,z_2)$ for any two points $z_1$ and $z_2$ in the bounded simply connected domain $G$ is given below.

\begin{lstlisting}
function dis = hypdist (et,etp,n,alpha,zo,z)
% This function computes the hyperbolic distance dis between a point zo
% and a row vector of points z, in a simply connected domain G where:
% et, etp:  the parametrization of the boundary of G and its derivative
% n: the number of discretization points
% alpha: a given point in G
A      = et-alpha;
gam    =-log(abs(et-alpha));
[mu,h] = fbie(et,etp,A,gam,n,5,[],1e-14,200);
fet    =(gam+h+i*mu)./A;
Phi    = exp(-mean(h(1:n))).*(et-alpha).*exp(gam+h+i.*mu);
Phio   = fcau(et,etp,Phi,zo);
Phiz   = fcau(et,etp,Phi,z);
dis    = 2*asinh(abs(Phiz-Phio)./sqrt((1-abs(Phiz).^2).*(1-abs(Phio).^2)));
end
\end{lstlisting}

In the remaining part of this section, we use the MATLAB function \verb|hypdist| to compute the hyperbolic metric $g(z_1,z_2)$ for several examples. In these examples, for a given point $z_1$ in $G$, we define the function $u(x,y)$ for all $x$ and $y$ such that $x+\i y$ is in $G$ by
\begin{equation}\label{eq:htp-u}
u(x,y)=g(z_1,x+\i y).
\end{equation}
Then we use the MATLAB function \verb|hypdist| to compute the values of the function $u(x,y)$ in the domain $G$ and plot the  contour lines for the function $u(x,y)$ corresponding to the several levels. 

\begin{figure}[ht] %
\centerline{
\scalebox{0.5}[0.5]{\includegraphics[trim=0 -1.0cm 0 0,clip]{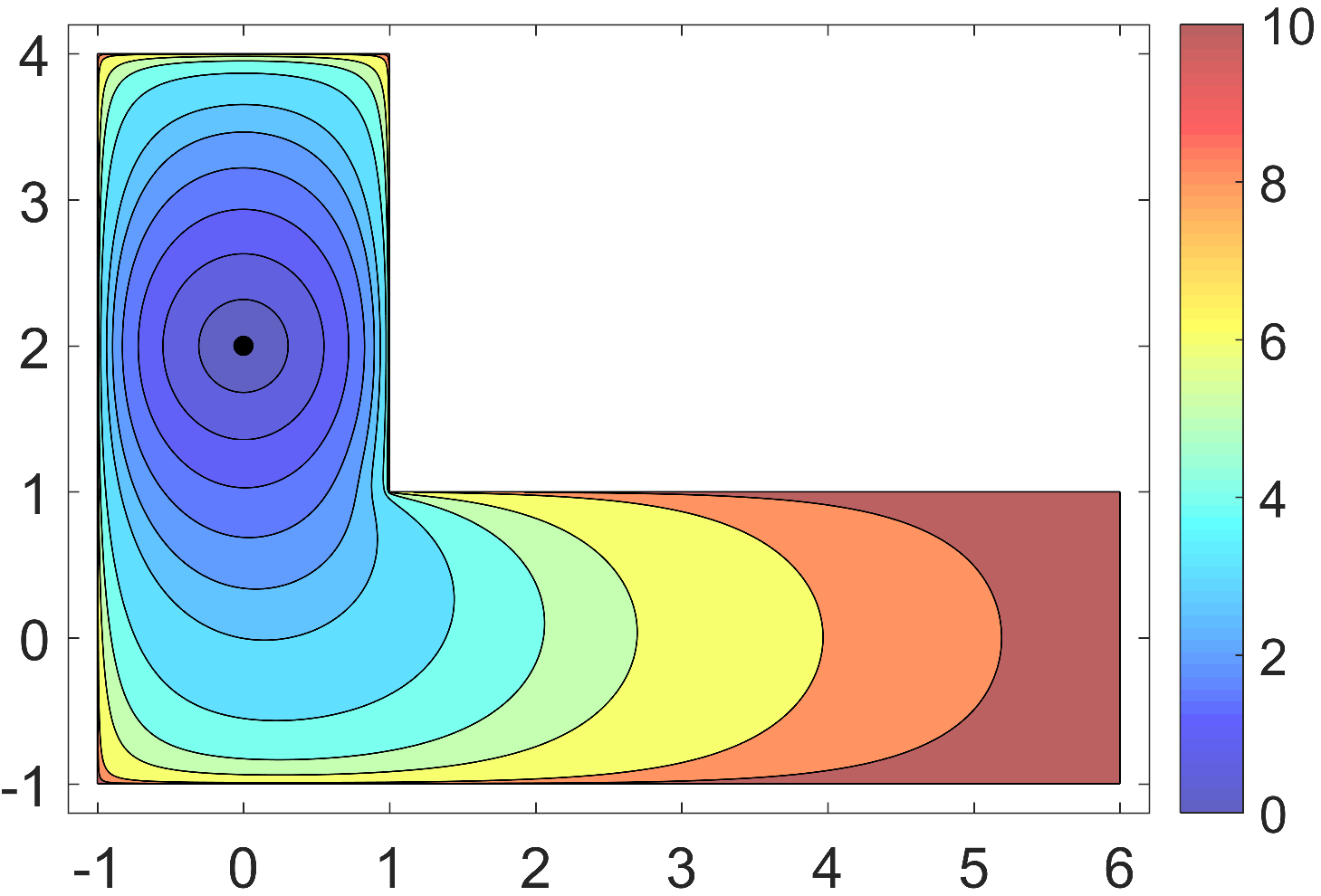}}
\hfill
\scalebox{0.5}[0.5]{\includegraphics[trim=0 -1.0cm 0 0,clip]{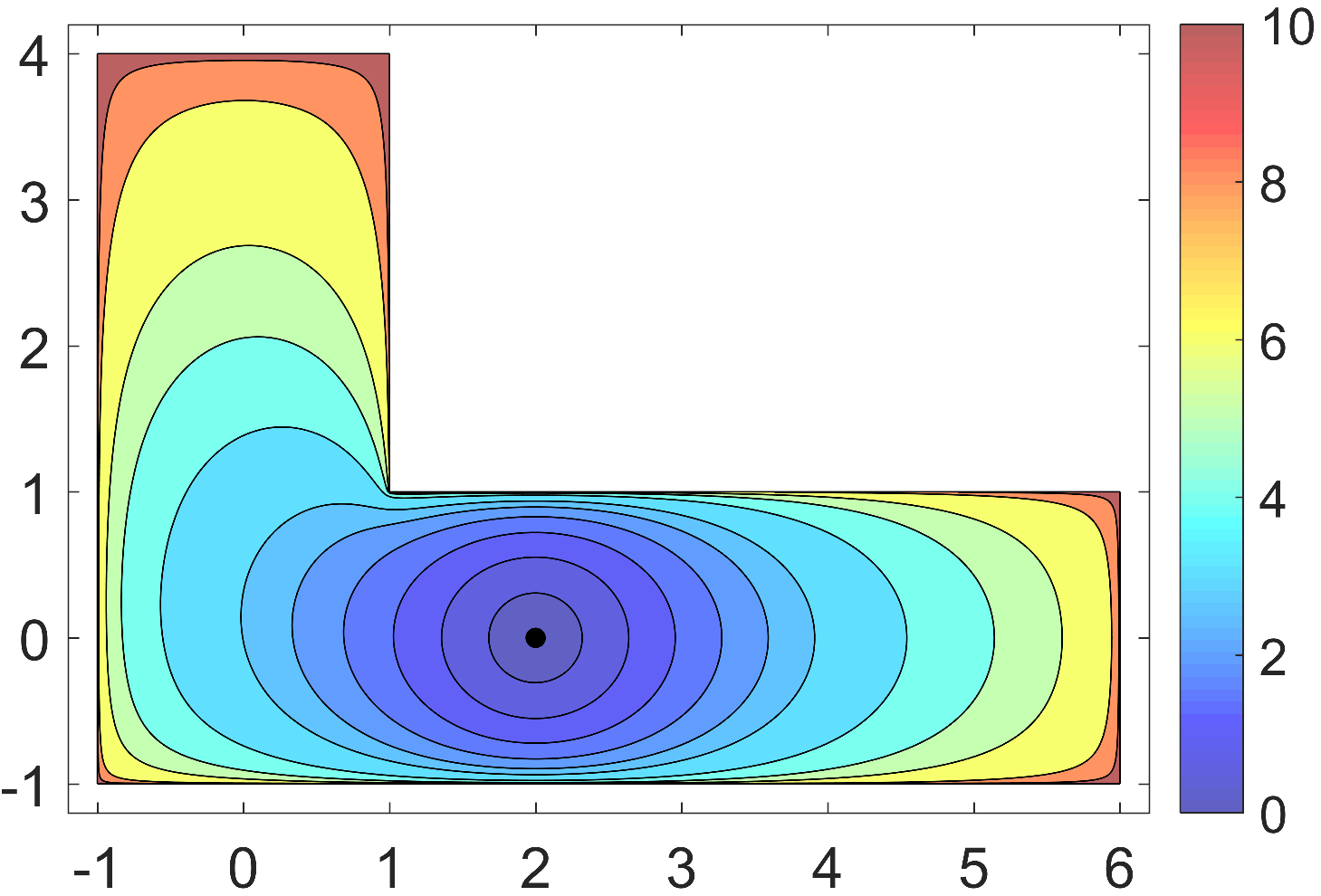}}
}
\caption{The contour lines of the function $u(x,y)$ for the L-shaped polygon for $z_1=2\i$ (left) and $z_1=2$  (right).}
\label{fig:hd-L}
\end{figure}

\nonsec{\bf L-shaped polygon.}
As our first example, we consider the simply connected domain $G$ in the interior of the L-shaped polygon with the vertices $6+\i$, $1+\i$, $1+4\i$, $-1+4\i$, $-1-\i$, and $6-\i$. 
The contour lines of the function $u(x,y)$ obtained with $n=6\times2^9$ discretization points on the boundary of the L-shaped polygon are shown in Figures~\ref{fig:hd-L} (left) for $z_1=2\i$ and in Figures~\ref{fig:hd-L} (right) for $z_1=2$. Table~\ref{tab:hyp-dist-L} presents the values of the hyperbolic metric $g(z_1,z_2)$ for $z_1=2\i$ (left), $z_1=2$ (right), and for several values of $z_2$.

\begin{table}[h]
\caption{The values of the hyperbolic metric $g(z_1,z_2)$.}
\label{tab:hyp-dist-L}%
\begin{tabular}{l|l||l|l}\hline
\multicolumn{2}{c||}{$z_1=2\i$}  & \multicolumn{2}{c}{$z_1=2$} \\ \hline
 $z_2$  & $g(z_1,z_2)$     & $z_2$  & $g(z_1,z_2)$\\ \hline
 1      & 3.50661554819086 & 0      & 2.99228771572299  \\
 2      & 4.91711064317017 & $\i$   & 3.50483278097652   \\
 3      & 6.47927360380709 & $2\i$  & 4.91711064317017  \\
 4      & 8.05147684115352 & $3\i$  & 6.52150321421451  \\
 5      & 9.66456147776192 &        &   \\
\hline
\end{tabular}
\end{table}

\begin{figure}[ht] %
\centerline{
\scalebox{0.5}[0.5]{\includegraphics[trim=0 -1.0cm 0 0,clip]{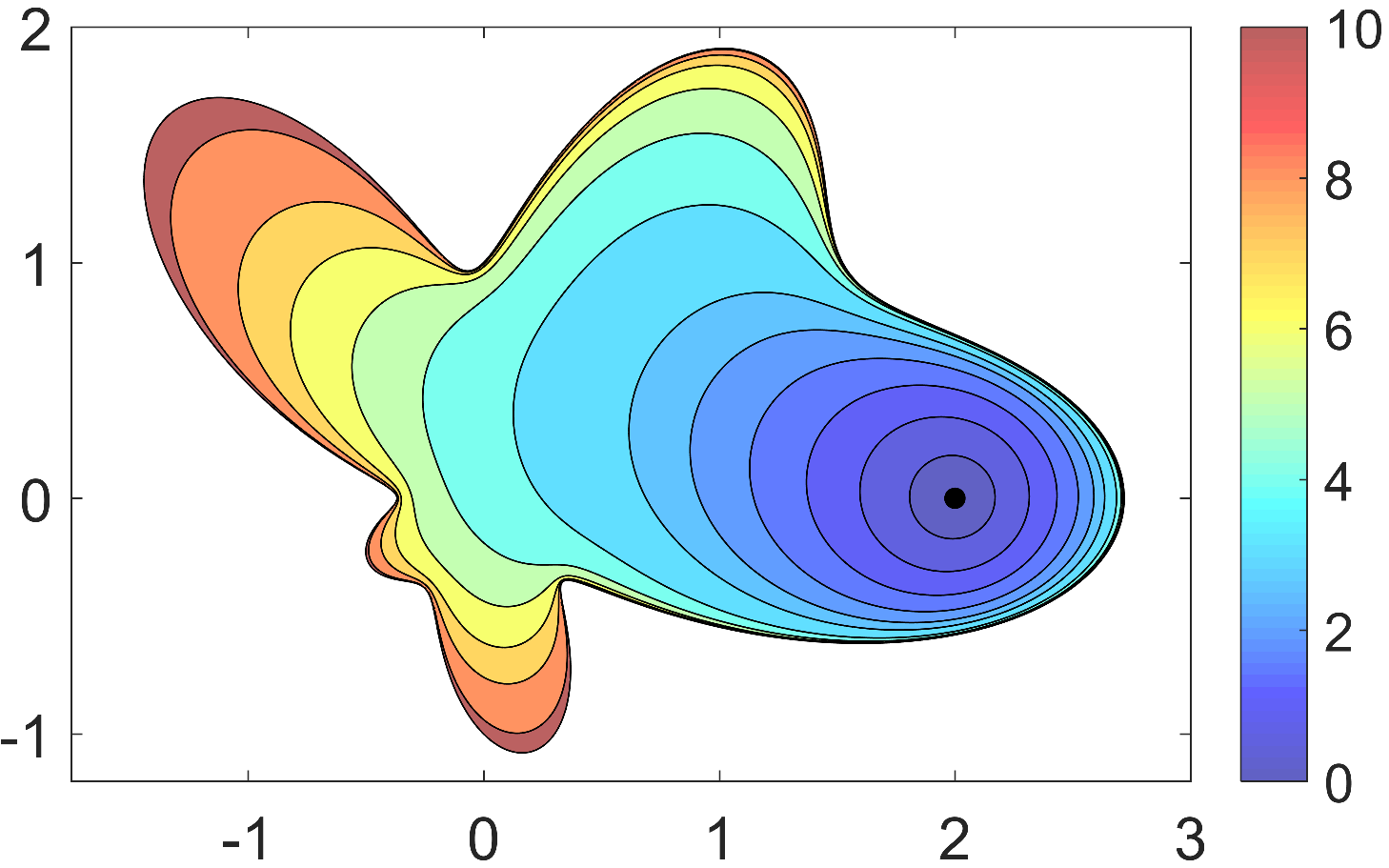}}
\hfill
\scalebox{0.5}[0.5]{\includegraphics[trim=0 -1.0cm 0 0,clip]{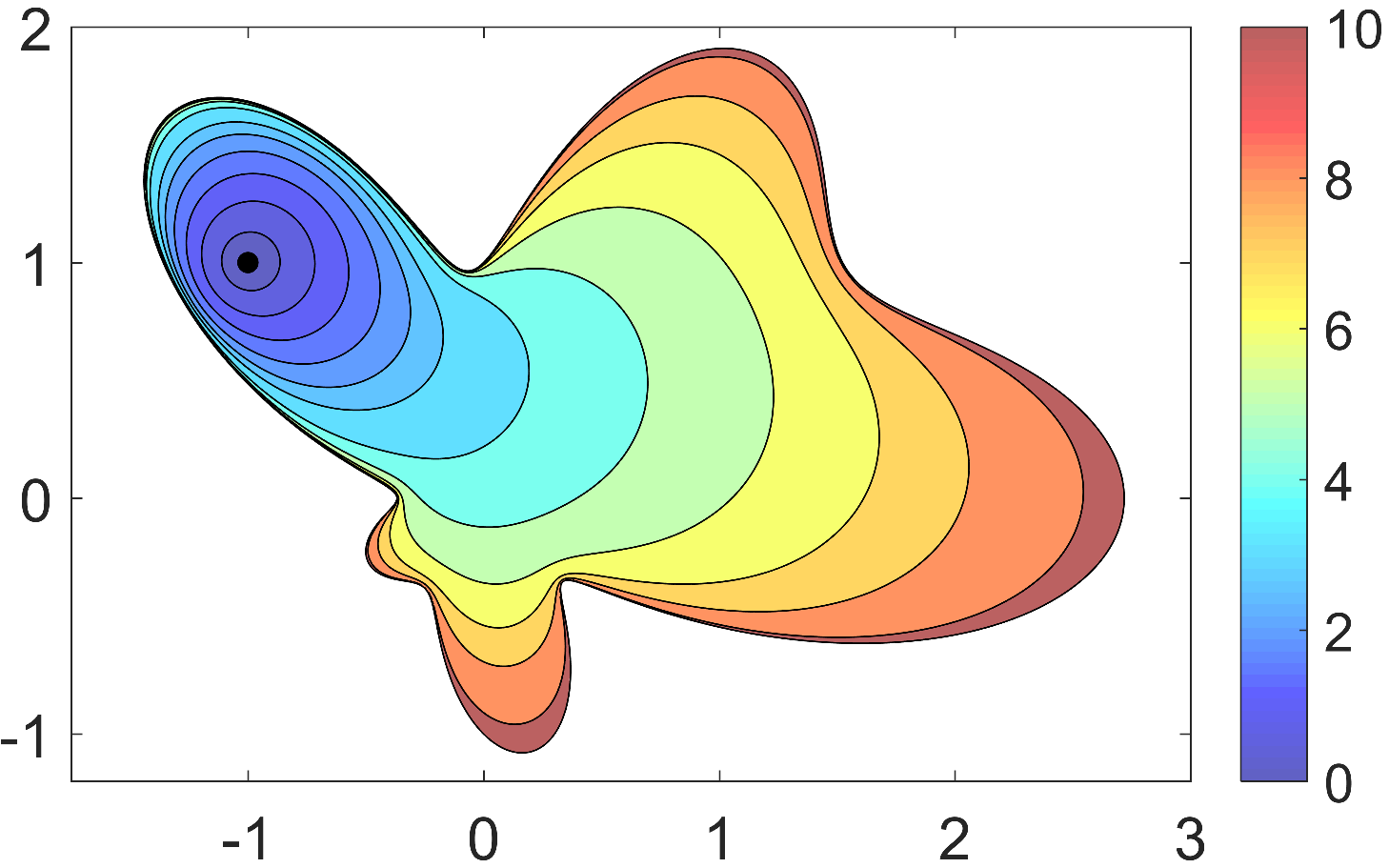}}
}
\caption{The contour lines of the function $u(x,y)$ for the amoeba-shaped boundary for $z_1=2$ (left) and $z_1=-1+\i$  (right).}
\label{fig:hd-s}
\end{figure}

\nonsec{\bf Amoeba-shaped boundary.}
In the second example, we consider the simply connected domain $G$ in the interior of the curve $\Gamma$ (amoeba-shaped boundary~\cite{AJ}) with the parametrization
\[
\eta(t)=\left(e^{\cos t}\cos^22t+e^{\sin t}\sin^22t\right)e^{\i t}, \quad 0\le t\le 2\pi.
\]
The contour lines for the function $u(x,y)$ computed with $n=2^{12}$ are shown in Figures~\ref{fig:hd-s} (left) for $z_1=2$ and in Figures~\ref{fig:hd-s} (right) for $z_1=-1+\i$.

\begin{figure}[ht] %
\centerline{
\scalebox{0.5}[0.5]{\includegraphics[trim=0 -1.0cm 0 0,clip]{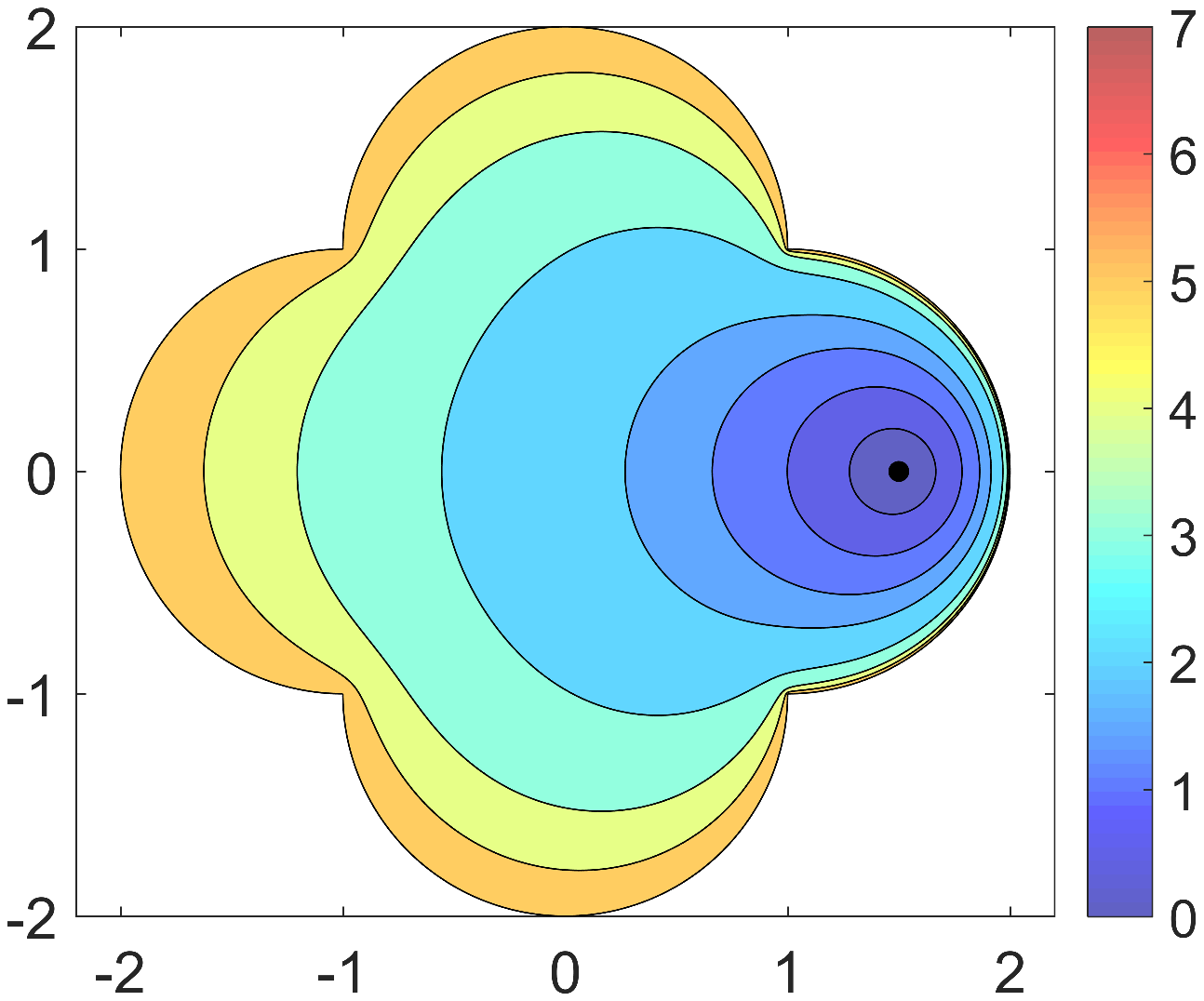}}
\hfill
\scalebox{0.5}[0.5]{\includegraphics[trim=0 -1.0cm 0 0,clip]{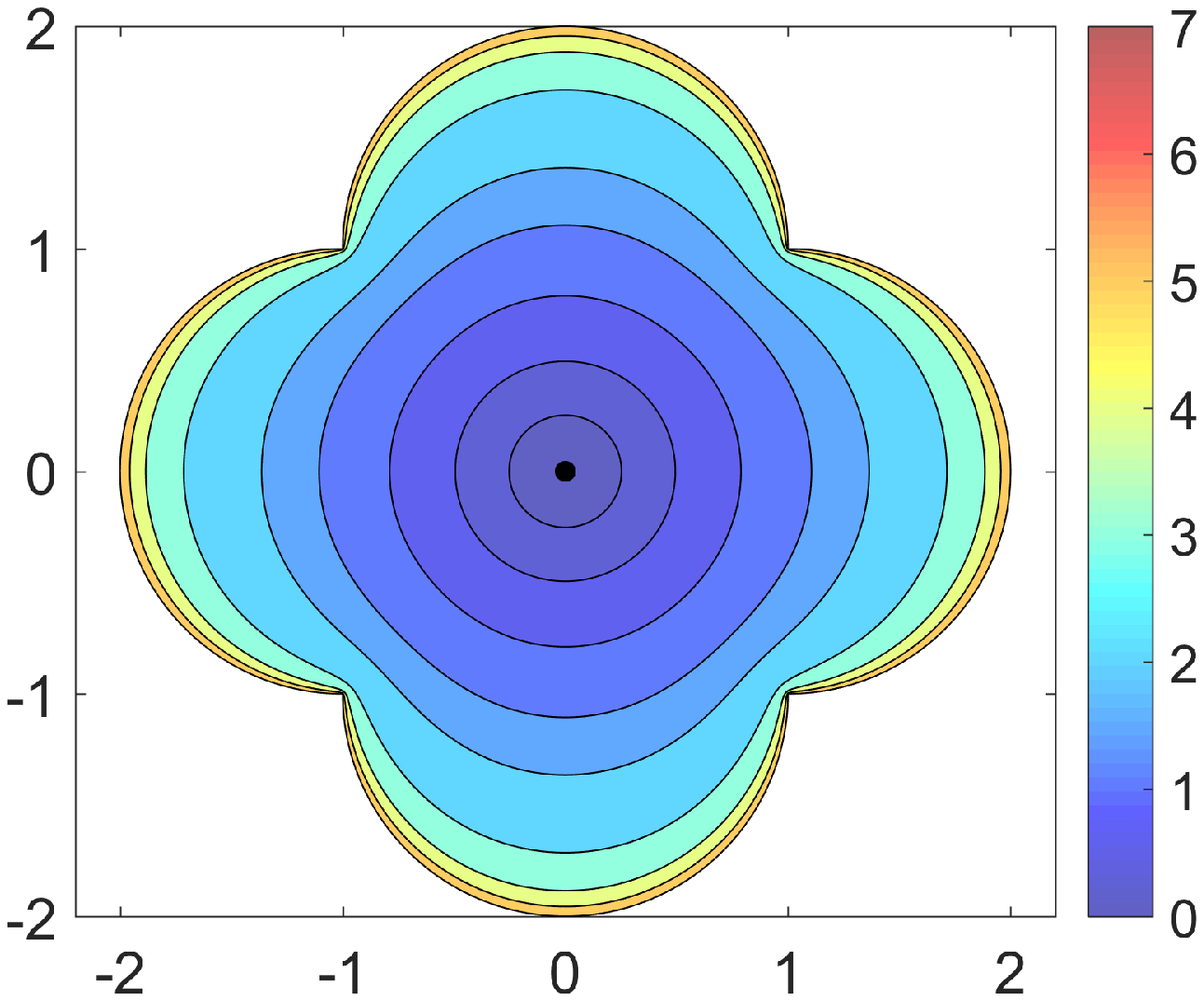}}
}
\caption{The contour lines of the function $u(x,y)$ for the circular arc quadrilateral $z_1=1.5$ (left) and $z_1=0$  (right).}
\label{fig:hd-caq}
\end{figure}

\nonsec{\bf Circular arc quadrilateral.}
For the third example, we consider the simply connected domain $G$ in the interior of the circular arc quadrilateral consists of the right-half of the circle with center $1$ and radius $1$,  the upper-half of the circle with center $\i$ and radius $1$, the left-half of the circle with center $-1$ and radius $1$, and the lower-half of the circle with center $-\i$ and radius $1$. The contour lines for the function $u(x,y)$ computed with $n=2^{12}$ are shown in Figures~\ref{fig:hd-caq} (left) for $z_1=1.5$ and in Figures~\ref{fig:hd-caq} (right) for $z_1=0$.

\begin{figure}[ht] %
\centerline{
\scalebox{0.5}[0.5]{\includegraphics[trim=0 -1.0cm 0 0,clip]{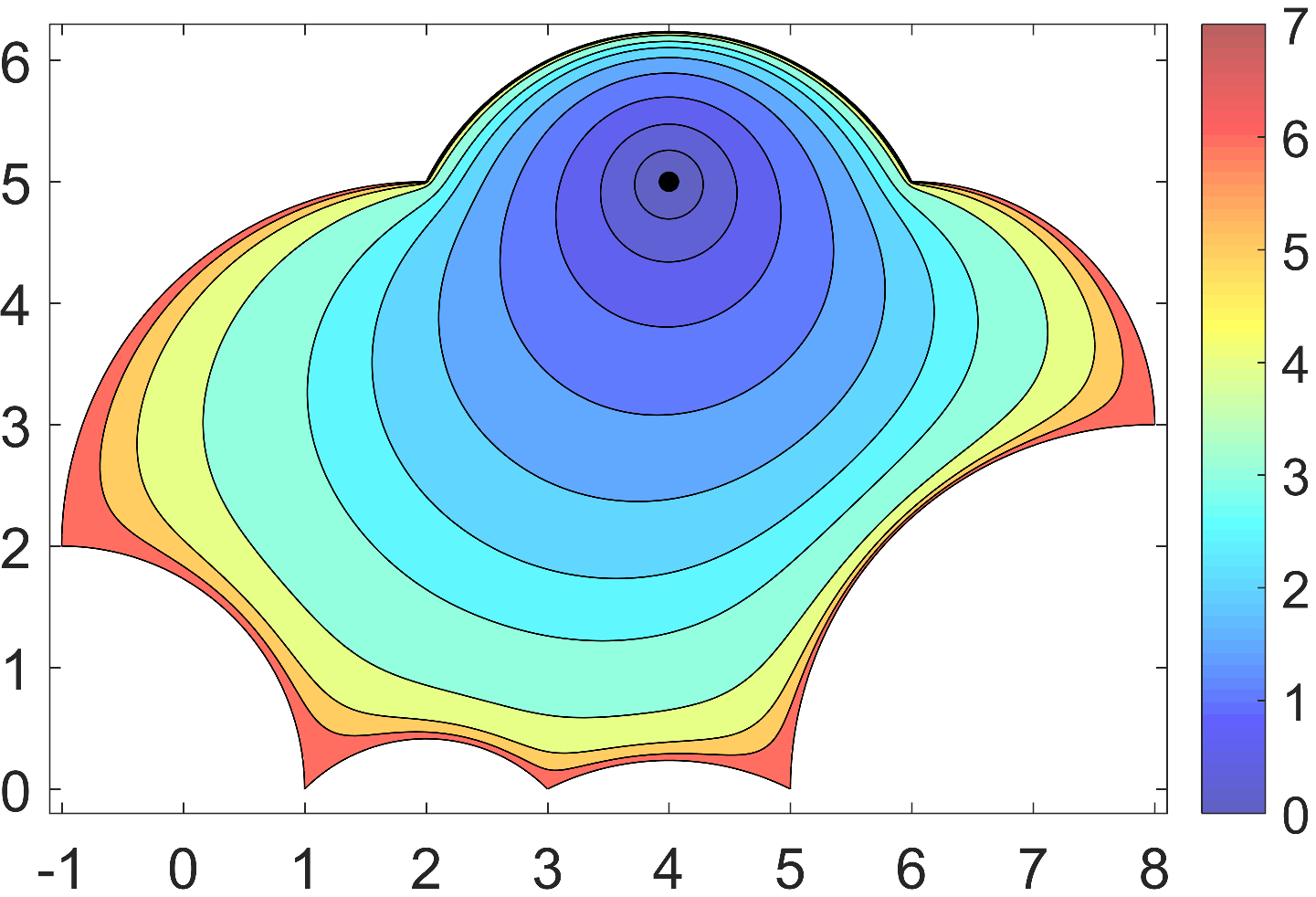}}
\hfill
\scalebox{0.5}[0.5]{\includegraphics[trim=0 -1.0cm 0 0,clip]{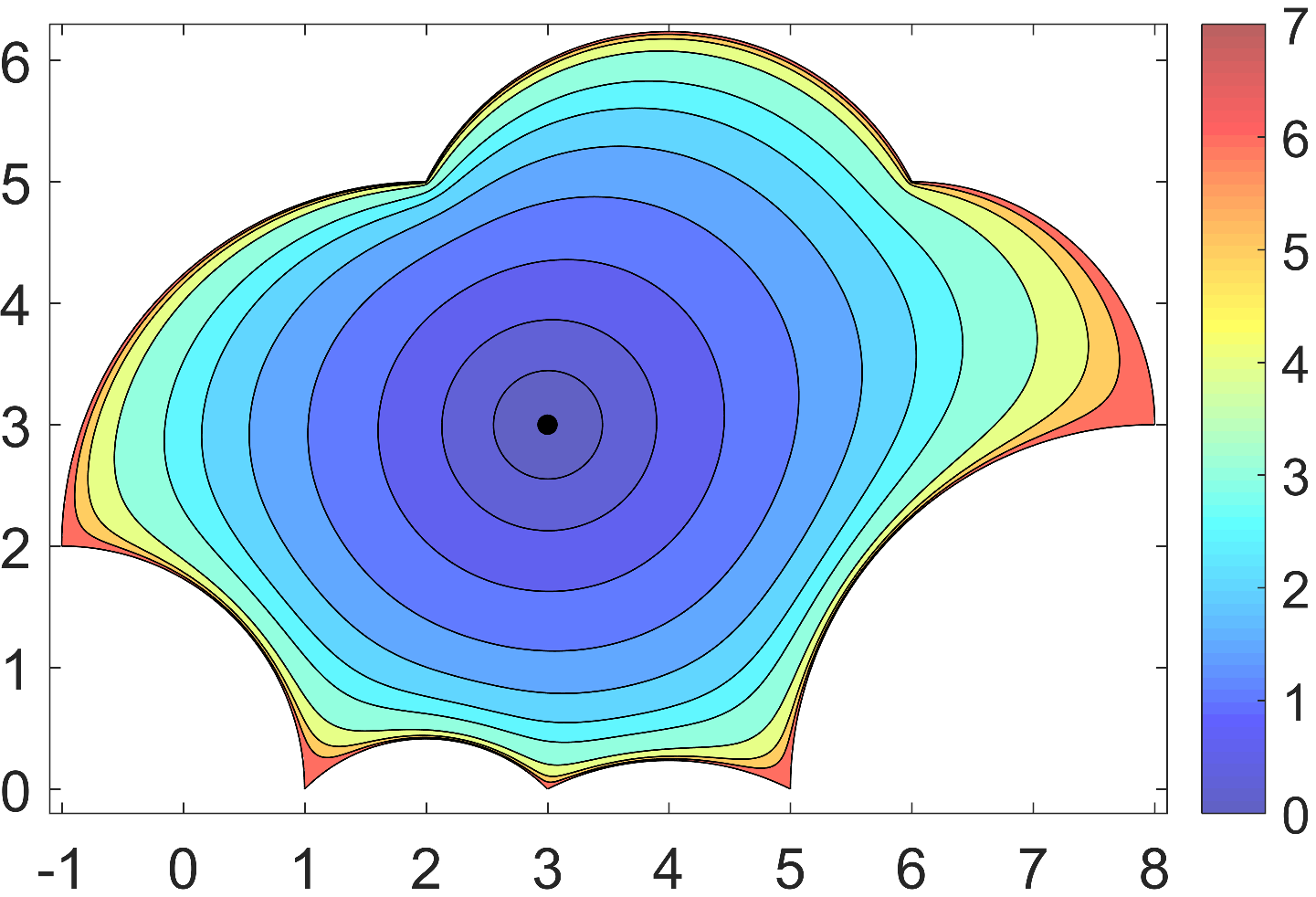}}
}
\caption{The contour lines of the function $u(x,y)$ for the circular arc polygon $z_1= 4+5\i$ (left) and  $z_1=3+3\i$  (right).}
\label{fig:hd-hyp}
\end{figure}

\nonsec{\bf Circular arc polygon.}
In the fourth example, we consider the simply connected domain $G$ in the interior of the circular arc polygon with seven sides. The contour lines for the function $u(x,y)$ computed with $n=7\times2^{10}$ are shown in Figures~\ref{fig:hd-hyp} (left) for $z_1= 4+5\i$ and in Figures~\ref{fig:hd-hyp} (right) for $z_1=3+3\i$.

%% file: sec408.tex

\section{Reduced modulus}
\label{sc:rm}

\nonsec{\bf Reduced modulus for simply connected domains.}
The \emph{reduced modulus} for simply connected domains are defined in terms of the conformal radius of simply connected domains introduced in Section~\ref{sec:cm}.

Let $G\subset\overline{\CC}$ be a bounded simply connected domain and $\alpha\in G$. The reduced modulus of the domain $G$ with respect to the point $\alpha$ is defined by~\cite[p.~16]{Vas02}, \cite[p.168, 560]{garmar},\cite[pp.26-27]{du}
\begin{equation}\label{eq:rm-b}
m(G,\alpha) = \frac{1}{2\pi}\log R(G,\alpha),
\end{equation}
where $R=R(G,\alpha)$ is the conformal radius of $G$ with respect to the point $\alpha$.
It follows from this definition that $m(G,\alpha)<0$ when $R(G,\alpha)<1$, $m(G,\alpha)=0$ when $R(G,\alpha)=1$, and $m(G,\alpha)>0$ when $R(G,\alpha)>1$.

For an unbounded simply connected domain $G\subset\overline{\CC}$ with $\infty\in G$, the reduced modulus of the domain $G$ with respect to $\infty$ is defined by~\cite[p.~17]{Vas02}
\begin{equation}\label{eq:rm-u}
m(G,\infty) = -\frac{1}{2\pi}\log R(G,\infty),
\end{equation}
where $R=R(G,\infty)$ is the conformal radius of $G$ with respect to $\infty$.

\nonsec{\bf Computing the reduced modulus of simply connected domains.}
\label{sc:cm-sim}
As was explained in Section~\ref{sec:cm}, the conformal radius of simply connected domains can be computed using the integral equation with the generalized Neumann kernel.
For bounded simply connected domains, it follows from~\eqref{eq:cr-b} that the reduced modulus of the domain $G$ with respect to the point $\alpha$ is given by
\[
m(G,\alpha) = \frac{h}{2\pi}
\]
where the constant $h$ is given by~\eqref{eq:h}. For unbounded simply connected domains, it follows from~\eqref{eq:cr-u} that the reduced modulus of the domain $G$ with respect to the point $\infty$ is given by
\[
m(G,\infty) = -\frac{h}{2\pi}.
\]

The above method for computing the conformal radius and the reduced modulus of bounded and unbounded simply connected domains can be implemented in MATLAB as in the following function.

\begin{lstlisting}
function [cr,m] = confrad (et,etp,n,alpha,type)
% This function computes the conformal radius cr=R(G,a) and the reduced
% modulus m=m(G,a) for a given simply connected domain G with respect to
% the point a=alpha for bounded G and a=inf for unbounded G where:
% et, etp:  the parametrization of the boundary of G and its derivative 
% n: the number of discretization points
% alpha: a given point in G for bounded G and alpha=beta (beta is an 
% auxiliary point in the exterior of G for unbounded G)
% type='b' for bounded G and type='u' for unbounded G
if type=='b' A = et-alpha; elseif type=='u' A = ones(size(et)); end
gam   = -log(abs(et-alpha));
[~,h] =  fbie(et,etp,A,gam,n,5,[],1e-14,200);
cr    =  exp(mean(h));
if type=='b' m = mean(h)/(2*pi); elseif type=='u' m = -mean(h)/(2*pi); end
end
\end{lstlisting}

\nonsec{\bf Domain exterior to an ellipse.}
As our first example, we consider the simply connected domain $G_r$ in the exterior of the ellipse
\[
\eta(t)=\cos t-\i r\sin t, \quad 0\le t\le 2\pi, \quad 0<r\le 1.
\]
For $r=0$, the ellipse reduces to the segment $[-1,1]$ and for $r=1$ to the unit circle.

We can easily show that the function
\[
z=\Psi(w)=w+\frac{1-r^2}{4}\frac{1}{w}
\]
maps the domain exterior to the circle $|w|=(1+r)/2$ onto the domain exterior of the ellipse. Hence, the inverse mapping
\begin{equation}\label{eq:map-to-ellipse}
w=\Phi(z)=z\left(\frac{1}{2}+\frac{1}{2}\sqrt{1-\frac{1-r^2}{z^2}}\right),
\end{equation}
maps the domain exterior to the ellipse onto the domain $|w|>(1+r)/2$, where the branch of the square root is chosen such that $\sqrt{1}=1$. It is clear that the function $\Phi$ satisfies $\Phi(\infty)=\infty$ and $\Phi'(\infty)=1$. Hence, $R(G,\infty)=(1+r)/2$ and
\[
m(G,\infty)=-\frac{1}{2\pi}\log\frac{1+r}{2}=\frac{1}{2\pi}\log\frac{2}{1+r}.
\]

\begin{figure}[ht] %
\centerline{
\scalebox{0.45}[0.45]{\includegraphics[trim=0 -1.0cm 0 0,clip]{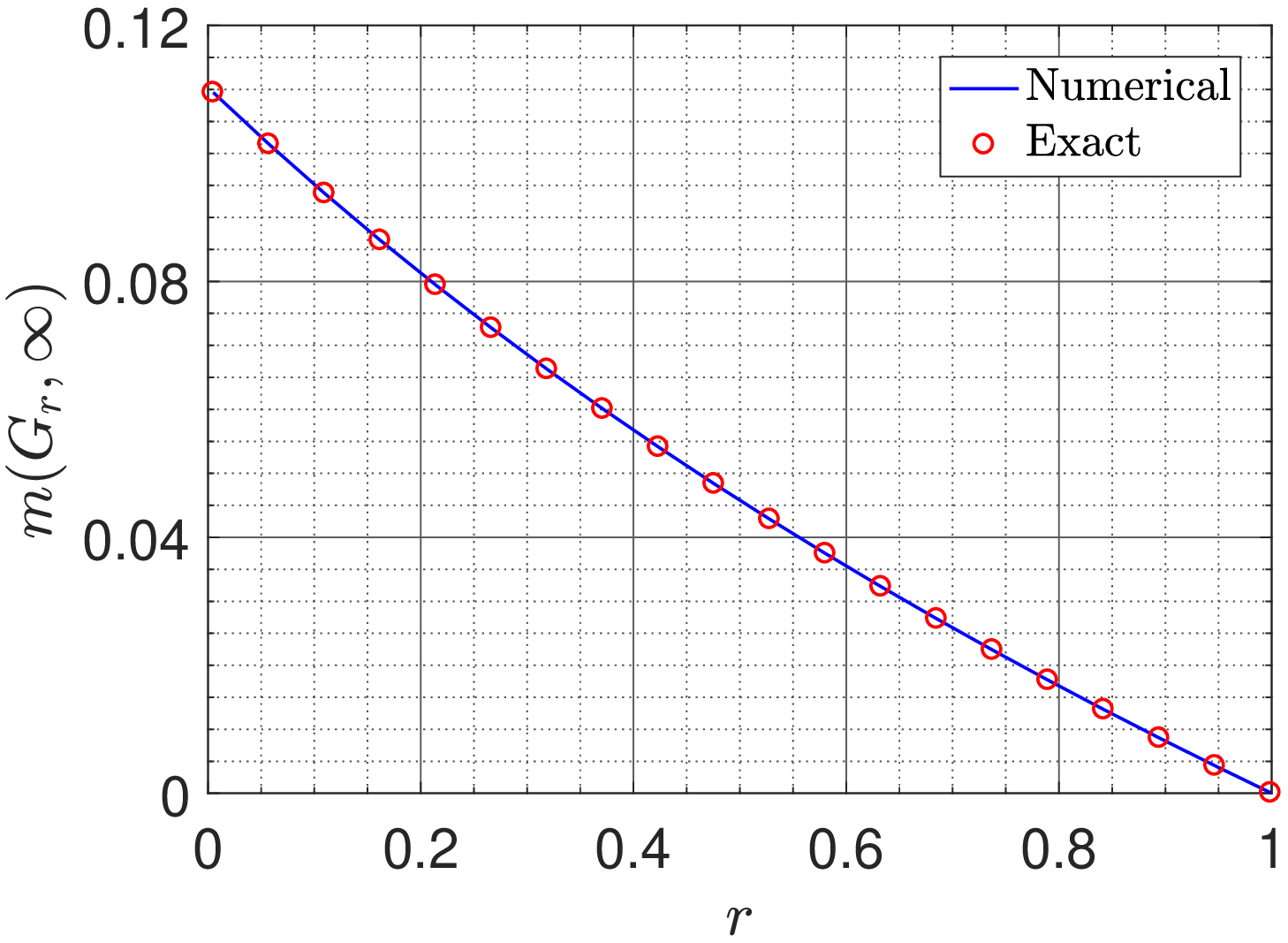}}
\hfill
\scalebox{0.45}[0.45]{\includegraphics[trim=0 -1.0cm 0 0,clip]{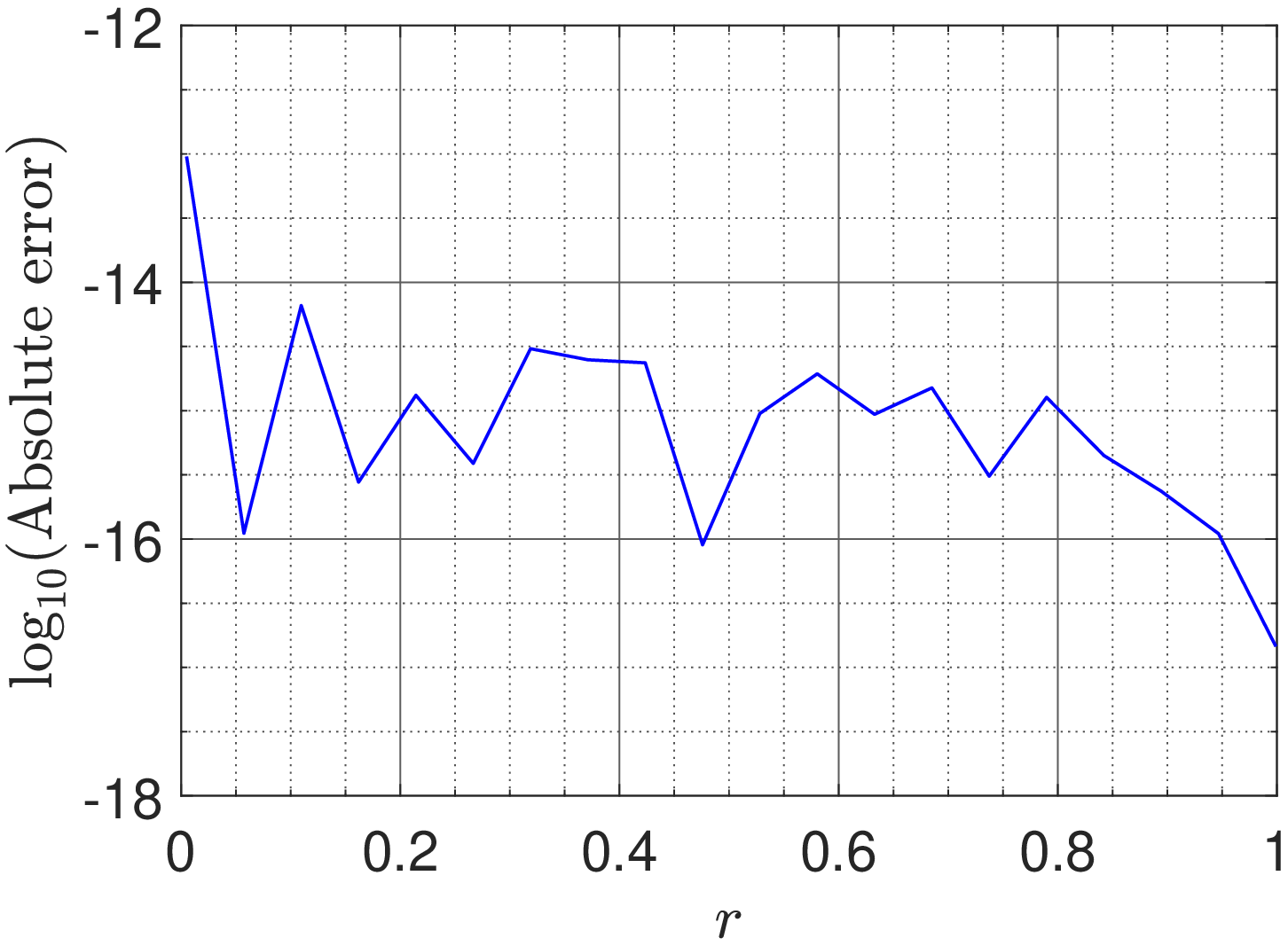}}
}
\caption{The computed and the exact reduced modulus of the domain $G_r$ exterior to an ellipse (left) and the absolute error in the computed values (right).}
\label{fig:rm-ue}
\end{figure}

We use the MATLAB function \verb|confrad| to compute the reduced modulus $m(G_r,\infty)$ with $n=2^{12}$ for $0.005\le r\le 1$. The obtained results are shown in Figure~\ref{fig:rm-ue}.

\nonsec{\bf Domain interior to an ellipse.}
For the second example, we consider the simply connected domain $G_r$ in the interior of the ellipse
\[
\eta(t)=\cosh r\cos t+\i\sinh r\sin t, \quad 0\le t\le 2\pi, \quad 0<r.
\]
Let $w=\Phi(z)$ be the unique conformal mapping from the interior of the ellipse onto the interior of the unit circle with the normalization $\Phi(0)=0$ and $\Phi'(0)>0$. The  exact form of the inverse conformal mapping $z=\Phi^{-1}(w)$ is given in~\cite{KS06}. In particular, it was shown in~\cite{KS06} that $(\Phi^{-1})'(0)=\pi/(2\sqrt{s}K(s))$ where $s=\mu^{-1}(2r)$. Hence, $\Phi'(0)=2\sqrt{s}K(s)/\pi$. Thus, the mapping function $\hat\Phi$ defined by
\[
w=\hat\Phi(z)=\frac{\Phi(z)}{\Phi'(0)}=\frac{\pi}{2\sqrt{s}K(s)}\Phi(z)
\]
is the unique conformal mapping from the interior of the ellipse onto the
\[
|w|<\frac{\pi}{2\sqrt{s}K(s)}
\]
with the normalization $\hat\Phi(0)=0$ and $\hat\Phi'(0)=1$. Thus, $R(G_r,0)=\pi/(2\sqrt{s}K(s))$ and hence
\[
m(G_r,0)=\frac{1}{2\pi}\log\frac{\pi}{2\sqrt{s}K(s)}, \quad s=\mu^{-1}(2r).
\]

\begin{figure}[ht] %
\centerline{
\scalebox{0.45}[0.45]{\includegraphics[trim=0 -1.0cm 0 0,clip]{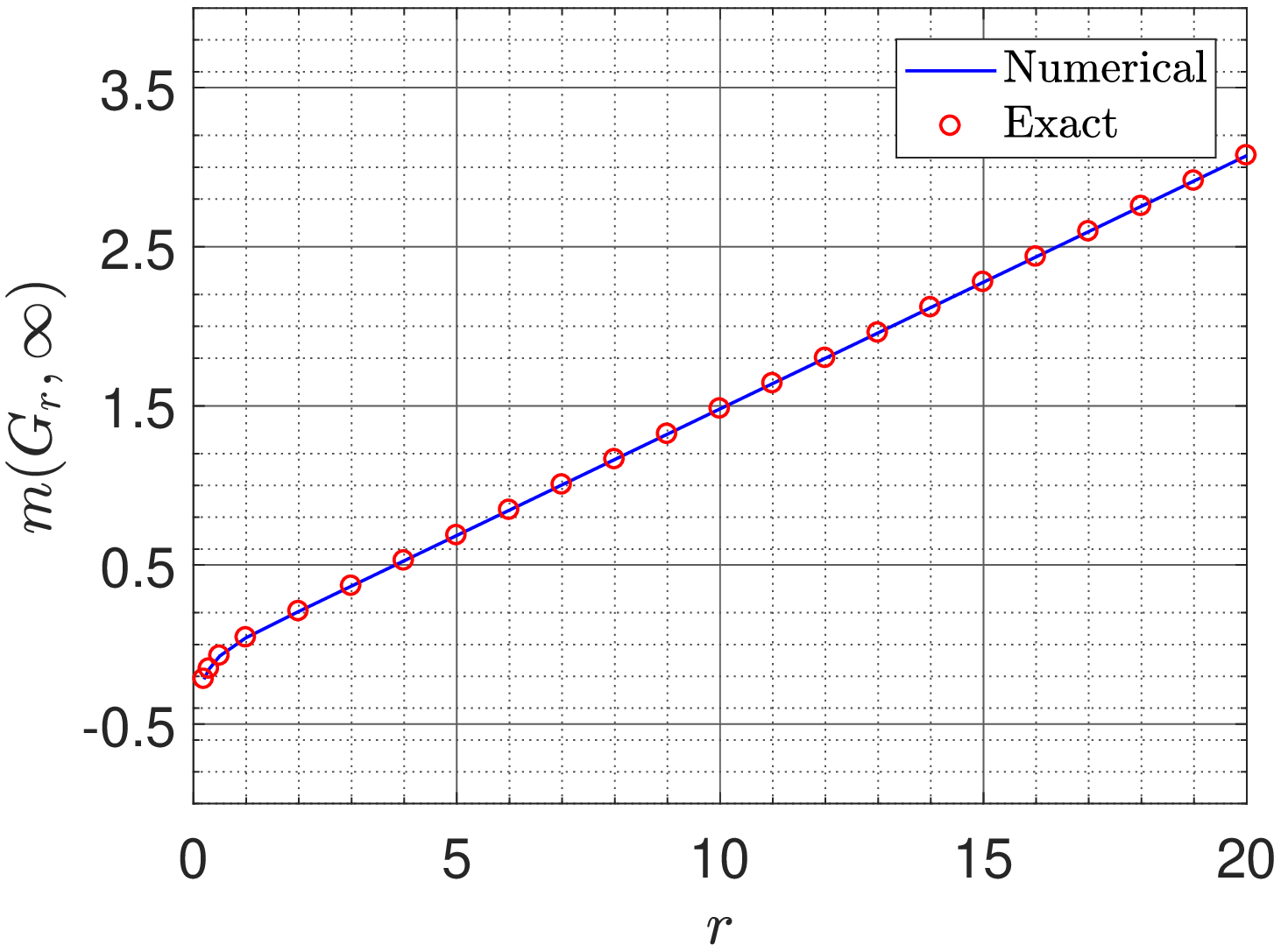}}
\hfill
\scalebox{0.45}[0.45]{\includegraphics[trim=0 -1.0cm 0 0,clip]{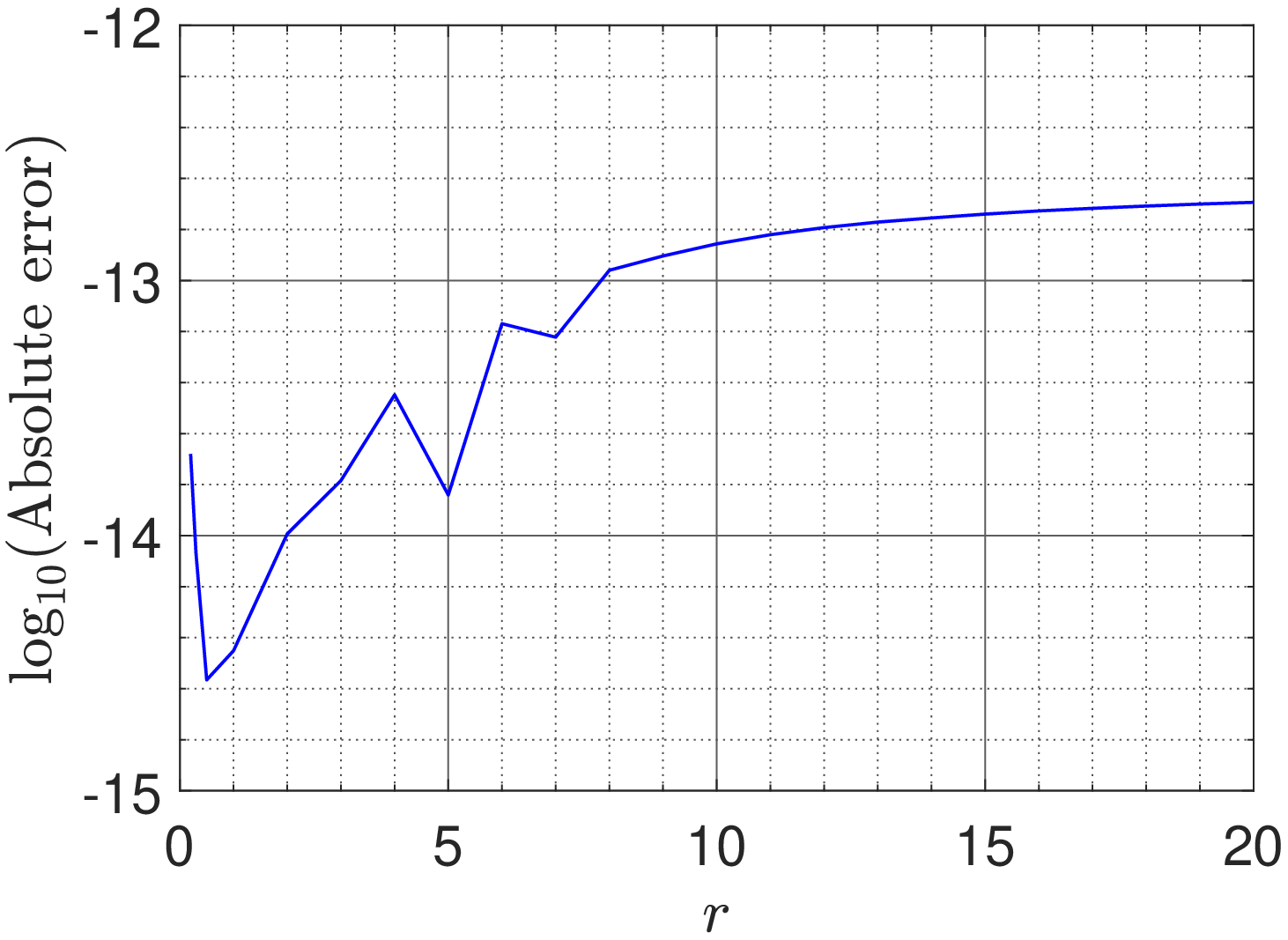}}
}
\caption{The computed and the exact reduced modulus of the domain $G_r$ interior to an ellipse  (left) and the absolute error in the computed values (right).}
\label{fig:rm-be}
\end{figure}

We use the MATLAB function \verb|confrad| to compute the reduced modulus $m(G_r,0)$ with $n=2^{12}$ for $0.2\le r\le 20$. The obtained results are shown in Figure~\ref{fig:rm-be}.

\nonsec{\bf Slitted unit disk.}
In the third example, we consider three types of slitted unit disks.
In each case the exact reduced moduli are given in~\cite[p.~33]{Vas02}.

	1) $G_1=\D\backslash(-1,0]$ where $\D$ is the unit disk. The exact value of the reduced modulus of $G_1$ with respect to $r\in(0,1)$ is given by~\cite{Vas02}
	\[
	m(G_1,r)=\frac{1}{2\pi}\log\frac{4r(1-r)}{1+r}.
	\]
	
\begin{figure}[ht] %
\centerline{
\scalebox{0.38}[0.38]{\includegraphics[trim=0 -1.0cm 0 0,clip]{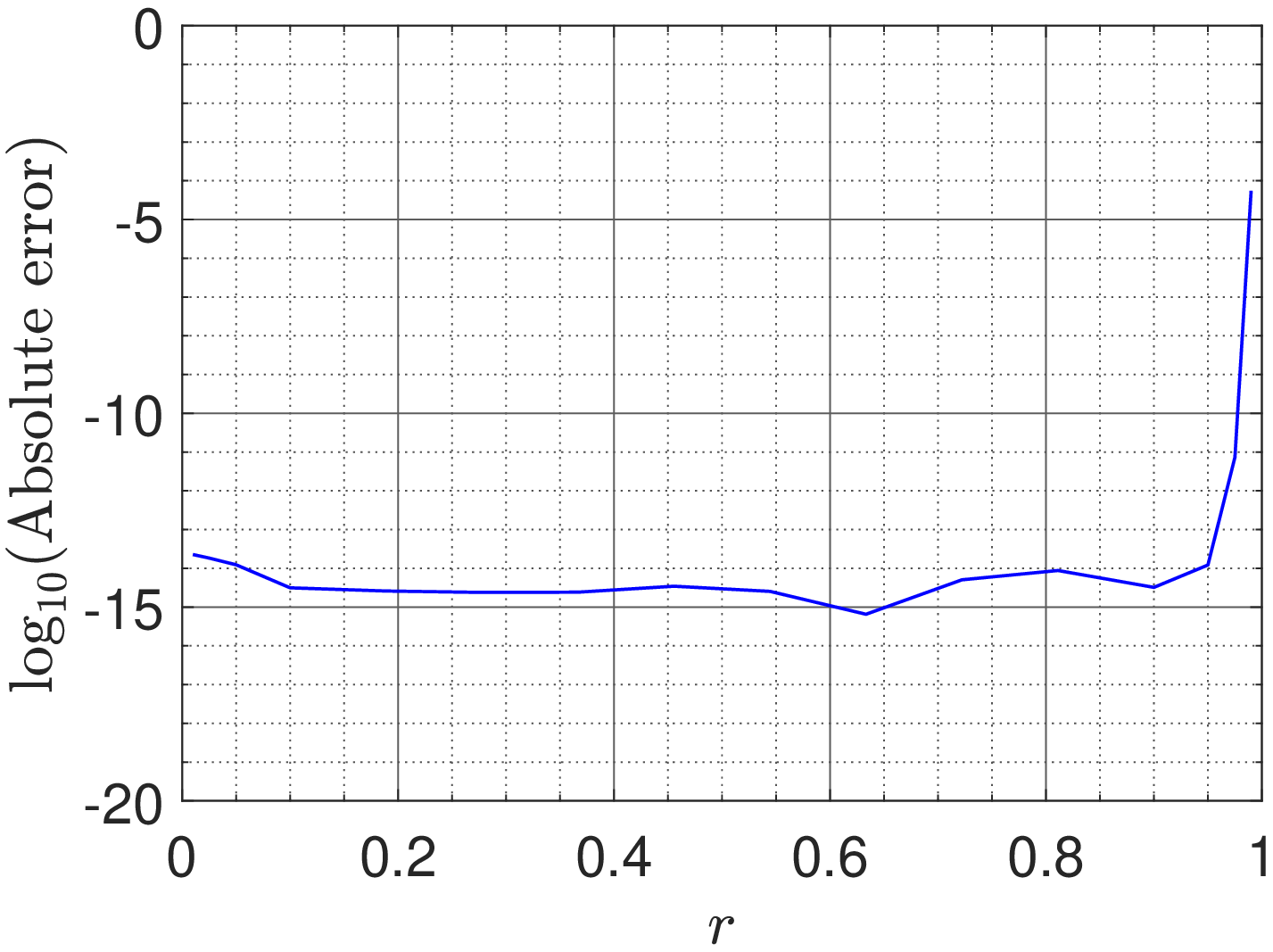}}
\hfill
\scalebox{0.38}[0.38]{\includegraphics[trim=0 -1.0cm 0 0,clip]{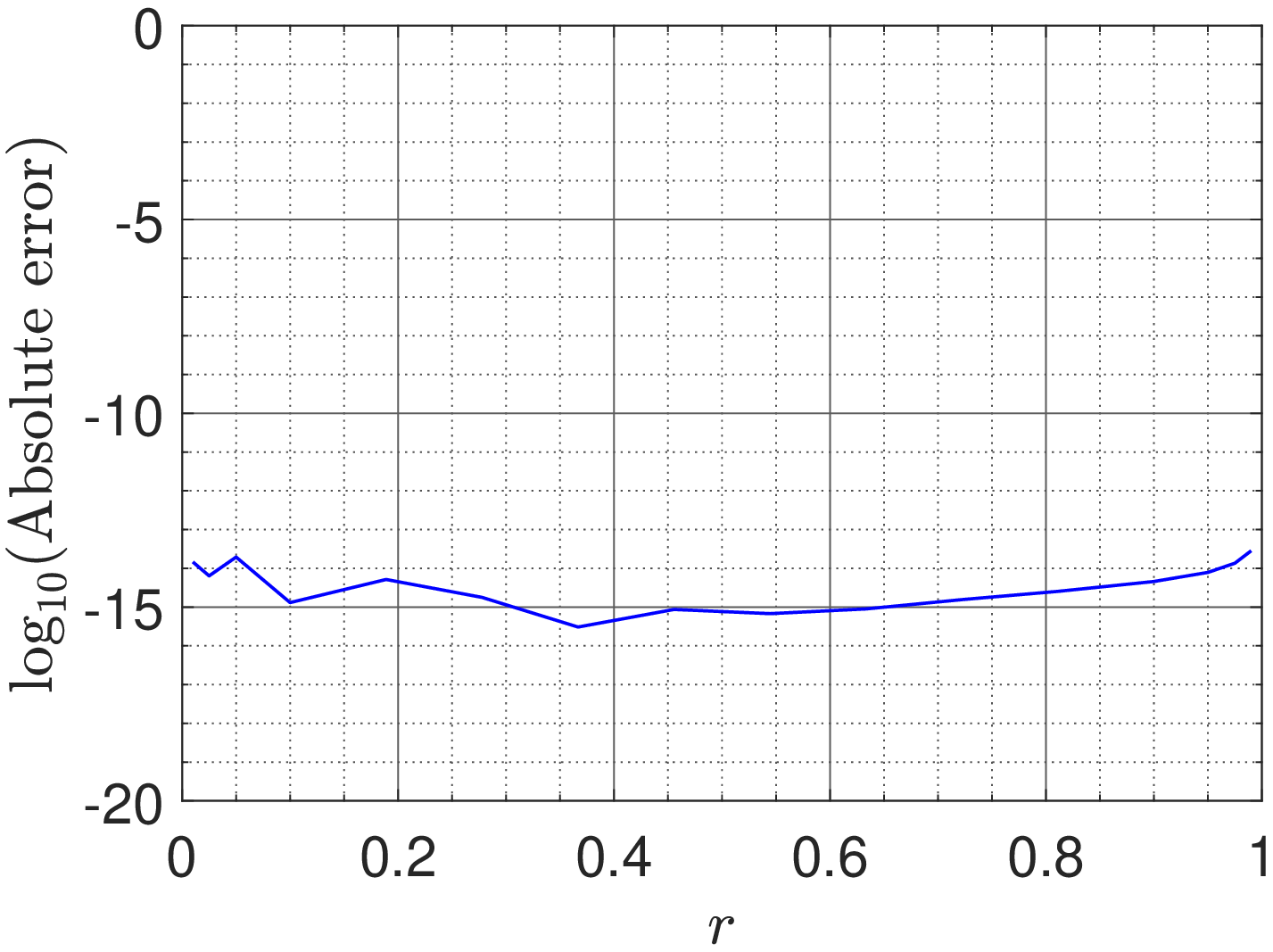}}
\hfill
\scalebox{0.38}[0.38]{\includegraphics[trim=0 -1.0cm 0 0,clip]{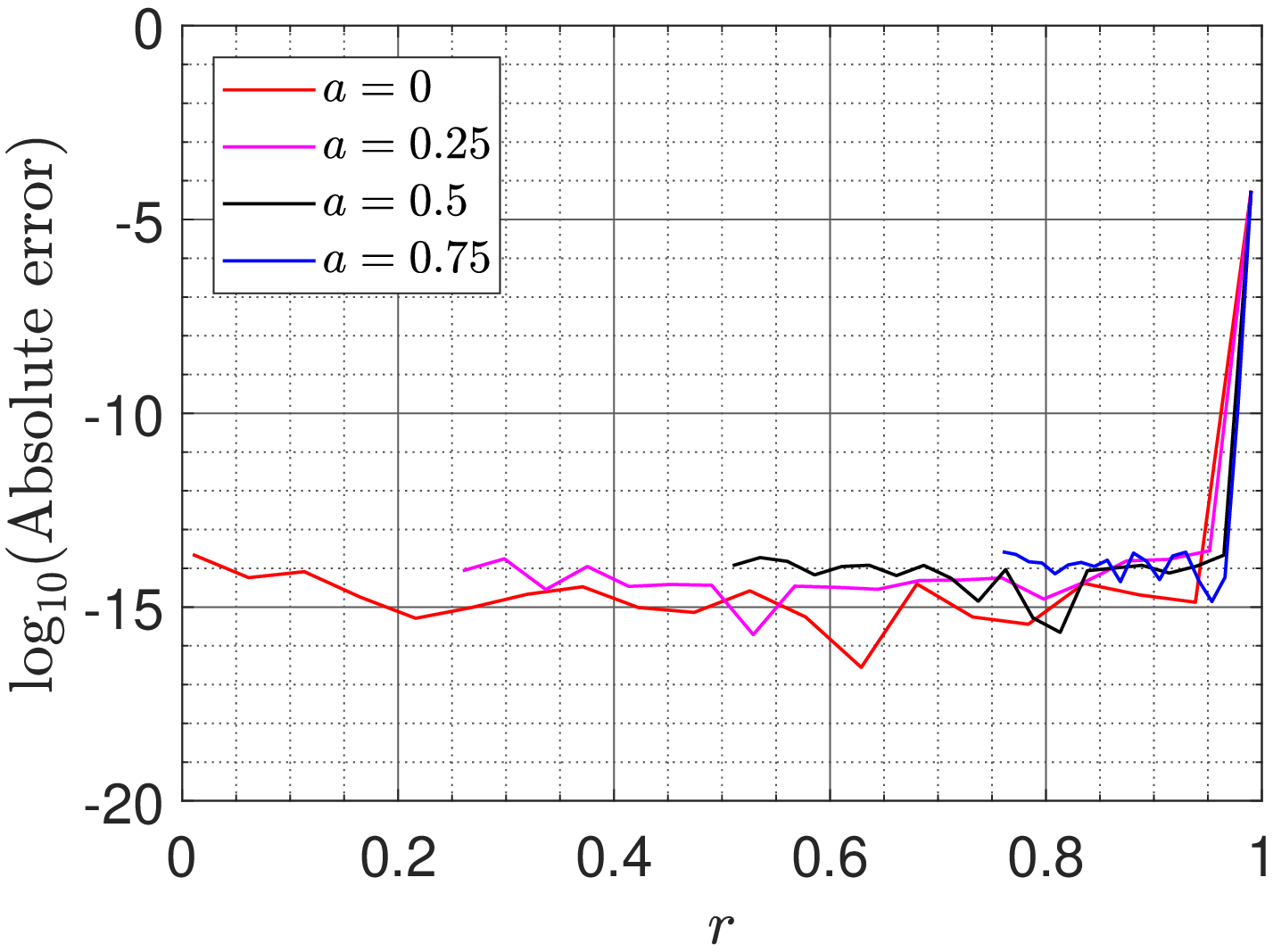}}
}
\caption{The absolute error  for $G_1$ (left),  $G_2$ (center),  and $G_3$ (right).}
\label{fig:rm-sd}
\end{figure}

	To use the integral equation to compute $m(G_1,r)$, we first use the auxiliary map
	\[
	\zeta=\Phi_1(z)=2\sqrt{r}\sqrt{z},
	\]
	where the branch of the square root is chosen on the negative real line, to open up the slit and map the region $G_1$ onto a region $\hat G_1$ bordered by piecewise smooth Jordan curve where $\Phi_1(r)=2r$ and $\Phi'_1(r)=1$.
Then, it follows from~\cite[Corollary 2.2.1]{Vas02} that $m(G_1,r)=m(\hat G_1,2r)$. We use the MATLAB function \verb|confrad| to compute the reduced modulus $m(\hat G_1,2r)$ with $n=2^{12}$ for $0.01\le r\le 0.99$. The obtained results are shown in Figure~\ref{fig:rm-sd} (left).
	
2) $G_2=\D\backslash[r,1)$ for $0<r<1$. The exact value of the reduced modulus of $G_2$ with respect to the origin is given by~\cite{Vas02}
	\[
	m(G_2,0)=\frac{1}{2\pi}\log\frac{4r}{(1+r)^2}.
	\]
	
	To compute $m(G_2,r)$, we first use the auxiliary map
	\[
	\zeta=\Phi_1(z)=2\i\sqrt{r}\sqrt{z-r},
	\]
	where the branch of the square root is chosen on the positive real line, to map the region $G_2$ onto a region $\hat G_2$ bordered by piecewise smooth Jordan curve where $\Phi_1(0)=-2r$ and $\Phi'_1(0)=1$. Then, $m(G_2,0)=m(\hat G_2,-2r)$. We use the MATLAB function \verb|confrad| to compute $m(\hat G_2,-2r)$ with $n=2^{12}$ for $0.01\le r\le 0.99$. The obtained results are shown in Figure~\ref{fig:rm-sd} (center).

3) $G_3=\D\backslash(-1,a]$ for $0\le a<1$. The exact value of the reduced modulus of $G_3$ with respect to $r\in(a,1)$ is given by~\cite{Vas02}
	\[
	m(G_3,r)=\frac{1}{2\pi}\log\frac{4(r-a)(1-ra)(1-r)}{(1+r)(1-a)^2}.
	\]
	
	To compute $m(G_3,r)$, we first use the auxiliary map
	\[
	\zeta=\Phi_1(z)=2\sqrt{r-a}\sqrt{z-a},
	\]
	where the branch of the square root is chosen on the negative real line, to map the region $G_3$ onto a region $\hat G_3$ bordered by a piecewise smooth Jordan curve where $\Phi_1(r)=2(r-a)$ and $\Phi'_1(r)=1$.
Hence, $m(G_3,r)=m(\hat G_3,2(r-a))$. We use the MATLAB function \verb|confrad| to compute $m(\hat G_3,2(r-a))$ with $n=2^{12}$ for $a=0,0.25,0.5,0.75$ and $a+0.01\le r\le 0.99$. The obtained results are shown in Figure~\ref{fig:rm-sd} (right).

\begin{figure}[ht] %
\centerline{
\scalebox{0.45}[0.45]{\includegraphics[trim=0 -1.0cm 0 0,clip]{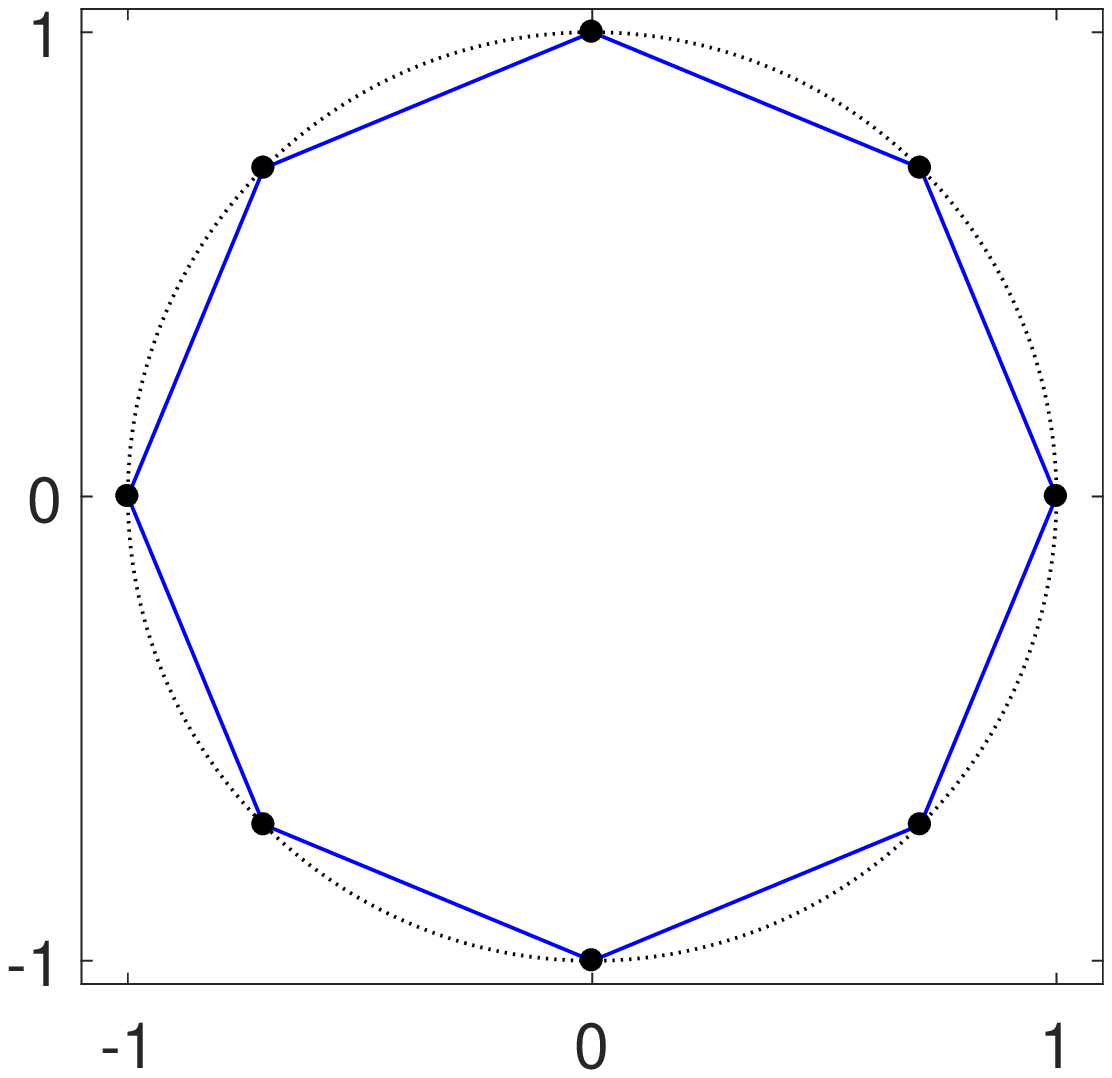}}
\hfill
\scalebox{0.45}[0.45]{\includegraphics[trim=0 -1.0cm 0 0,clip]{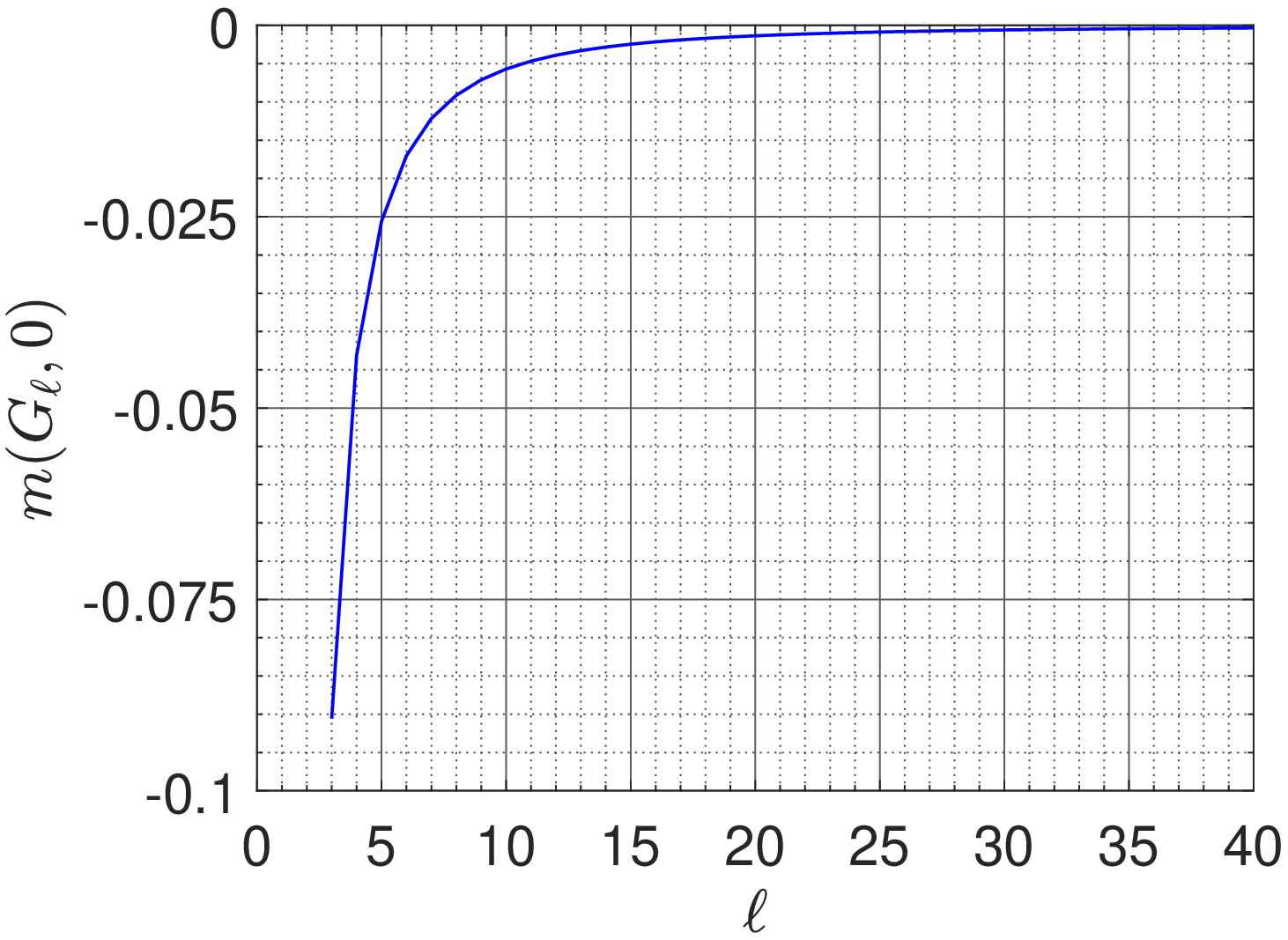}}
}
\caption{The domain $G_\ell$ for $\ell=8$ (left) and the computed reduced modulus for the domain $G_\ell$ for $3\le\ell\le40$ (right).}
\label{fig:rm-pg}
\end{figure}

\nonsec{\bf Polygon.}
For the fourth example, we consider the simply connected domain $G_\ell$ in the interior of the a polygon with $\ell$ vertices where $\ell\ge3$ (see Figure~\ref{fig:rm-pg} (left) for $\ell=8$). We assume that the vertices of the polygon are given by
\[
v_k=e^{\frac{2k\pi\i}{\ell}}, \quad k=0,1,2,\ldots,\ell-1.
\]

In this example, the exact value of the reduced modulus is unknown. We use the MATLAB function \verb|confrad| to compute the reduced modulus $m(G_\ell,0)$ with $n=\ell\times2^{9}$ for $\ell=3,4,\ldots,40$. The obtained results are shown in Figure~\ref{fig:rm-pg}. It is clear from this figure that $m(G_\ell,0)<0$ which means that $R(G_\ell,0)<1$ for the above values of $\ell$. In other words, the conformal mapping $\Phi$ with the normalization~\eqref{eq:cm-cond-b2-1} maps the domain $G_\ell$ onto a disk interior to the unit disk.

%% file: sec508.tex

\section{Harmonic measure}


\nonsec{\bf Harmonic measure.} 
Let $G$ be a Jordan domain in $\overline{\CC}$ and $\Gamma$ be its boundary. Let also $L$ be a boundary arc on $\Gamma$ such that $L\ne\emptyset$ and $\Gamma\backslash L\ne\emptyset$. 
The harmonic measure of $L$ with respect o $G$ is the $C^2(G)$ function $u: G \to (0,1)$ satisfying the Laplace equation
\[
\Delta u = 0
\]
in $G$ and $u(z) \to 1$ when $z \to L $ and $u(z) \to 0$ when $z\to\Gamma\setminus L\,$.
The harmonic measure is one of the key notions of potential theory and it has numerous applications to geometric function theory \cite{garmar}.
The harmonic measure of $L$ with respect to $G$ will be denoted by $\omega(z,L)$ (see e.g.,~\cite[p.~123]{avv}, \cite[Ch I]{garmar}, and~\cite[p.~111]{Tsu59}).

\nonsec{\bf Harmonic measure for the unit disk.}
Assume that $G$ is the unit disk $|z|<1$, $\Gamma$ is the unit circle $|z|=1$, and $L$ is the right half of the unit circle. It is clear that the M\"obius transformation
\[
z \mapsto \frac{z-\i}{\i z-1}
\]
maps the unit circle onto the real line and the interior of the unit circle onto the upper half plane. More precisely, it maps the right half of the unit circle onto the negative real line, maps the point $\i$ onto $0$, maps the left half of the unit circle onto the positive real line, and maps $-\i$ onto $\infty$. Hence, the harmonic measure of $L$ with respect to $G$ is given by~\cite[p.~123]{avv}, \cite[Ch I]{garmar},
\begin{equation}\label{eq:hm-disk}
\omega(z,L)=\frac{1}{\pi}\Im\log\frac{z-\i}{\i z-1},
\end{equation}
where the branch with $\log 1=0$ is chosen.

\begin{figure}[ht] %
\centerline{
\scalebox{0.45}[0.45]{\includegraphics[trim=0cm 0cm 0cm 0cm,clip]{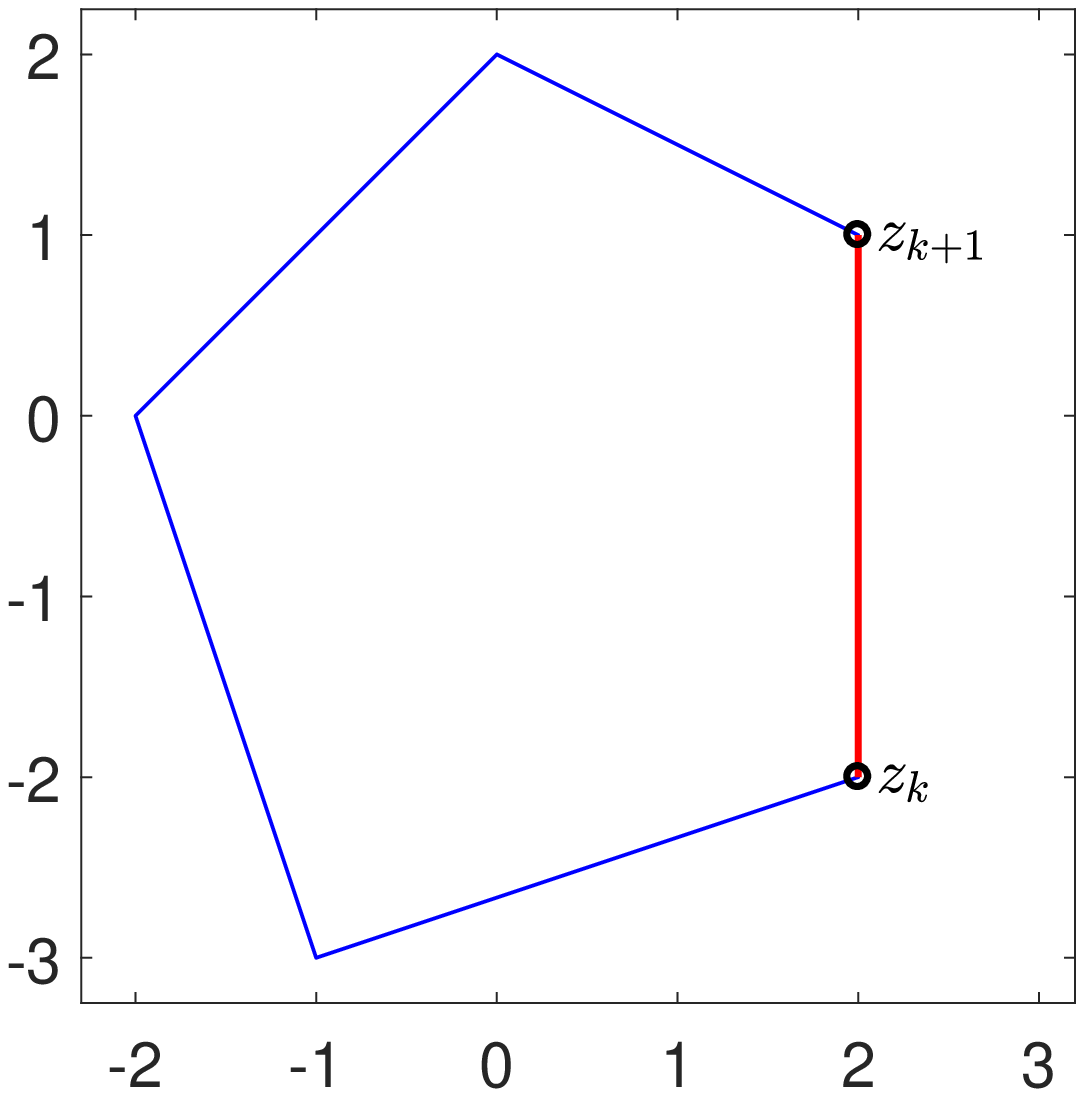}}
\hfill
\scalebox{0.45}[0.45]{\includegraphics[trim=0cm 0cm 0cm 0cm,clip]{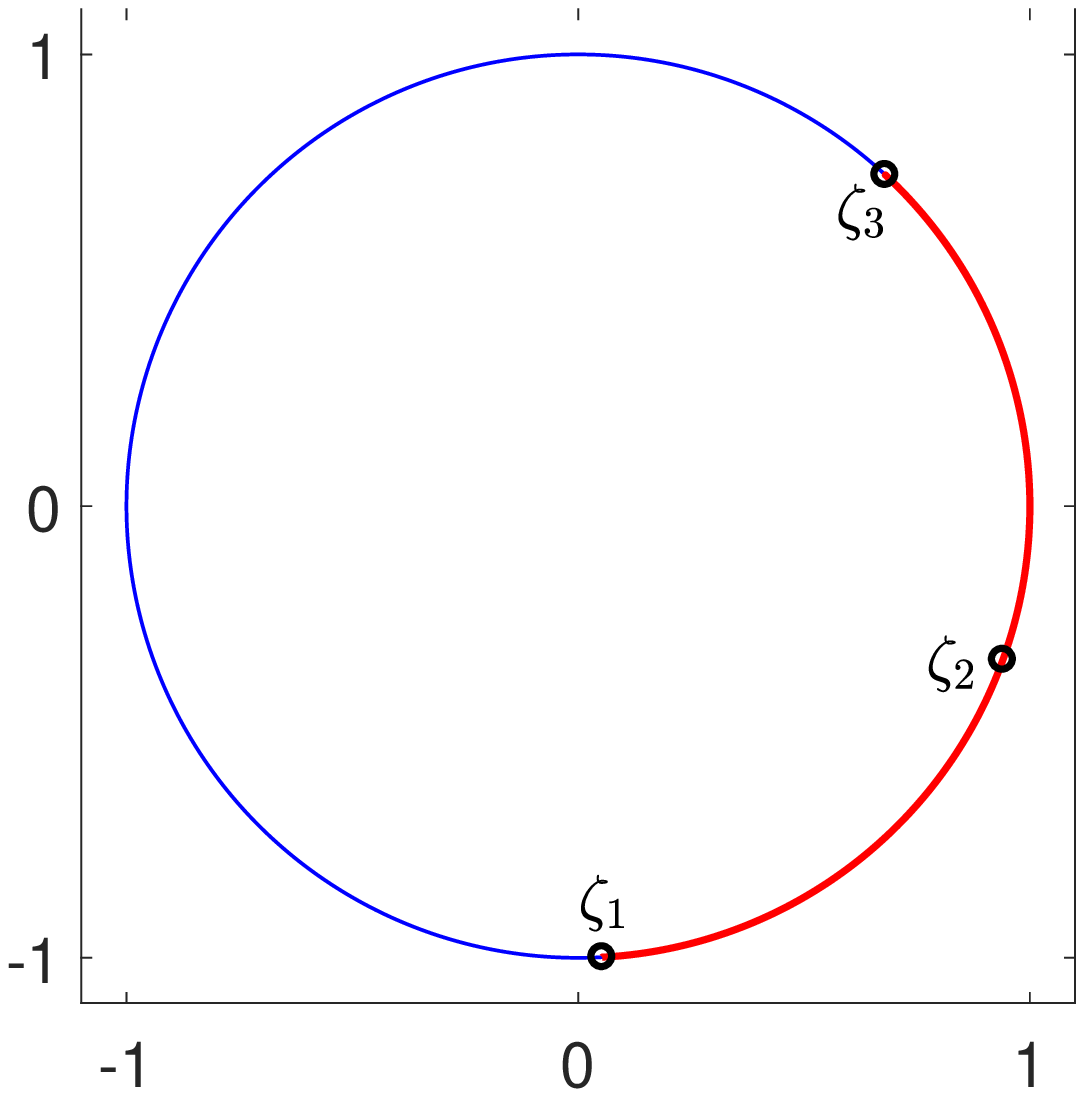}}
\hfill
\scalebox{0.45}[0.45]{\includegraphics[trim=0cm 0cm 0cm 0cm,clip]{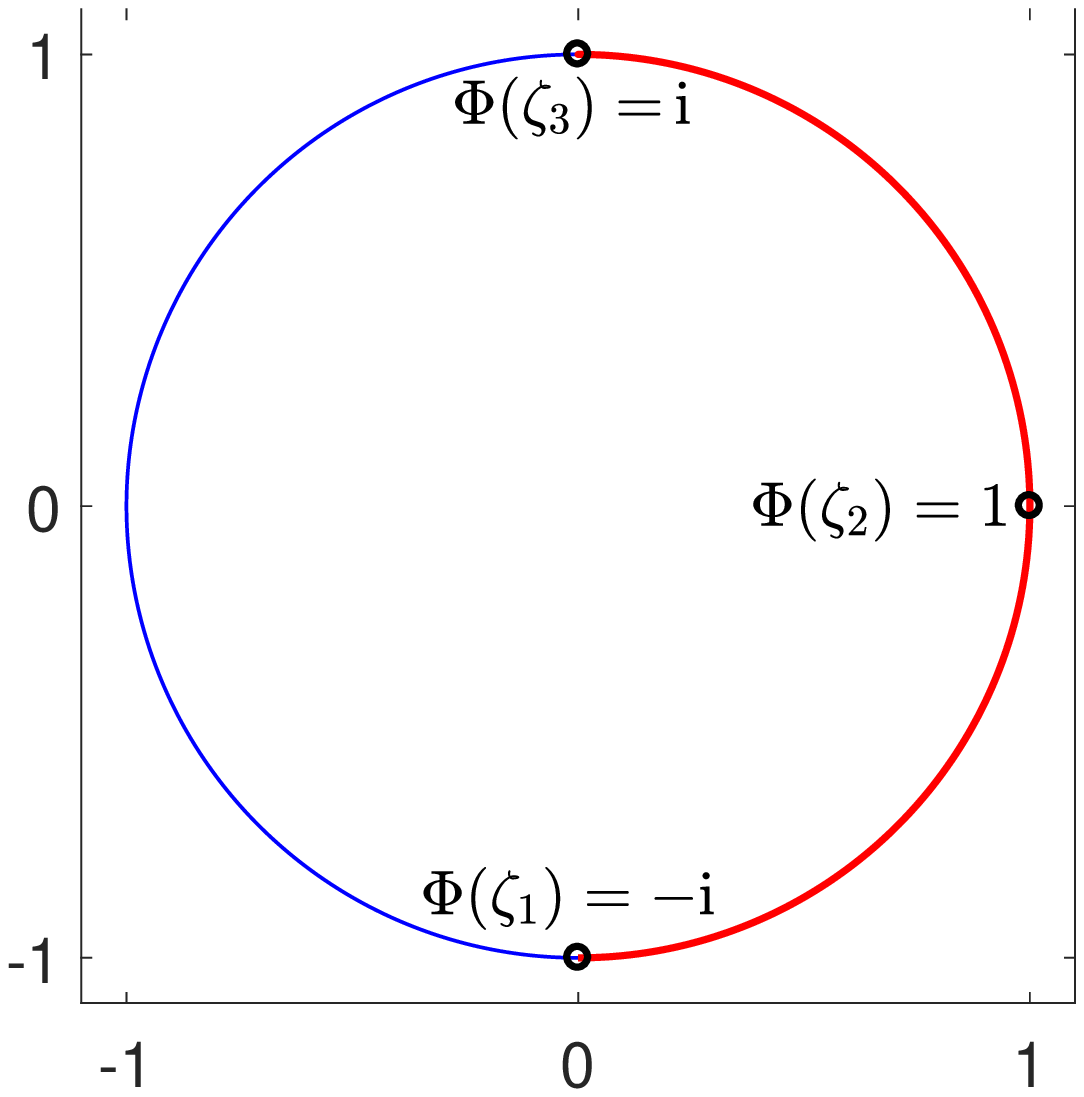}}
}
\caption{The arc $L$ between $z_1$ and $z_3$ and}
\label{fig:hm-cir}
\end{figure}

\nonsec{\bf Harmonic measure for a polygon.}
Assume that $G$ is the interior domain of a polygon $\Gamma$ with $m$ vertices $\{z_1,z_2,\ldots,z_m\}$, labelled in counterclockwise orientation, and $L$ is the segment $[z_k,z_{k+1}]$ for $k=1,2,\ldots,m$ (we define $z_{m+1}=z_1$) (see Figure~\ref{fig:hm-cir} (left)).

To compute $\omega(z,L)$, we discretize the parametrization $\eta(t)$, $0\le t\le 2\pi$, of the polygon $\Gamma$ on each segment $[z_k,z_{k+1}]$ by $n_s$ graded points on $[2(k-1)\pi/m,2k\pi/m]$. Thus, the whole polygon $\Gamma$ is discretized by $n=mn_s$ point $t_i$, $i=1,2,\ldots,n$ in $[0,2\pi]$ such that $z_k=\eta(t_{1+(k-1)n_s})$ for $k=1,2,\ldots,m$.
Then we use the method presented in Section~\ref{sc:hyp} to compute the conformal mapping $\zeta=\Phi(z)$ from the interior of $\Gamma$ onto the unit disk $|\zeta|<1$. The mapping function $\Phi$ maps the two points $z_k$ and $z_{k+1}$ onto two points $\zeta_1$ and $\zeta_3$, respectively, on the unit circle $|\zeta|=1$. The segment $L$ is then mapped onto the arc $\hat L$ on the unit circle $|\zeta|=1$ from $\zeta_1$ to $\zeta_3$. Let $\zeta_2$ be the point on the middle of $\hat L$ between $\zeta_1$ and $\zeta_3$ so that $\zeta_1$, $\zeta_2$ and $\zeta_3$ arranged in counterclockwise orientation (see Figure~\ref{fig:hm-cir} (center)). Then the M\"obius transformation
\[
w=\Psi(\zeta)=\frac{(\zeta-\zeta_1)(\zeta_2-\zeta_3)-\i(\zeta-\zeta_3)(\zeta_2-\zeta_1)}{(\zeta-\zeta_3)(\zeta_2-\zeta_1)-\i(\zeta-\zeta_1)(\zeta_2-\zeta_3)}
\]
maps the unit disk $|\zeta|<1$ onto the unit disk $|w|<1$ and maps the unit circle $|\zeta|=1$ onto the unit circle $|w|=1$ such that the points $\zeta_1$, $\zeta_2$ and $\zeta_3$ are mapped onto the points $-\i$, $1$ and $\i$, respectively. Thus, the mapping function $\Psi$ maps the arc $\hat L$ on $|\zeta|=1$ onto the right half of the unit circle $|w|=1$ (see Figure~\ref{fig:hm-cir} (right)).

Finally, the mapping function
\[
w=\Psi(\Phi(z))
\]
maps the domain $G$ onto the disk $|w|<1$ and maps the segment $L$ on $\Gamma$ onto the right half of the unit circle $|w|=1$. Hence, by~\eqref{eq:hm-disk}, the harmonic measure of $L$ with respect to $G$ is given by
\[
\omega(z,L)=\frac{1}{\pi}\Im\log\frac{\Psi(\Phi(z))-\i}{\i\Psi(\Phi(z))-1}.
\]

The above method for computing the harmonic measure of a segment $L=[z_k,z_{k+1}]$ with respect to the polygon domain $G$ can be implemented in MATLAB as in the following function \verb|hm.m| where the discretization of the parametrization of the polygon is computed using the MATLAB function \verb|polygonp.m| (see Appendix *.*).

\begin{lstlisting}
function   Uz=hm(ver,L,alpha,z,ns)
% Compute the harmonic measure w(z,L) of a side L with respect to the
% polygon domain G with the vertices ver
% Input:
% ver=[z1,z2,...,zm] (the vertices of the polygon)
% L=[z_k,z_(k+1)] (a side of the polygon)
% alpha (a given point in the domain G)
% z (a vector of points z in G)
% ns (the graded points on each side of the polygon)
% Output:
% Uz (the values of the harmonic measure w(z,L) at the points z).
% compute the the parametrization of the polygon
%
[et,etp]=polygonp(ver,ns);
% Compute the values of the mapping function \Phi(z)
A     =  et-alpha;
gam   = -log(abs(et-alpha));
[mu,h]=  fbie(et,etp,A,gam,length(ver)*ns,5,[],0.5e-14,100);
fet   = (gam+h+i*mu)./A;   c = exp(-mean(h));
f_z   =  fcau(et,etp,fet,z(:).');
Phi_z = c.*(z-alpha).*exp((z-alpha).*f_z);
% Compute the arc \hat L
zet   = c.*(et-alpha).*exp((et-alpha).*fet);
for k=1:length(ver), iver(L(1)==ver(k)) = k; end
ver1 = et((iver-1)*ns+1); cver1 = zet((iver-1)*ns+1);
if iver==length(ver) ver3 = et(1); else ver3 = et(iver*ns+1); end
if iver==length(ver) cver3 = zet(1); else cver3 = zet(iver*ns+1); end
ang1 = angle(cver1); ang3 = angle(cver3); ang3(ang3<ang1)=ang3+2*pi;
ang2 = (ang1+ang3)/2; cver2 = exp(i*ang2); hL = [cver1;cver2;cver3];
% Compute the values of the Mobius transform \Psi(\Phi(z))
Psi = @(z,v)(((z-v(1)).*(v(2)-v(3))-i*(z-v(3)).*(v(2)-v(1)))./...
             ((z-v(3)).*(v(2)-v(1))-i*(z-v(1)).*(v(2)-v(3))));
Psi_Phi_z = Psi(Phi_z,hL);
% compute Uz=w(z,L)
H = @(z)((1/pi)*imag(log((i-z)./(1-i.*z))));
Uz = H(Psi_Phi_z);
end
\end{lstlisting}

\begin{figure}[ht] %
\centerline{
\scalebox{0.5}{\includegraphics[trim=1.75cm 0.0cm 1.75cm 0.0cm,clip]{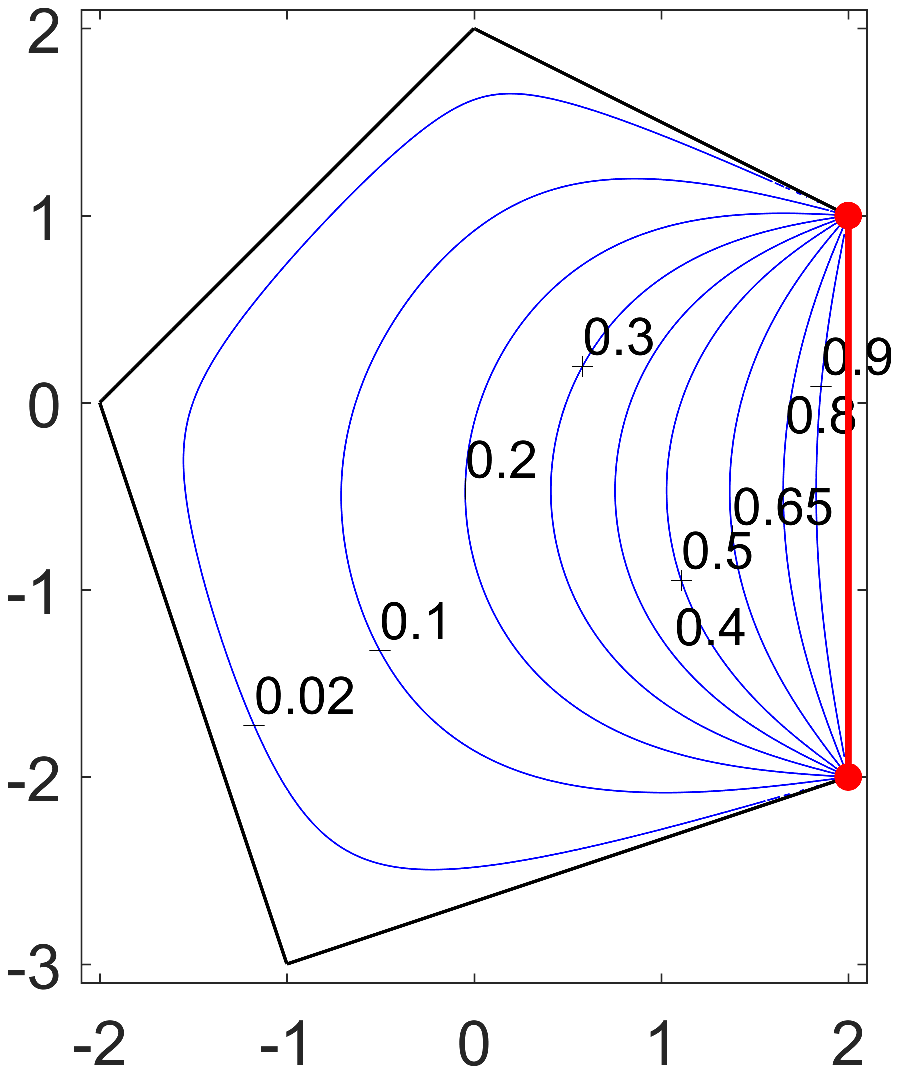}}
\hfill
\scalebox{0.5}{\includegraphics[trim=1.75cm 0.0cm 1.75cm 0.0cm,clip]{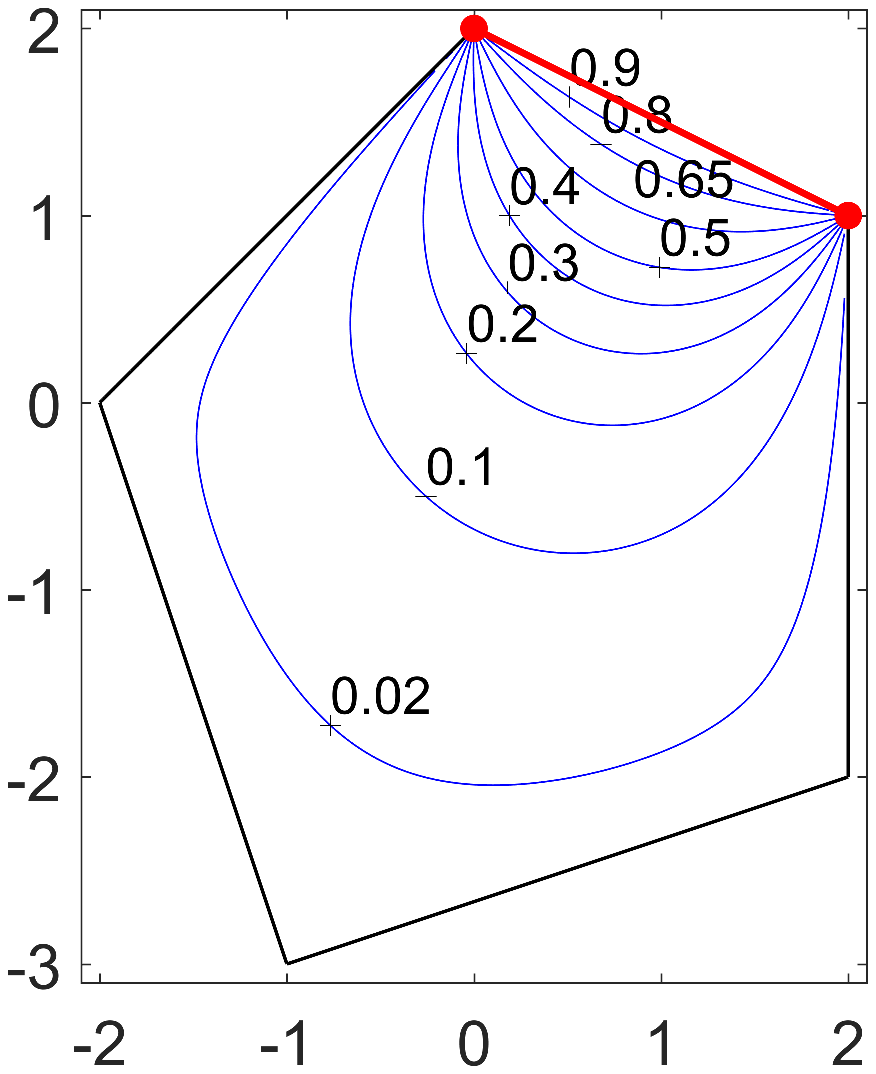}}
\hfill
\scalebox{0.5}{\includegraphics[trim=1.75cm 0.0cm 1.75cm 0.0cm,clip]{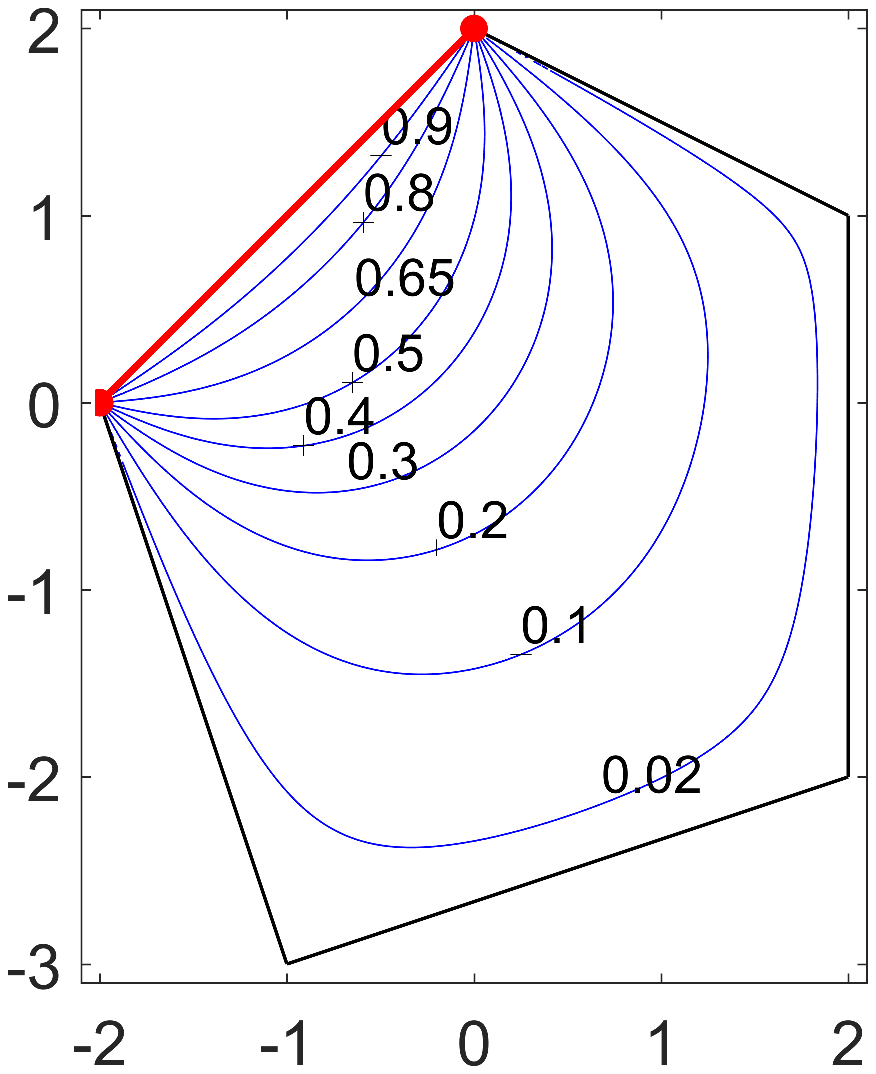}}
}
\centerline{
\hfill
\scalebox{0.5}{\includegraphics[trim=1.5cm 0.0cm 1.5cm 0.0cm,clip]{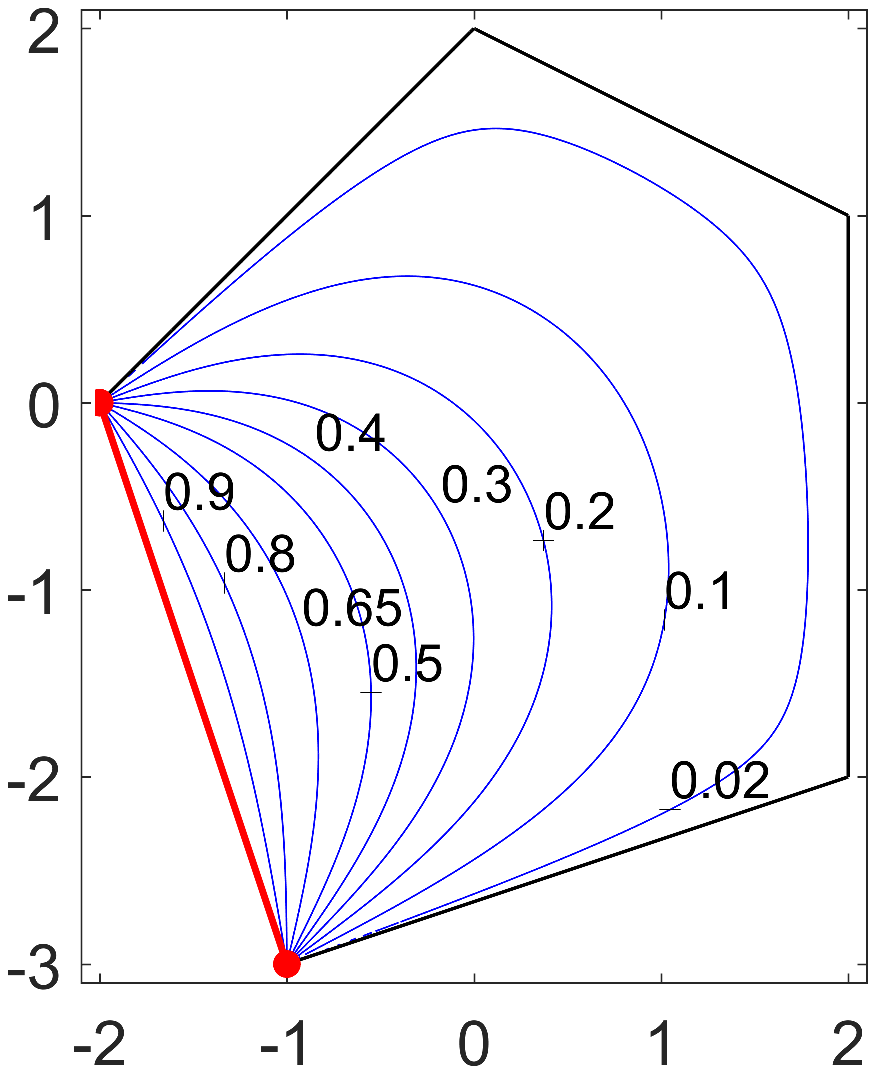}}
\hfill
\scalebox{0.5}{\includegraphics[trim=1.5cm 0.0cm 1.5cm 0.0cm,clip]{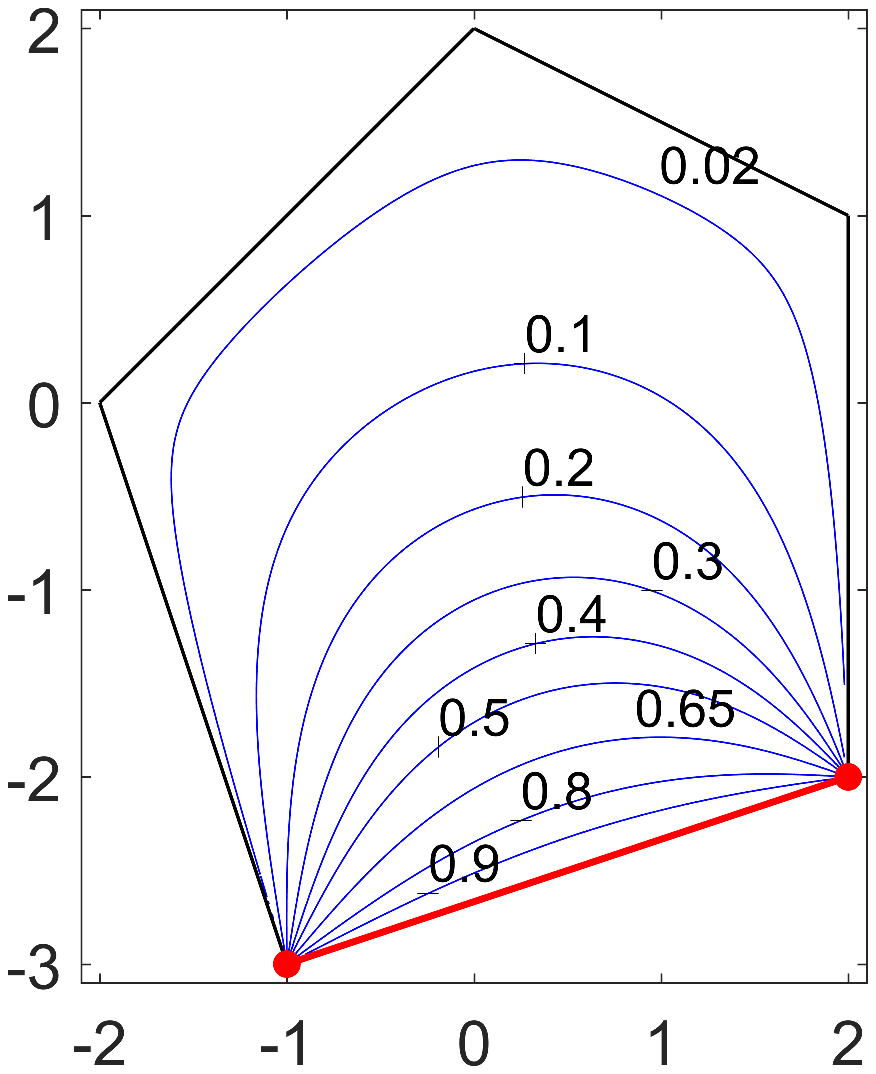}}
\hfill
}
\caption{The level curves of the function $\omega(z,L)$ for the polygon with $5$ sides.}
\label{fig:hm-5s}
\end{figure}

\nonsec{\bf Polygon with $5$ sides.}
As our first example, we consider the simply connected domain $G$ in the interior of the polygon with $5$ sides (the polygon shown in Figure~\ref{fig:hm-5s} and the vertices of the polygon are $2-2\i$, $2+\i$, $2\i$, $-2$, and $-1-3\i$).
We use the MATLAB function \verb|hm| with $n_s=2^{9}$ to compute the harmonic measure $\omega(z,L)$  of each side $L$ of the polygon with respect to the polygon domain $G$. The level curves of the function $\omega(z,L)$ are shown in Figure~\ref{fig:hm-5s}.

\begin{figure}[ht] %
\centerline{
\scalebox{0.45}{\includegraphics[trim=1.0cm 0.0cm 1.0cm 0.0cm,clip]{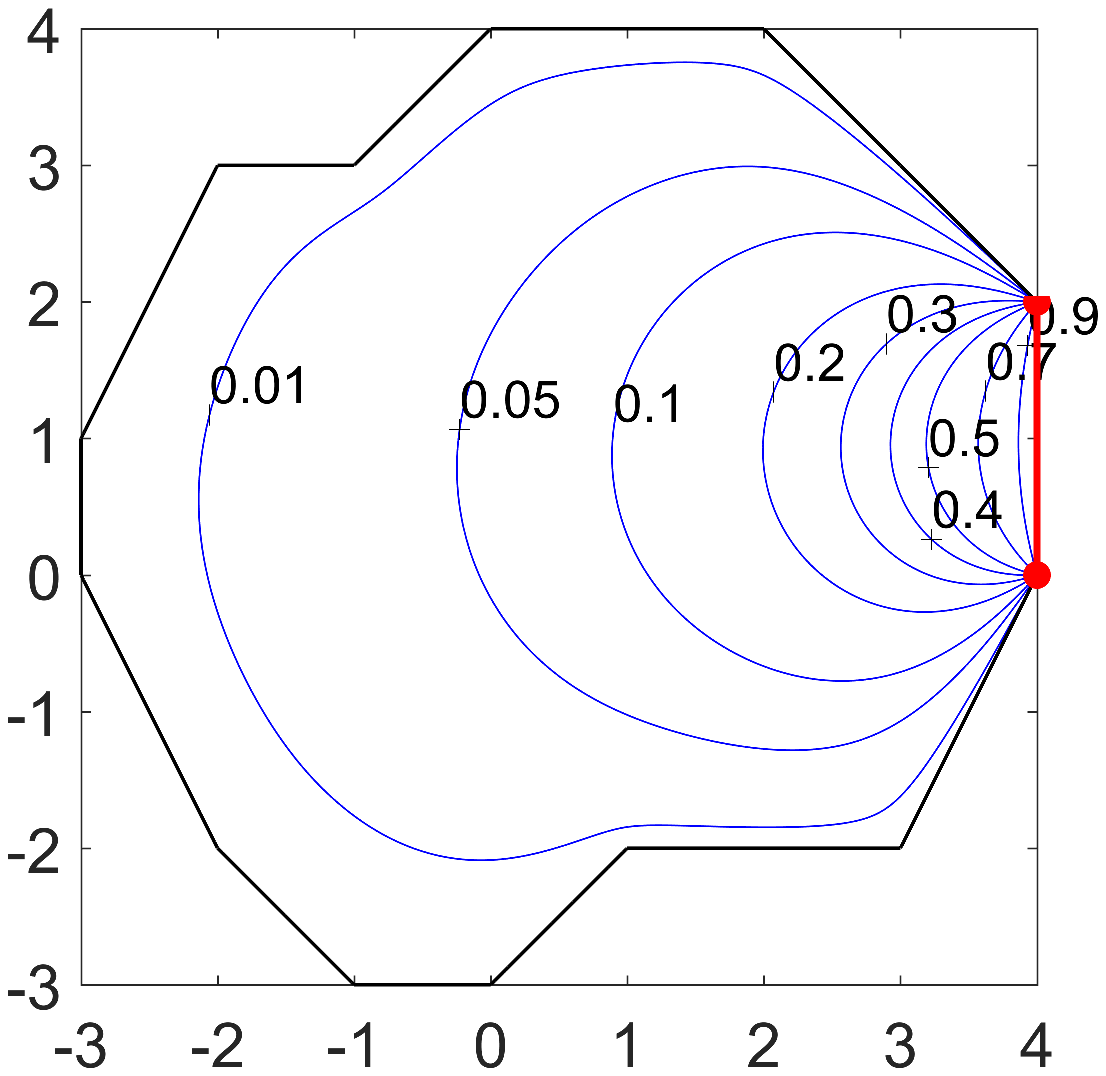}}
\hfill
\scalebox{0.45}{\includegraphics[trim=1.0cm 0.0cm 1.0cm 0.0cm,clip]{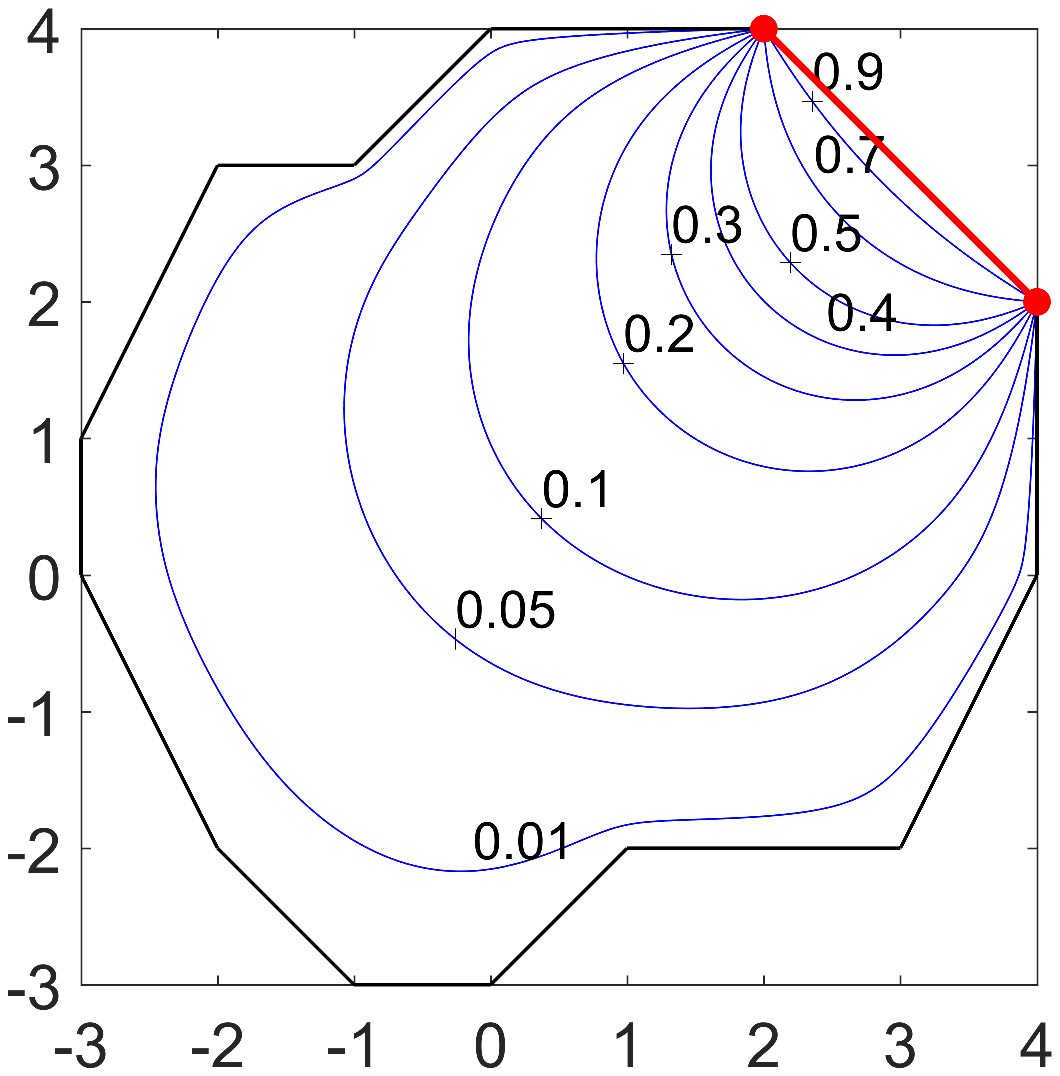}}
\hfill
\scalebox{0.45}{\includegraphics[trim=1.0cm 0.0cm 1.0cm 0.0cm,clip]{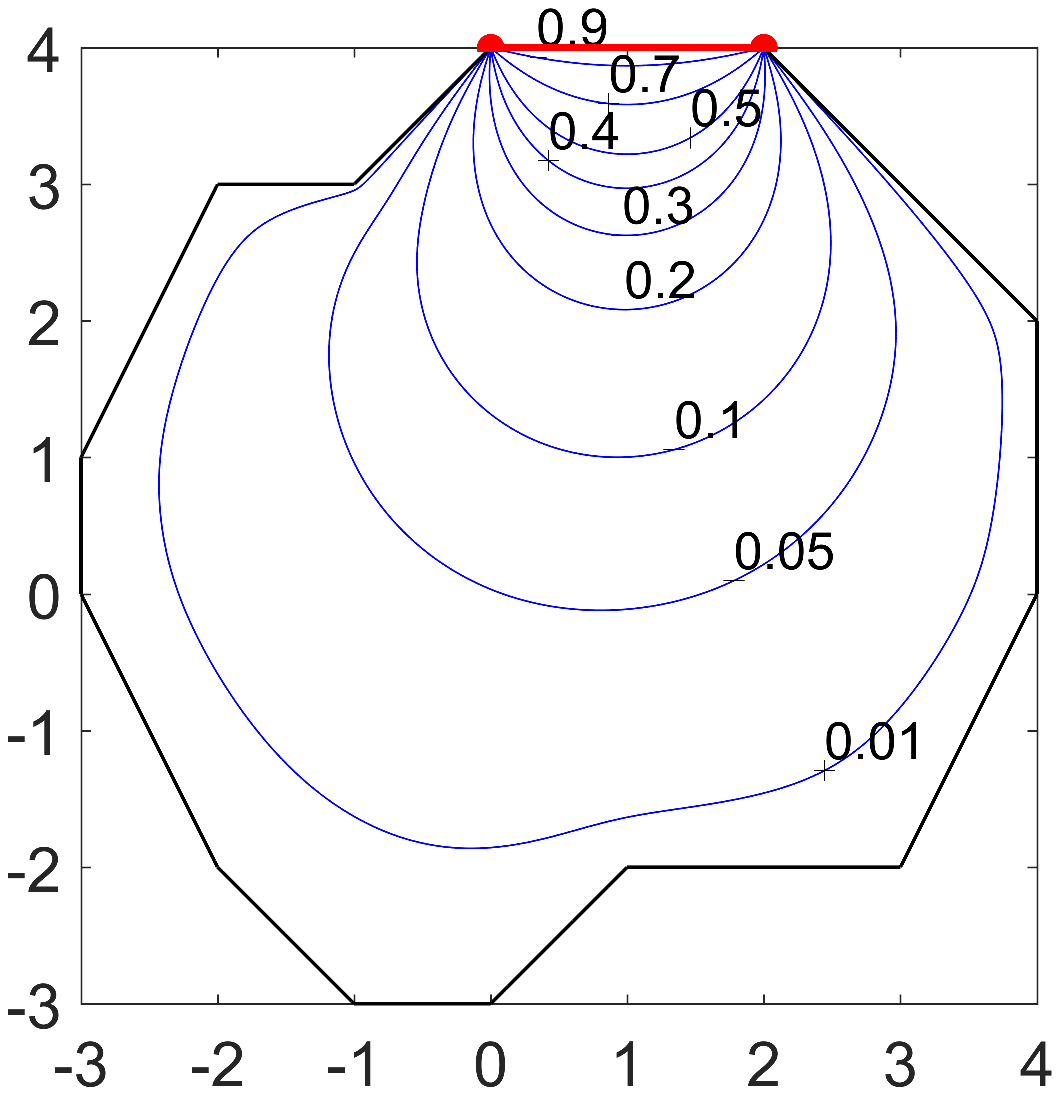}}
}
\centerline{
\hfill
\scalebox{0.45}{\includegraphics[trim=1.0cm 0.0cm 1.0cm 0.0cm,clip]{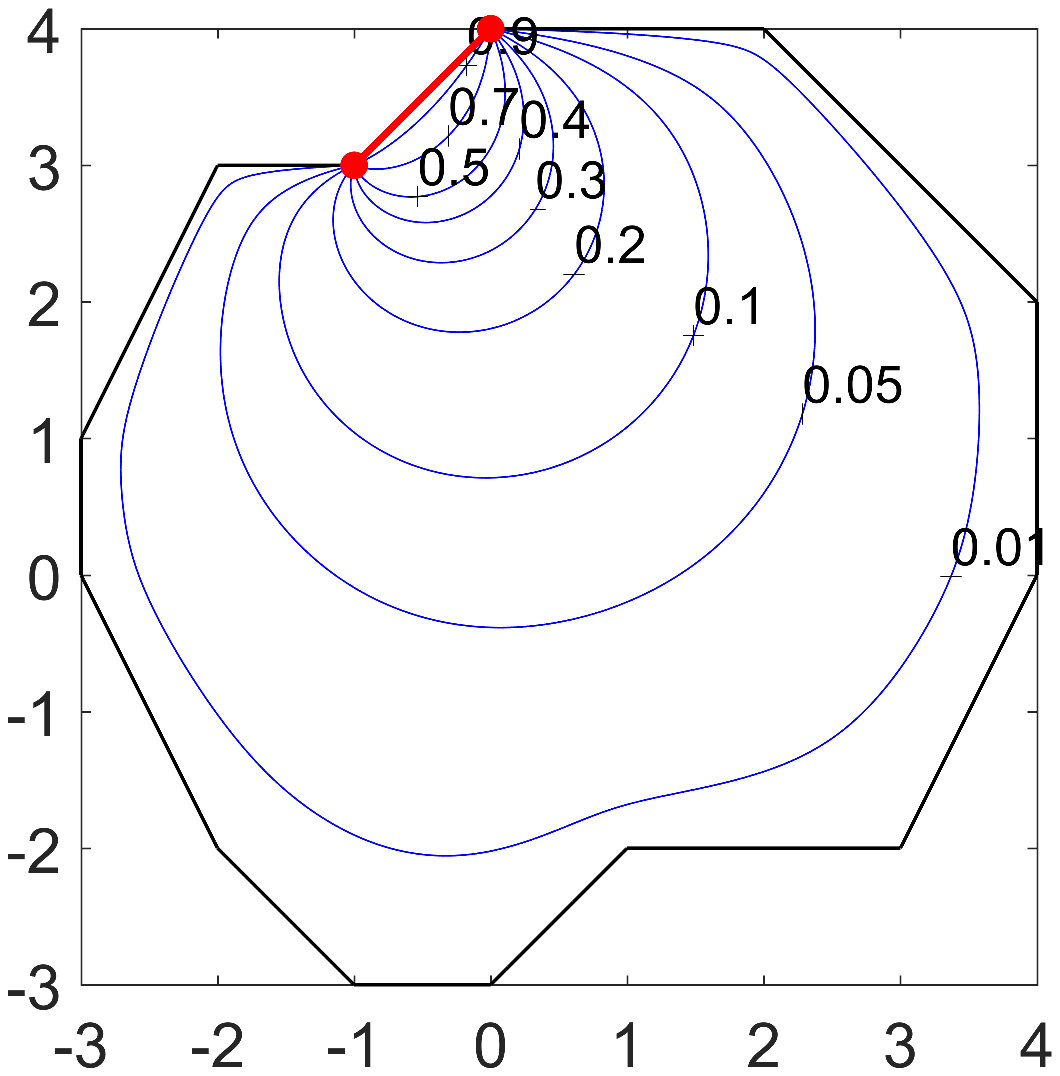}}
\hfill
\scalebox{0.45}{\includegraphics[trim=1.0cm 0.0cm 1.0cm 0.0cm,clip]{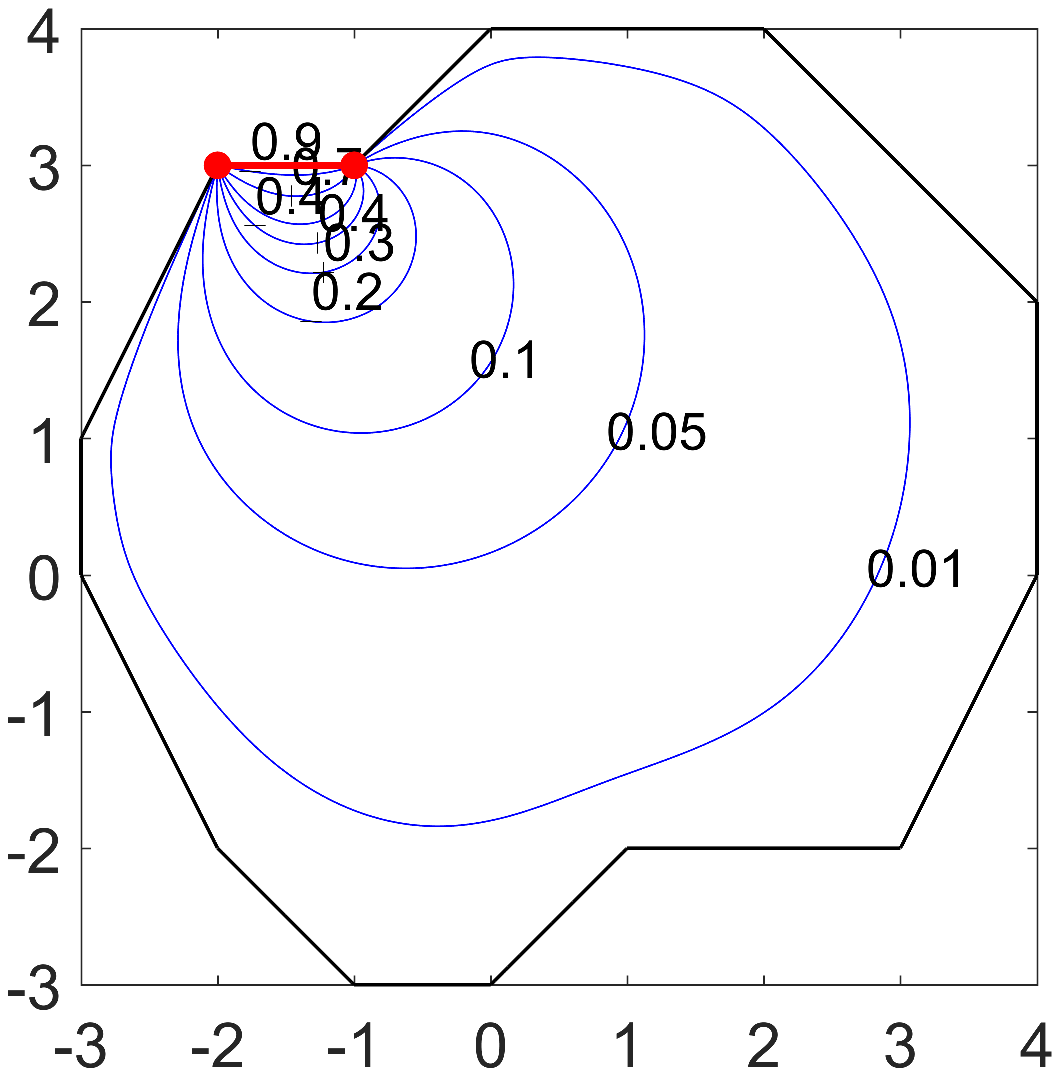}}
\hfill
\scalebox{0.45}{\includegraphics[trim=1.0cm 0.0cm 1.0cm 0.0cm,clip]{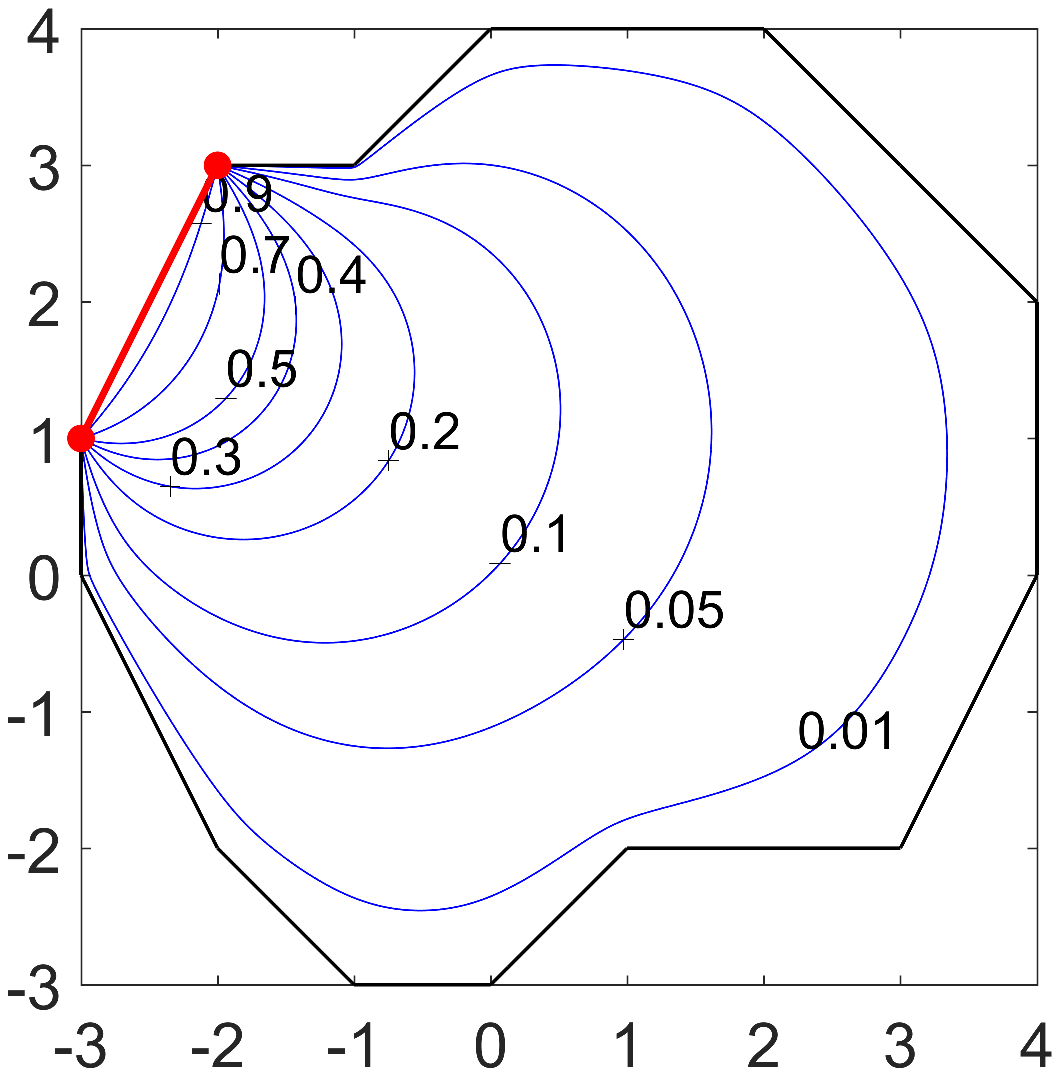}}
\hfill
}
\caption{The level curves of the function $\omega(z,L)$ for the polygon with $13$ sides.}
\label{fig:hm-13s}
\end{figure}

\nonsec{\bf Polygon with $13$ sides.}
For the second example, we consider the simply connected domain $G$ in the interior of the polygon with $13$ sides where the vertices of the polygon are $4$, $4+2\i$, $2+4\i$, $4\i$, $-1+3\i$, $-2+3\i$, $-3+\i$, $-3$, $-2-2\i$, $-1-3\i$, $-3\i$, $1-2\i$, and $3-2\i$. The MATLAB function \verb|hm| with $n_s=2^{9}$ is used to compute the harmonic measure $\omega(z,L)$ for each side $L$ of the polygon with respect to the polygon domain $G$. The level curves of the function $\omega(z,L)$ for the first $6$ sides are shown in Figure~\ref{fig:hm-13s}.

%% file: sec608.tex

\section{Quadrilateral Domains}

\nonsec{\bf Iterative method.}
Let $w=\Phi(z)$ be the conformal mapping from the interior of the unit circle $D=\{z\in\CC:|z|=1\}$ onto the interior of the rectangle 
\begin{equation}\label{eq:R_r}
R_r=\{w\: :\: 0<\Re w<1, \; 0<\Im w< r\}
\end{equation}
such that
\[
\Phi(z_1)=0,\quad \Phi(z_2)=1,\quad \Phi(z_3)=1+\i r,\quad \Phi(z_4)=\i r,
\]
where $z_1$, $z_2$, $z_3$, and $z_4$ are points on $\partial D$ (in the counterclockwise orientation) and $r>0$ is an undetermined positive real constant. The constant $r$ is known as the modulus of the quadrilateral $(D;z_1,z_2,z_3,z_4)$ and denoted by $M(D;z_1,z_2,z_3,z_4)$.
The modulus of the quadrilateral domains is invariant under conformal mappings, and hence
general bounded simply connected domains can be handled by mapping them onto the unit disk with the help of method presented in Section 2. 

If the domain $R_r$ is known (i.e., if $r$ is known), then we can map $R_r$ onto the unit disk using the method described in Section~\ref{sec:cm}. Let $\zeta=\Psi_1(w)$ be the conformal mapping from $R_r$ onto the disk $|\zeta|<1$ such that 
\[
\Psi_1(\alpha)=0, \quad \Psi'_1(\alpha)>0
\]
where $\alpha$ is a given point in $R_r$. The mapping function $z=\Psi_1(w)$ maps the points $0$, $1$, and $1+\i r$ on $\partial R_r$ onto points $\zeta_1$, $\zeta_2$, and $\zeta_3$ on the unit circle. Then the M\"obius transform
\[
z=\Psi_2(\zeta;\zeta_1,\zeta_2,\zeta_3)=z_3+\frac{(z_3-z_1)(z_2-z_3)(\zeta_2-\zeta_1)(\zeta-\zeta_3)}{(z_2-z_1)(\zeta_2-\zeta_3)(\zeta-\zeta_1)-(z_2-z_3)(\zeta_2-\zeta_1)(\zeta-\zeta_3)}
\]
maps the unit disk $|\zeta|<1$ onto the unit disk $|z|<1$ such that the points $\zeta_1$, $\zeta_2$ and $\zeta_3$ are mapped onto the points $z_1$, $z_2$ and $z_3$, respectively. 
Thus, the mapping function 
\begin{equation}\label{eq:Psi-quad}
z=\Psi(w)=\Psi_2(\Psi_1(w))
\end{equation}
maps the domain $R_r$ onto the unit disk $D$ such that the points $0$, $1$, $1+\i r$ are mapped onto the points $z_1$, $z_2$ and $z_3$, respectively. If $z=\Psi(w)$ maps also the point $\i r$ onto the point $z_4$, then $\Psi^{-1}$ will be the required map, i.e.,
\[
w=\Phi(z)=\Psi^{-1}(z).
\]

\begin{figure}[ht] %
\centerline{
\hfill
\scalebox{0.45}[0.45]{\includegraphics[trim=0cm 0cm 0cm 0cm,clip]{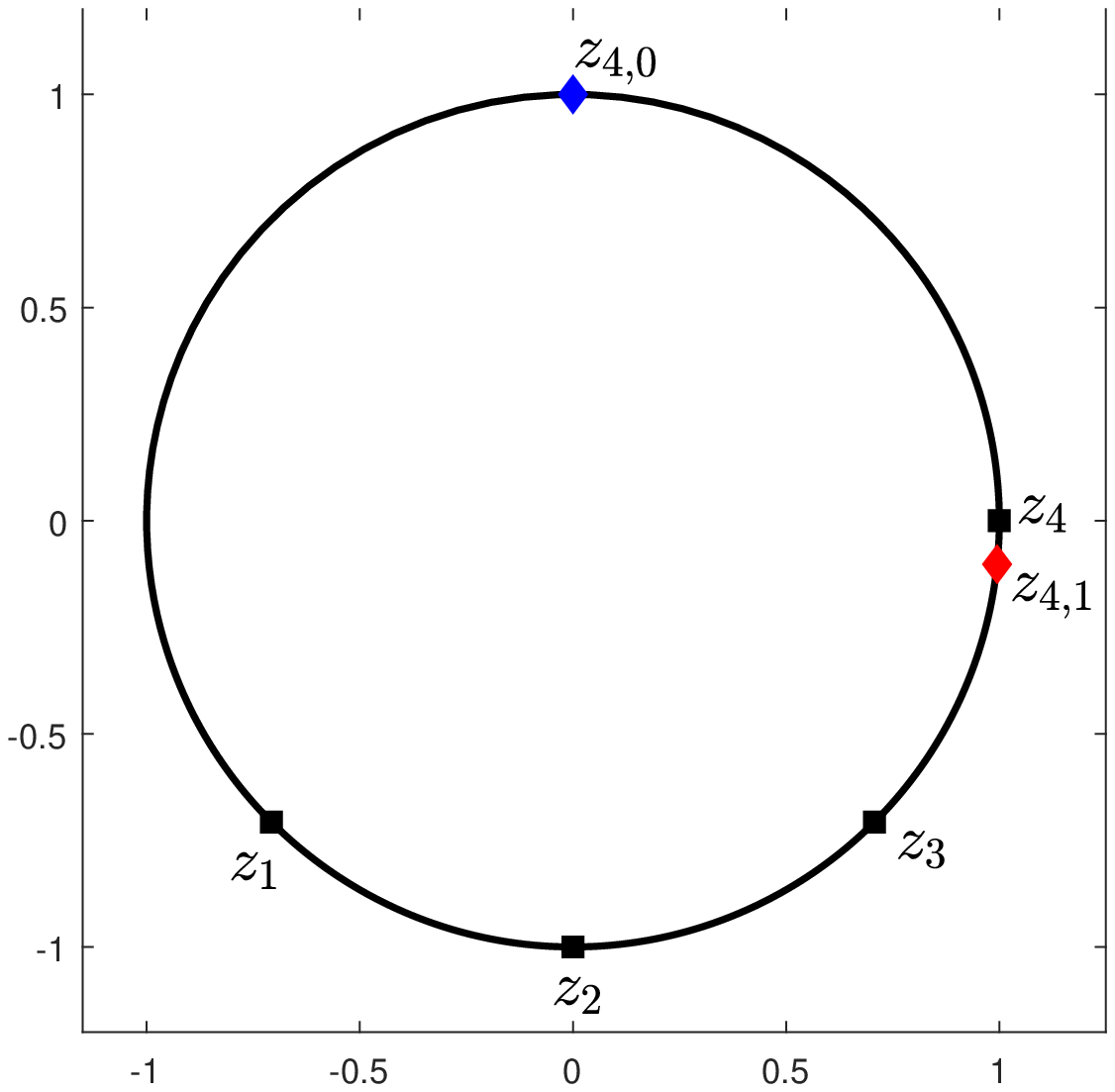}}
\hfill
\scalebox{0.45}[0.45]{\includegraphics[trim=0cm 0cm 0cm 0cm,clip]{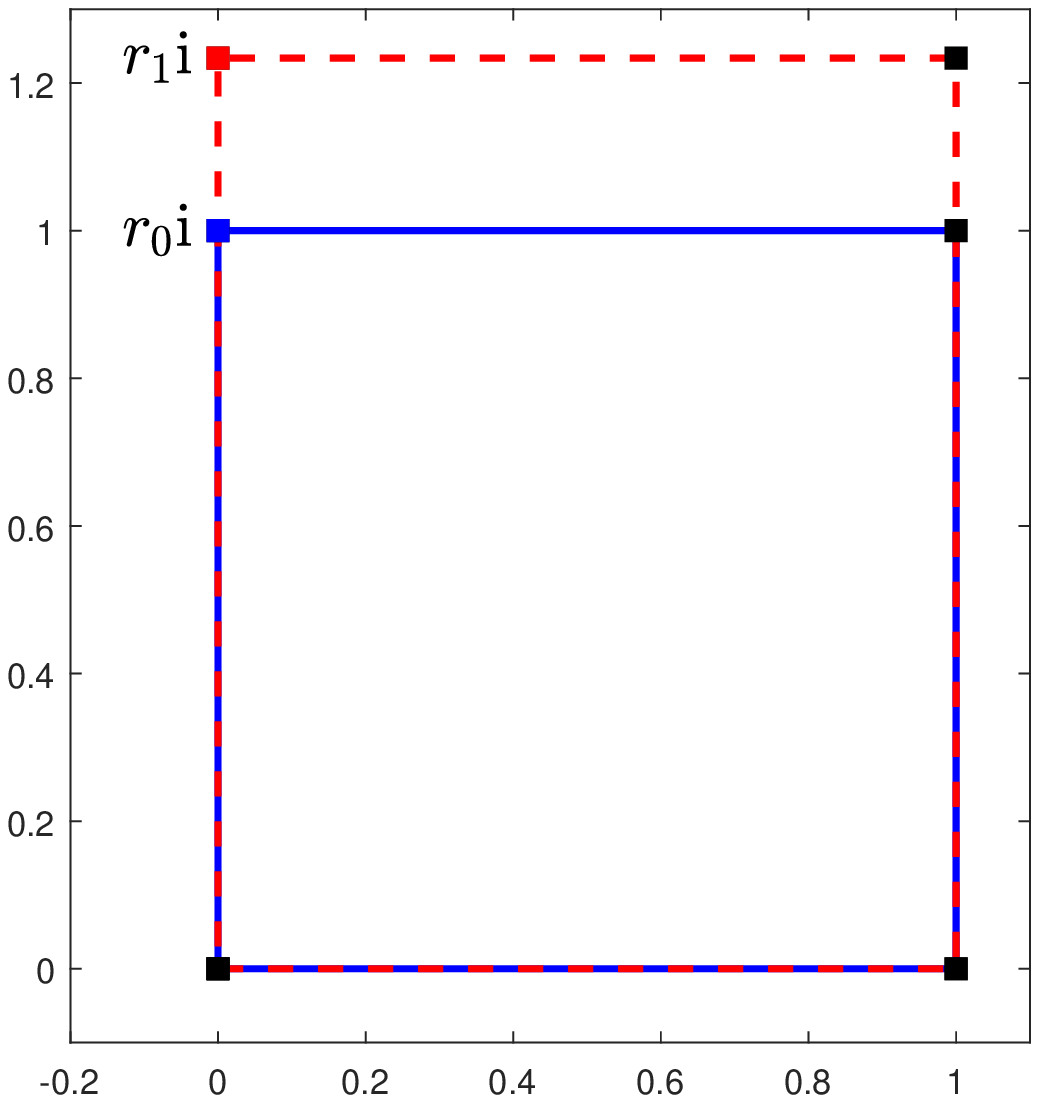}}
\hfill
}
\caption{The first two iterations.}
\label{fig:quad-0}
\end{figure}

In this section, for a given quadrilateral $(D;z_1,z_2,z_3,z_4)$, we present an iterative method for computing the unknown constant $r$ and the mapping function $z=\Psi(w)$ such that $\Psi(0)=z_1$, $\Psi(1)=z_2$, $\Psi(1+\i r)=z_3$, and $\Psi(\i r)=z_4$. 
First we choose an initial value $r_0=1$, then we use the function $\Psi$ to map $R_{r_0}$ to a quadrilateral $(D;z_1,z_2,z_3,z_{4,0})$ where $z_{4,0}$ is a point on the arc $[z_3,z_1]$ containing $z_4$ (see Figure~\ref{fig:quad-0}). 
The point $z_{4,0}$ could be on either side of $z_4$ on the arc $[z_3,z_1]$. 
We add a correction $\Delta_0$ to $r_0$ to get a new approximation $r_1$. 
Then we map $R_{r_1}$ to a quadrilateral $(D;z_1,z_2,z_3,z_{4,1})$ using the function $\Psi$. 
The point $z_{4,1}$ is now close to the point $z_4$. 
We continue with this iterative method to generate a sequence of approximation ${r_0,r_1,r_2,\ldots}$ and the mapping function $\Psi$ maps the rectangle $R_{r_{k}}$ to a quadrilateral $(D;z_1,z_2,z_3,z_{4,k})$. 
We stop the iteration when the distance (on the unit circle) between the two points $z_4$ and $z_{4,k}$ is small. Then, we consider $r_{k}$ as an approximation to $r$. 

Since, for each iteration $k$, the point $z_{4,k}$ is on the arc $[z_3,z_1]$, we can always choose suitable corrections $\Delta_{k}$ to ensure the convergence of the iterative method. In this paper, for $k\ge1$, we choose
\begin{equation}\label{eq:Delta-k}
\Delta_k=\frac{1}{2\pi}\arg\left(\frac{z_{4,k}}{z_4}\right).
\end{equation}

To accelerate the convergence of the iterative method, we introduce a factor $\delta_k$ and we calculate $r_k$ using the formula
\begin{equation}\label{eq:r-k}
r_{k+1}=r_{k}+\delta_{k}\Delta_{k}, \quad k\ge 0,
\end{equation}
where $r_0=1$, $\delta_0=\delta_1=1$, and for $k\ge2$,
\begin{equation}\label{eq:delta-k}
\delta_k =
\begin{cases}
2\delta_{k-1},  & \mbox{if $\arg\left(\frac{z_{4,k-2}}{z_4}\right)\arg\left(\frac{z_{4,k-1}}{z_4}\right)>0$ and $\arg\left(\frac{z_{4,k-1}}{z_4}\right)\arg\left(\frac{z_{4,k}}{z_4}\right)>0$}, \\
\delta_{k-1}/2, & \mbox{if $\arg\left(\frac{z_{4,k-2}}{z_4}\right)\arg\left(\frac{z_{4,k-1}}{z_4}\right)<0$ and $\arg\left(\frac{z_{4,k-1}}{z_4}\right)\arg\left(\frac{z_{4,k}}{z_4}\right)<0$}, \\
\delta_{k-1}, & \mbox{otherwise}.
\end{cases}
\end{equation}
In other words, when the three points $z_{4,k-2}$, $z_{4,k-1}$ and $z_{4,k}$ are in the same side of $z_4$, we double $\delta_{k-1}$ to increase the correction added to $r_{k}$ and so push $z_{4,k}$ toward $z_4$. However, when the three points $z_{4,k-2}$, $z_{4,k-1}$ and $z_{4,k}$ oscillate around $z_4$, we bisect $\delta_{k-1}$ to reduce the correction added to $r_{k}$.
To avoid getting very long rectangle or very narrow rectangle during the iterations, we do not allow $\delta_k\Delta_k$ to be more than $0.2r_k$ or less than $-0.2r_k$. 

\nonsec{\bf Algorithm.}\label{sc:algorithm}
The above iterative method is summarized as follows.

\noindent{\bf Initialization:}\\
Set $r_0=1$, $\delta_0=\delta_1=1$. \\
\noindent{\bf Iterations:} \\
For $k=1,2,3,\ldots$, where $k$ denotes the iteration number:
\begin{itemize}
	\item Map the domain interior to the rectangle with the vertices $\{0,1,1+\i r_{k-1},\i r_{k-1}\}$ onto the unit disk $D$ by the function $\Psi$ in~\eqref{eq:Psi-quad} such $\Psi(0)=z_1$, $\Psi(1)=z_2$, $\Psi(1+\i r_{k-1})=z_3$.
	\item Let $z_{4,k-1}=\Psi(\i r_{k-1})$. 
	\item Compute $\Delta_{k-1}$ from~\eqref{eq:Delta-k}.
	\item For $k\ge2$, compute $\delta_{k-1}$ from~\eqref{eq:delta-k}.
	\item If $\delta_{k-1}\Delta_{k-1}>0.2r_{k-1}$, then set $\delta_{k-1}\Delta_{k-1}=0.2r_{k-1}$ and $\delta_{k-1}=\delta_{k-1}/2$.
	\item If $\delta_{k-1}\Delta_{k-1}<-0.2r_{k-1}$, then set $\delta_{k-1}\Delta_{k-1}=-0.2r_{k-1}$ and $\delta_{k-1}=\delta_{k-1}/2$.
	\item The approximate value of $r$ is updated through $r_k=r_{k-1}+\delta_{k-1}\Delta_{k-1}$.
  \item Stop the iteration if $|r_k-r_{k-1}|<\varepsilon$ or $k>{\tt Max}$ where ${\tt Max}$ is the maximum number of allowed iterations and $\varepsilon$ is a given tolerance. 
\end{itemize}

In our numerical experiments, we choose ${\tt Max}=50$ and $\varepsilon=0.5\times 10^{-13}$. The iterative method produces a sequence of numbers $r_0,r_1,r_2,r_3,\ldots$ which converges to the required constant $r$. The iterative method provides us also with a conformal map $z=\Psi(w)$ from $R_r$ onto the given domain $D$. Then the required map $\Phi$ is given by
\[
w=\Phi(z)=\Psi^{-1}(z).
\]
The numerical examples presented in this section show that the iterative method converges for several examples. However, no proof of convergence is available so far.

\nonsec{\bf Examples.}
We consider the computing of the conformal mapping from the quadrilateral domains
\[
(D;e^{\i\theta_1},e^{\i\theta_2},e^{\i\theta_3},e^{\i\theta_4}),
\]
onto rectangular domains for the following values of $\theta_1$, $\theta_2$, $\theta_3$, and $\theta_4$:
\begin{enumerate}
	\item $Q_1:$ $\theta_1=-\pi$, $\theta_2=-0.5\pi$, $\theta_3=0$, $\theta_4=0.5\pi$.
	\item $Q_2:$ $\theta_1=-0.5\pi$, $\theta_2=-0.25\pi$, $\theta_3=0.25\pi$, $\theta_4=0.5\pi$.
	\item $Q_3:$ $\theta_1=-0.5005\pi$, $\theta_2=-0.4995\pi$, $\theta_3=0.4995\pi$, $\theta_4=0.5005\pi$.
	\item $Q_4:$ $\theta_1=-\pi$, $\theta_2=-0.0001\pi$, $\theta_3=0$, $\theta_4=0.5\pi$.
\end{enumerate}

\begin{figure}[ht] %
\centerline{
\hfill
\scalebox{0.45}[0.45]{\includegraphics[trim=0cm 0cm 0cm 0cm,clip]{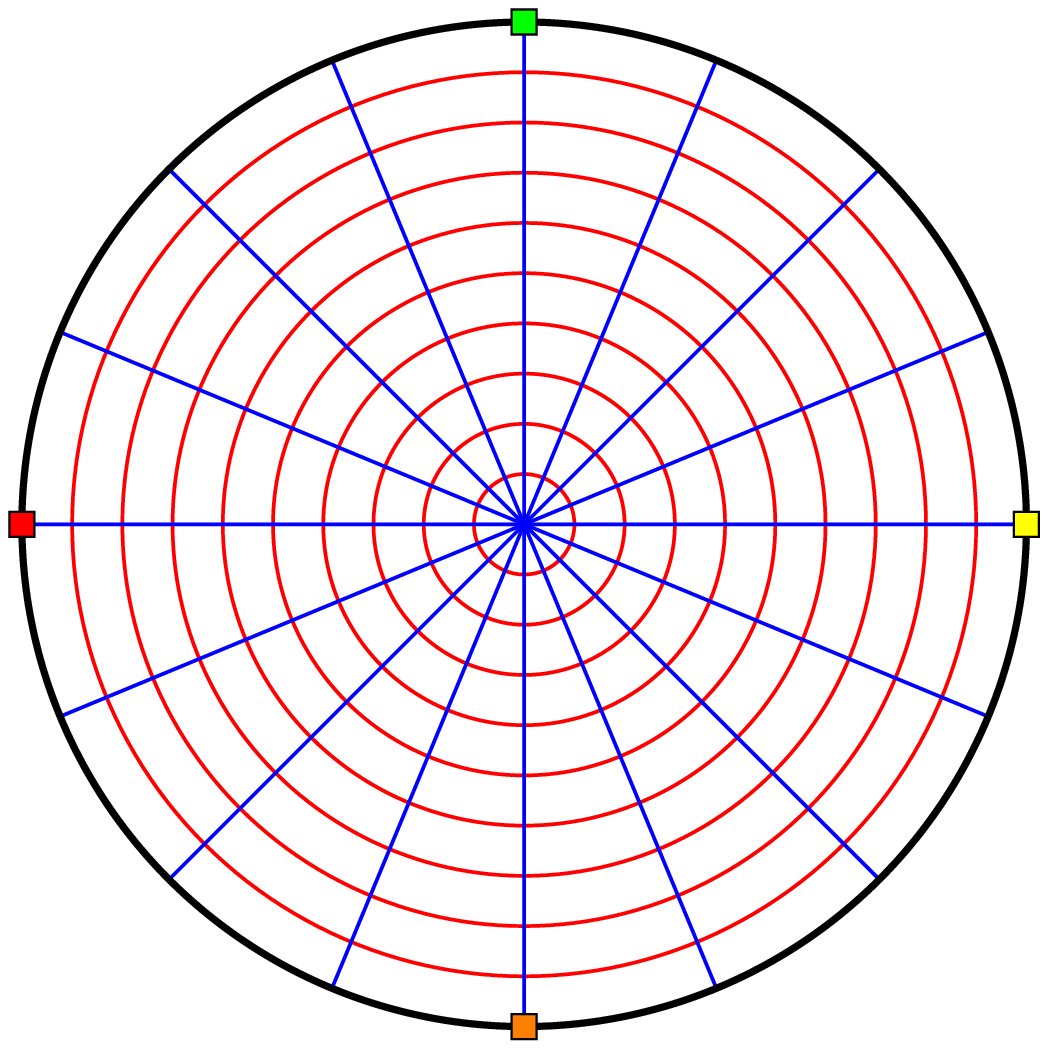}}
\hfill
\scalebox{0.45}[0.45]{\includegraphics[trim=0cm 0cm 0cm 0cm,clip]{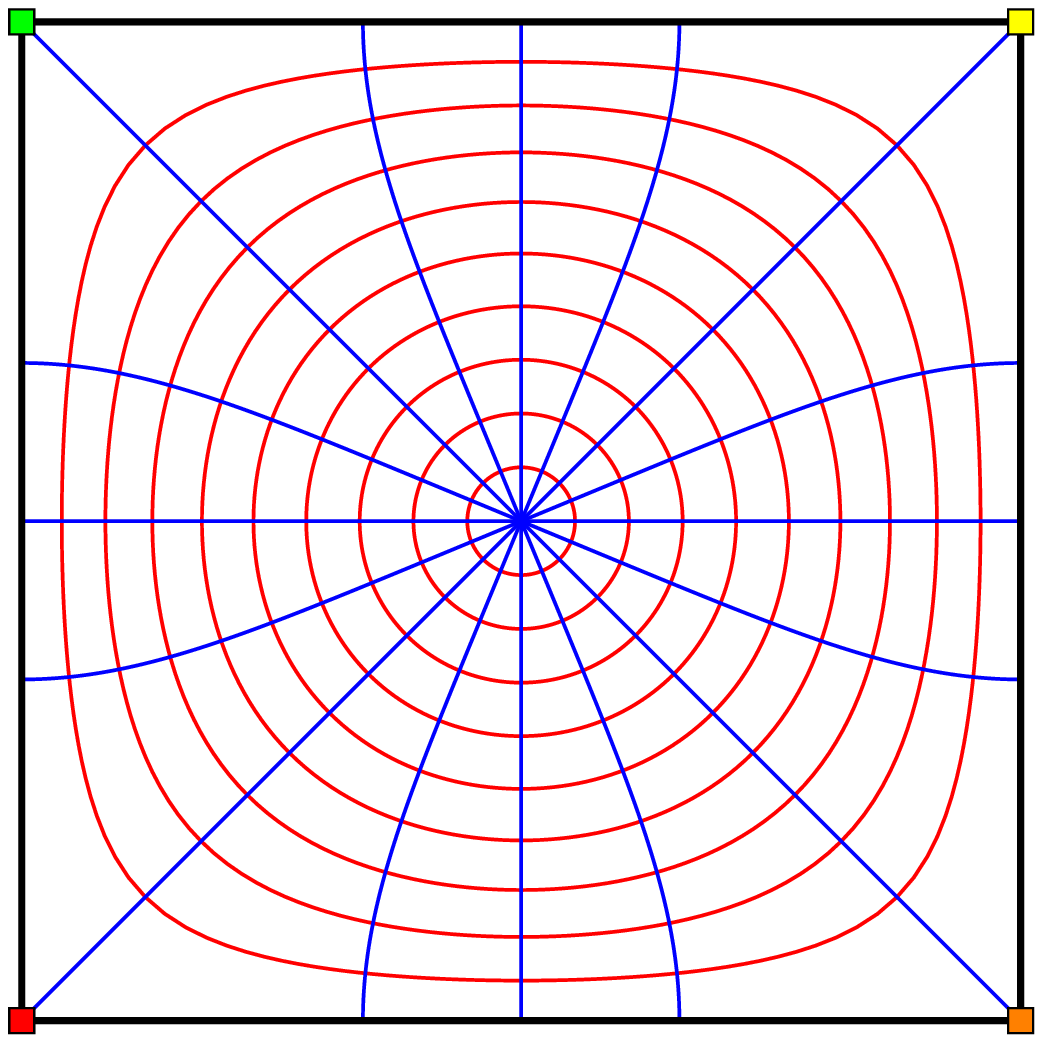}}
\hfill
}
\caption{The quadrilateral domain $Q_1$ and its image.}
\label{fig:quad-1}
\end{figure}

\begin{figure}[ht] %
\centerline{
\hfill
\scalebox{0.45}[0.45]{\includegraphics[trim=0cm 0cm 0cm 0cm,clip]{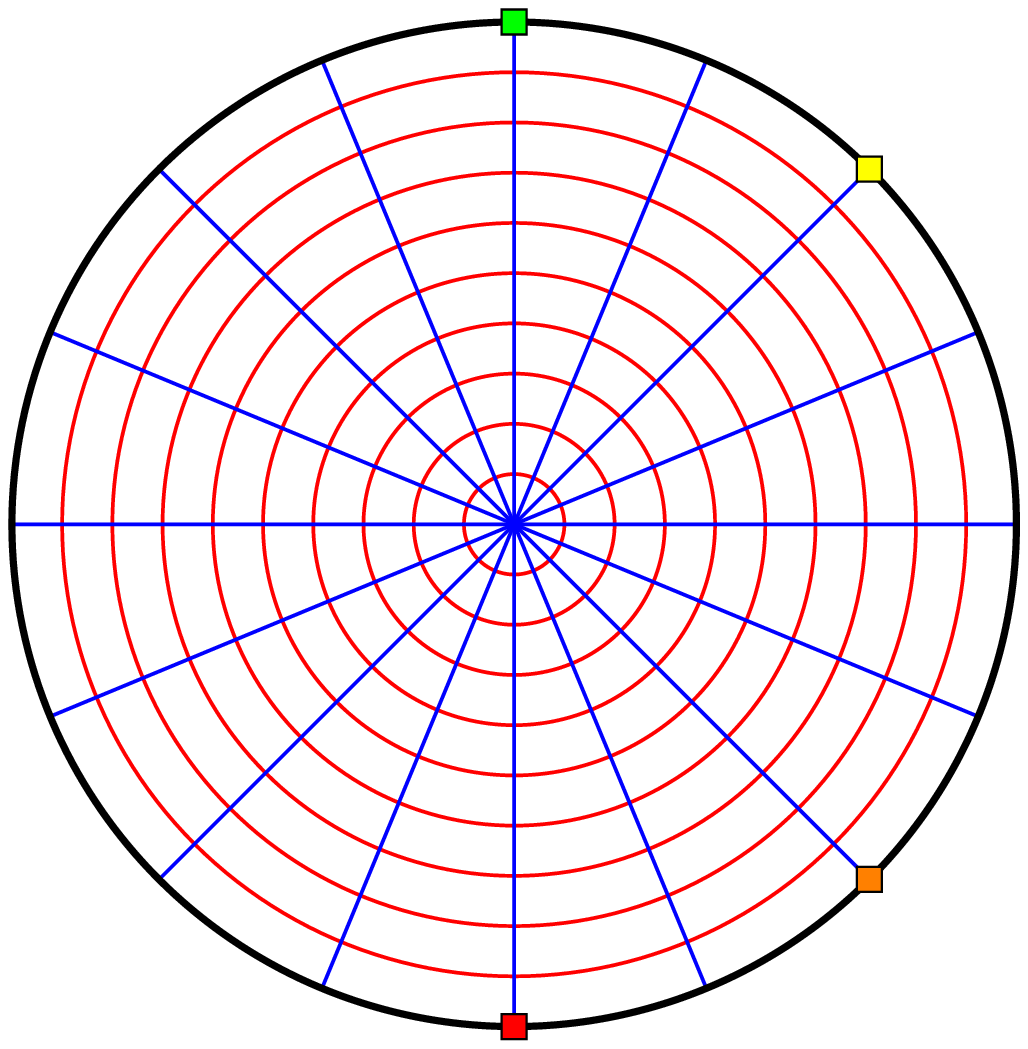}}
\hfill
\scalebox{0.45}[0.45]{\includegraphics[trim=0cm 0cm 0cm 0cm,clip]{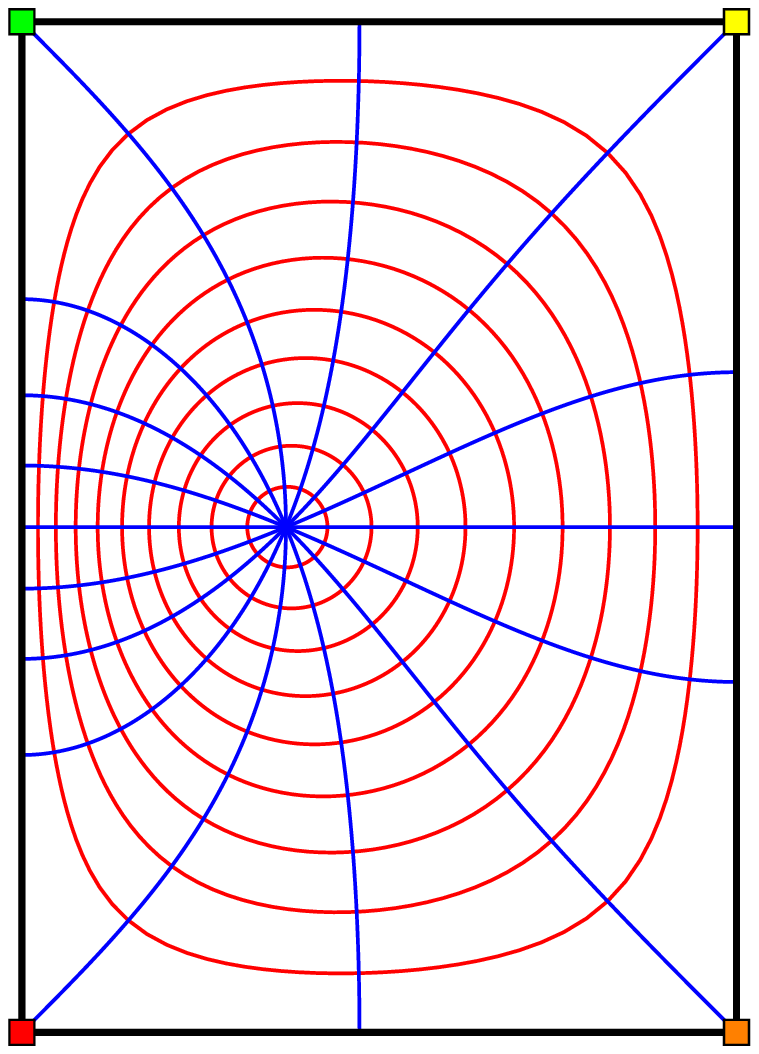}}
\hfill
}
\caption{The quadrilateral domain $Q_2$ and its image.}
\label{fig:quad-2}
\end{figure}

\begin{figure}[ht] %
\centerline{
\hfill
\scalebox{0.45}[0.45]{\includegraphics[trim=0cm 0cm 0cm 0cm,clip]{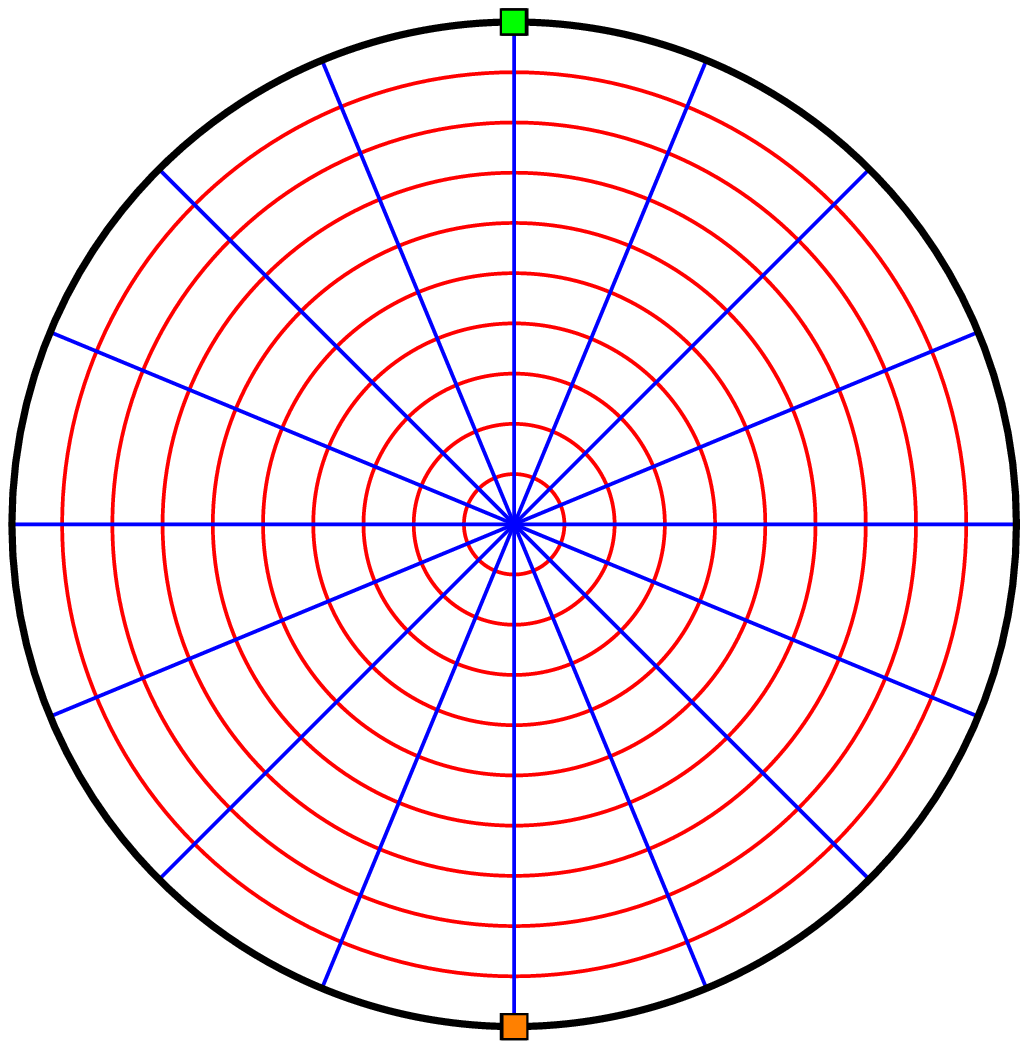}}
\hfill
\scalebox{0.45}[0.45]{\includegraphics[trim=0cm 0cm 0cm 0cm,clip]{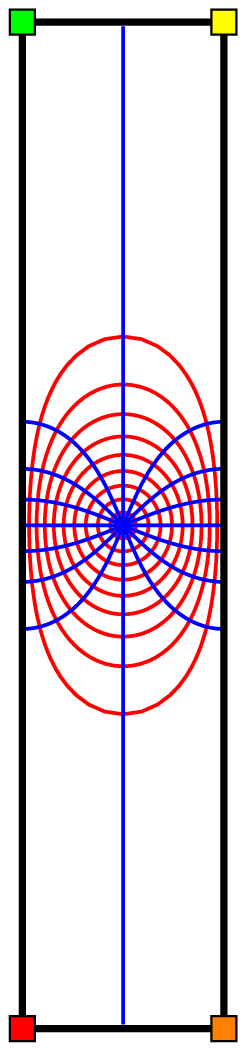}}
\hfill
}
\caption{The quadrilateral domain $Q_3$ and its image.}
\label{fig:quad-3}
\end{figure}

\begin{figure}[ht] %
\centerline{
\hfill
\scalebox{0.45}[0.45]{\includegraphics[trim=0cm 0cm 0cm 0cm,clip]{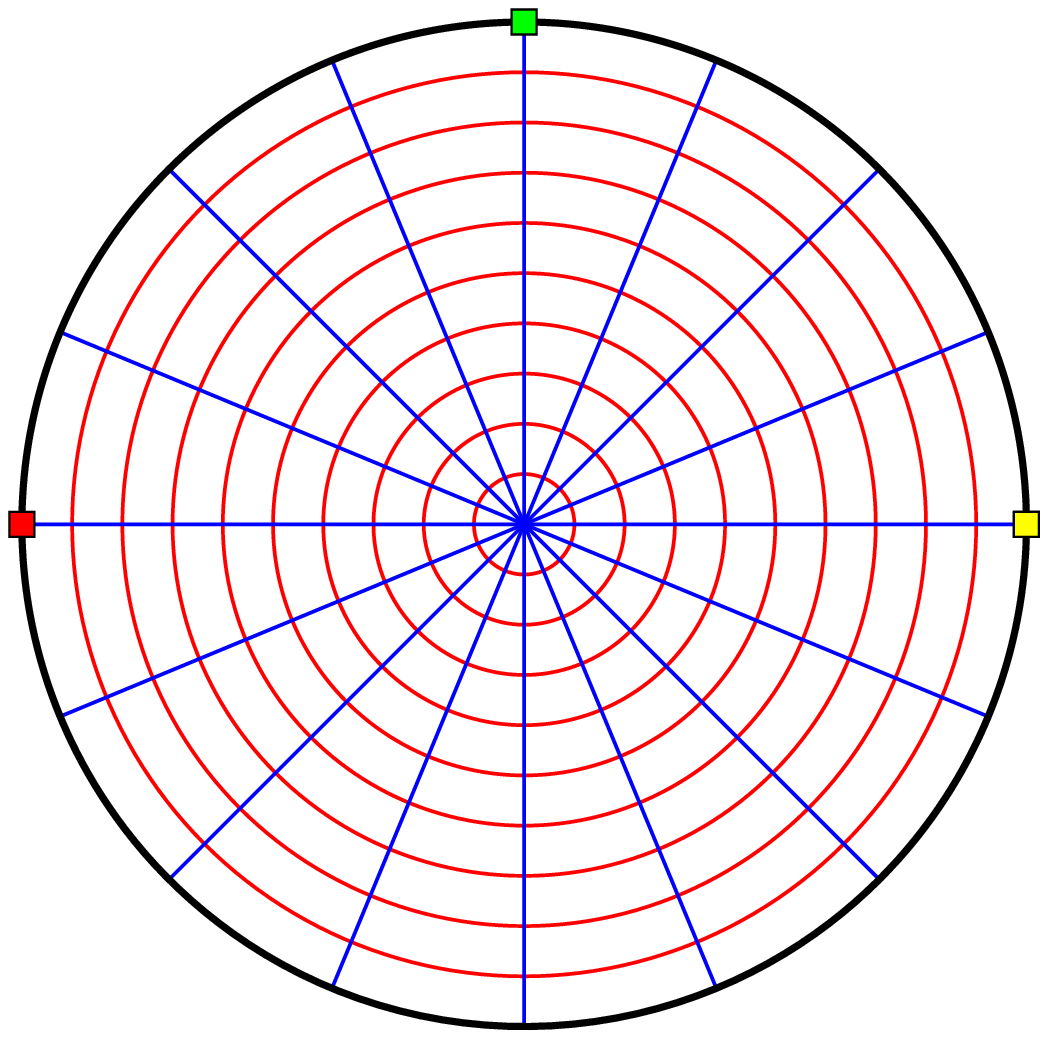}}
\hfill
\scalebox{0.45}[0.45]{\includegraphics[trim=0cm -3cm 0cm 0cm,clip]{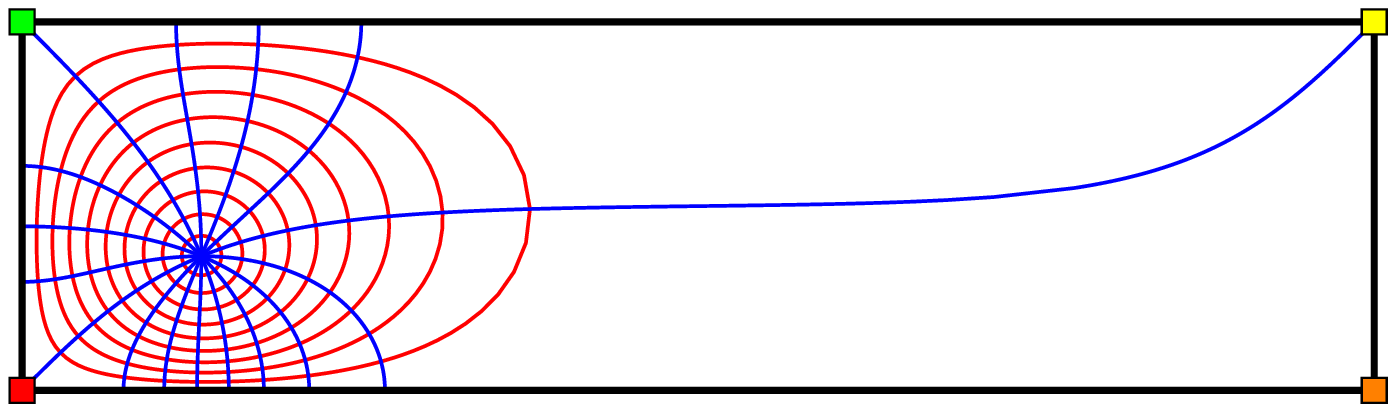}}
\hfill
}
\caption{The quadrilateral domain $Q_4$ and its image.}
\label{fig:quad-4}
\end{figure}

The values of the modulus $r=M(Q_j)$, $j=1,2,3,4$, the number of iterations required for convergence, and the total CPU time (sec) required for convergence are listed in Table~\ref{tab:quad-map}. For the four domains, we use $n=2^{11}$. Orthogonal polar grids in the circular domains and their images under the conformal mapping are shown in Figures~\ref{fig:quad-1}--\ref{fig:quad-4}. The points $z_1,z_2,z_3,z_4$ on the unit circle and their images on the rectangle are shown as small colored squares such that a point $z_k$ and its image has the same color. 
For $Q_3$, $z_1=e^{-0.5005\pi}$ and $z_2=e^{-0.4995\pi}$ which are very close to each other. Similarly, $z_3=e^{0.4995\pi}$ and $z_4=e^{0.5005\pi}$ are very close to each other. The length of the arcs between $z_1$ and $z_2$ and between $z_3$ and $z_4$ is $0.001\pi$. Thus, we can not distinguish between $z_1$ and $z_2$ and between $z_3$ and $z_4$ in Figure~\ref{fig:quad-3} (left). 
The small arc between $z_1$ and $z_2$ is mapped by the conformal mapping to the lower side of the rectangle. Similarly, the small arc between $z_3$ and $z_4$ is mapped by the conformal mapping to the upper side of the rectangle. 
Although these arcs are too small, the proposed iterative method converges after only $36$ iterations.
In $Q_4$, the two points $z_2=e^{-0.0001\pi}$ and $z_3=1$ are very close to each other where the length of the arcs between them is $0.0001\pi$, and hence we can not distinguish between $z_2$ and $z_3$ in Figure~\ref{fig:quad-4} (left). The small arc between $z_2$ and $z_3$ is mapped by the conformal mapping to the right side of the rectangle. The proposed iterative method converges after only $40$ iterations.

\begin{figure}[ht] %
\centerline{
\scalebox{0.45}[0.45]{\includegraphics[trim=0cm 0cm 0cm 0cm,clip]{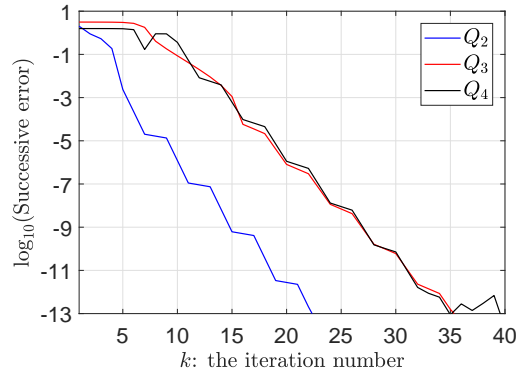}}
}
\caption{The successive error $|r_{k}-r_{k-1}|$ for the three domains $Q_2$, $Q_3$, and $Q_4$ vs the number of iteration $k$.}
\label{fig:quad-error}
\end{figure}

For the three domains $Q_2$, $Q_3$, and $Q_4$, the error per iteration is shown in Figure~\ref{fig:quad-error}. For $Q_3$ and $Q_4$, we have points on the unit circle that are very close to each other. This explains why the number of iterations for $Q_3$ and $Q_4$ is larger than the number of iterations for $Q_2$. For $Q_1$, the method converges after only one iteration since the exact value of $r$ is $1$ which is the same as our initial value $r_0$.

\begin{table}[h]
\caption{The numerical results.}
\label{tab:quad-map}%
\begin{tabular}{l|l|l|l}\hline
 Domain & $r$               & number of iterations & total CPU time \\ \hline
 $Q_1$  & 1                 & 1                    & 0.6  \\
 $Q_2$  & 1.41421356237738  & 23                   & 6.5  \\
 $Q_3$  & 4.99266938932358  & 36                   & 10.4  \\
 $Q_4$  & 0.272437506734334 & 40                   & 12.1  \\
\hline
\end{tabular}
\end{table}

\begin{figure}[ht] %
\centerline{
\scalebox{0.45}[0.45]{\includegraphics[trim=0cm 0cm 0cm 0cm,clip]{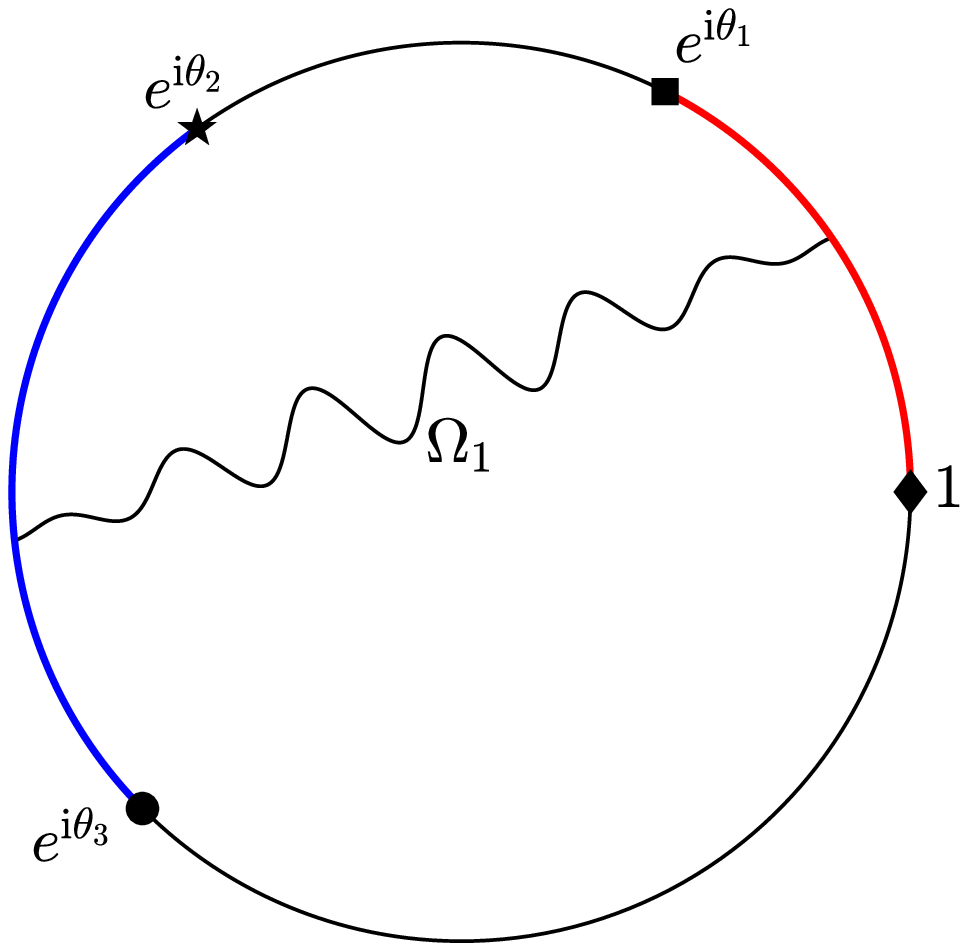}}
\hfill
\scalebox{0.45}[0.45]{\includegraphics[trim=0cm 0cm 0cm 0cm,clip]{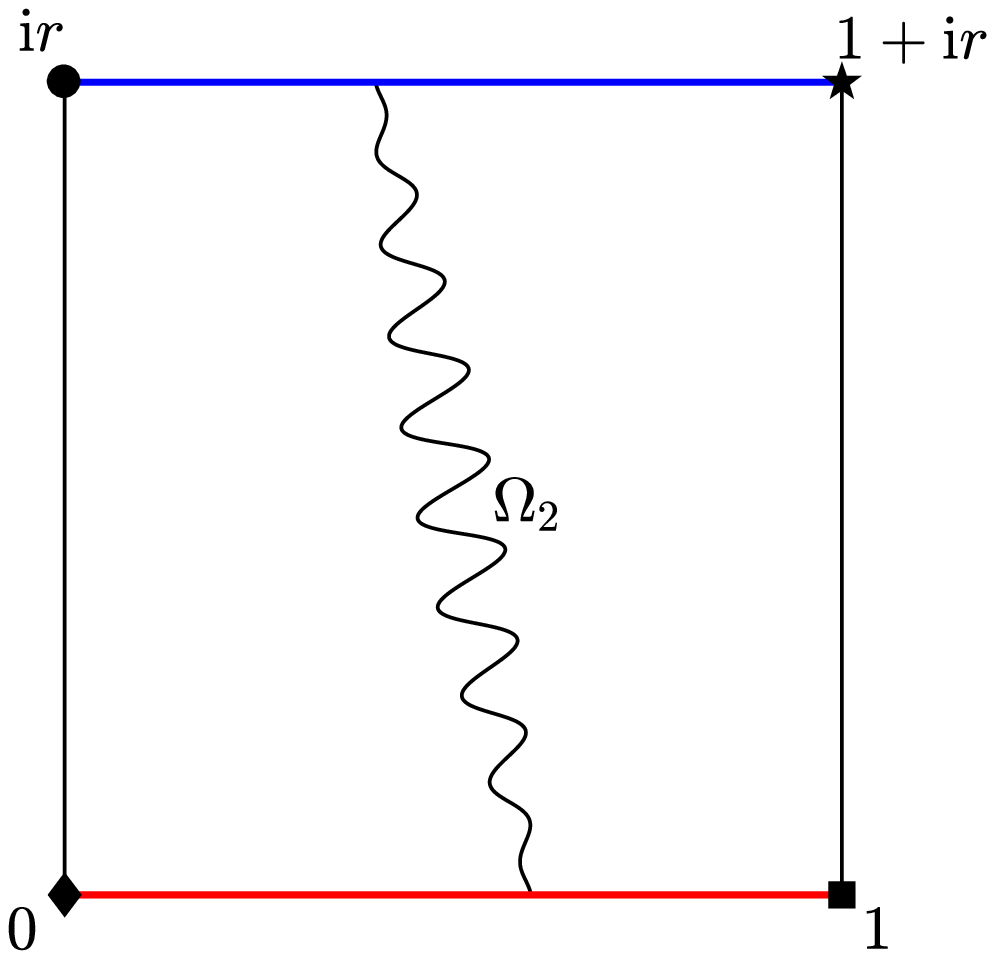}}
\hfill
\scalebox{0.45}[0.45]{\includegraphics[trim=0cm 0cm 0cm 0cm,clip]{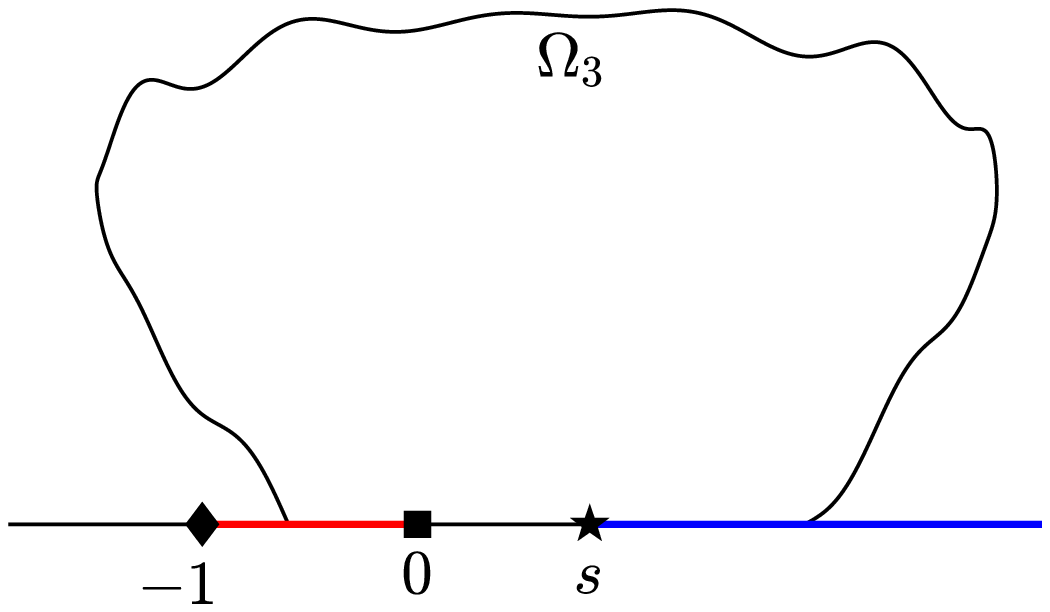}}
}
\caption{The quadrilateral domain (left), the rectangular domain (center), and the upper half-plane (right).}
\label{fig:quad-exact}
\end{figure}

\nonsec{\bf Explicit formula for the modulus.}\label{ex:explicit}
Consider the quadrilateral domain
\[
(D;1,e^{\i\theta_1},e^{\i\theta_2},e^{\i\theta_3}),
\]
which can be mapped conformally onto the rectangular domain $R_r=\{w\: :\: 0<\Re w<1, \; 0<\Im w< r\}$ such that the point $1$ is mapped to $0$, $e^{\i\theta_1}$ is mapped to $1$, $e^{\i\theta_2}$ is mapped to $1+r\i$, and $e^{\i\theta_3}$ is mapped to $r\i$. The quadrilateral domain can be mapped also by M\"obius transform onto the upper half-plane such that the point $1$ is mapped to $-1$, $e^{\i\theta_1}$ is mapped to $0$, $e^{\i\theta_2}$ is mapped to a positive real number $s$, and $e^{\i\theta_3}$ is mapped to $\infty$. See Figure~\ref{fig:quad-exact}. 

The exact value of the positive constant $s$ can be obtained in terms of the values of $\theta_1$, $\theta_2$, and $\theta_3$. For distinct points $z_1$, $z_2$, $z_3$, and $z_4$ in $\CC$, we define the absolute (cross) ratio by~\cite{be}
\[
|z_1,z_2,z_3,z_4|=\frac{|z_1-z_3||z_2-z_4|}{|z_1-z_2||z_3-z_4|}.
\]
This definition can be extended if $z_4=\infty$ by taking the limit, i.e.,
\[
|z_1,z_2,z_3,\infty|=\frac{|z_1-z_3|}{|z_1-z_2|}.
\]
Thus, for the four points $1$, $e^{\i\theta_1}$, $e^{\i\theta_2}$, and $e^{\i\theta_3}$ on the unit circle, we have
\[
|1,e^{\i\theta_1},e^{\i\theta_2},e^{\i\theta_3}|=\frac{\sin(\theta_2/2)}{\sin(\theta_1/2)}\frac{\sin((\theta_3-\theta_1)/2)}{\sin((\theta_3-\theta_2)/2)}.
\]
Similarly, for the four points $-1$, $0$, $s$, and $\infty$ on the real line, we have
\[
|-1,0,s,\infty|=1+s.
\]
An important property of M\"obius transformations is that they preserve the absolute ratios~\cite{be}, thus
\[
|-1,0,s,\infty|=|1,e^{\i\theta_1},e^{\i\theta_2},e^{\i\theta_3}|,
\]
and hence the exact value of $1+s$ is given by the formula
\begin{equation}\label{eq:r-sp1}
1+s=\frac{\sin(\theta_2/2)}{\sin(\theta_1/2)}\frac{\sin((\theta_3-\theta_1)/2)}{\sin((\theta_3-\theta_2)/2)}.
\end{equation}

Let $\Omega_1$ be the family of curves lying in $D$ and joining the sets $(1,e^{\i\theta_1})$ and $(e^{\i\theta_2},e^{\i\theta_3})$ (see Figure~\ref{fig:quad-exact}). 
Similarly, let $\Omega_2$ be the family of curves lying in $R_r$ and joining the sets $(1,1+\i r)$ and $(\i r,0)$, and let $\Omega_3$ be the family of curves lying in the upper half-plane and joining the sets $(-1,0)$ and $(s,\infty)$. The modulus is invariant under conformal mapping~\cite{ah,LV,vu88} and hence~\cite[p.~20]{du}
\[
M(\Omega_1)=M(\Omega_2)=M(\Omega_3)=1/r.
\]
The exact value of $M(\Omega_3)$ can be obtained also in terms of the real constant $s$ as in the following formula from~\cite[Eq.~(5.52) and Exercise 5.60(1)]{vu88}, 
\[
M(\Omega_3)=\frac{\pi}{2\mu(1/\sqrt{1+s})}
\]
were $\mu(s)$ is the Gr\"otzsch modulus function~\cite[Chapter~5]{avv}
\[
\mu(s)=\frac{\pi}{2}\frac{K(s')}{K(s)}, \quad K(s)=\int_{0}^{1}\frac{dx}{\sqrt{(1-x^2)(1-s^2x^2)}}, \quad s'=\sqrt{1-s^2}.
\]
Consequently, the exact value of the constant $r$ is given by
\begin{equation}\label{eq:r-mu}
r=2\mu(1/\sqrt{1+s})/\pi
\end{equation}
where the value of the $s+1$ is given by~\eqref{eq:r-sp1}. In this paper, the values of the function $\mu$ are computed as described in~\cite{NV-ccc}.

\begin{figure}[ht] %
\centerline{
\scalebox{0.35}{\includegraphics[trim=0cm 0cm 0cm 0cm,clip]{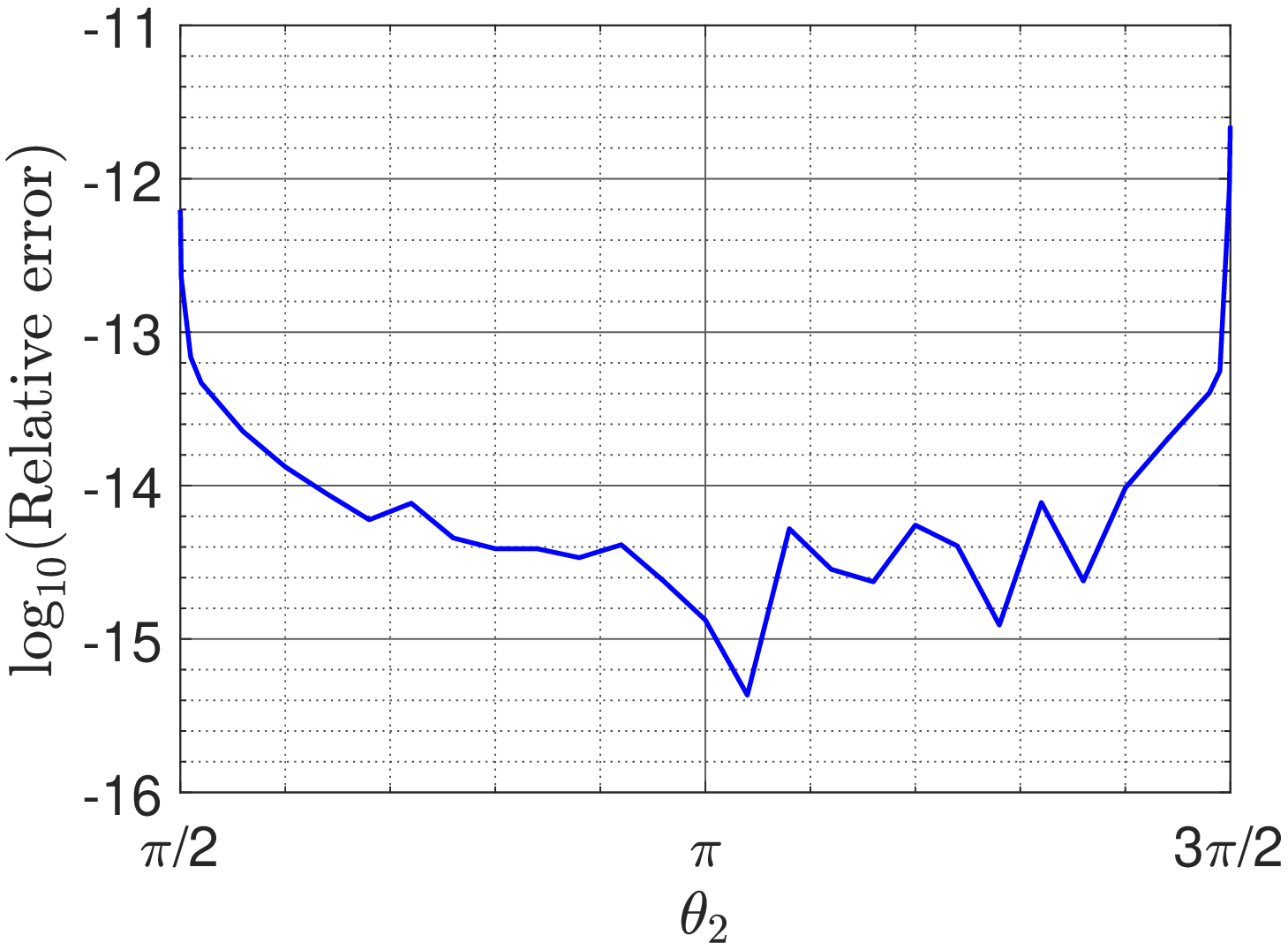}}
\hfill
\scalebox{0.35}{\includegraphics[trim=0cm 0cm 0cm 0cm,clip]{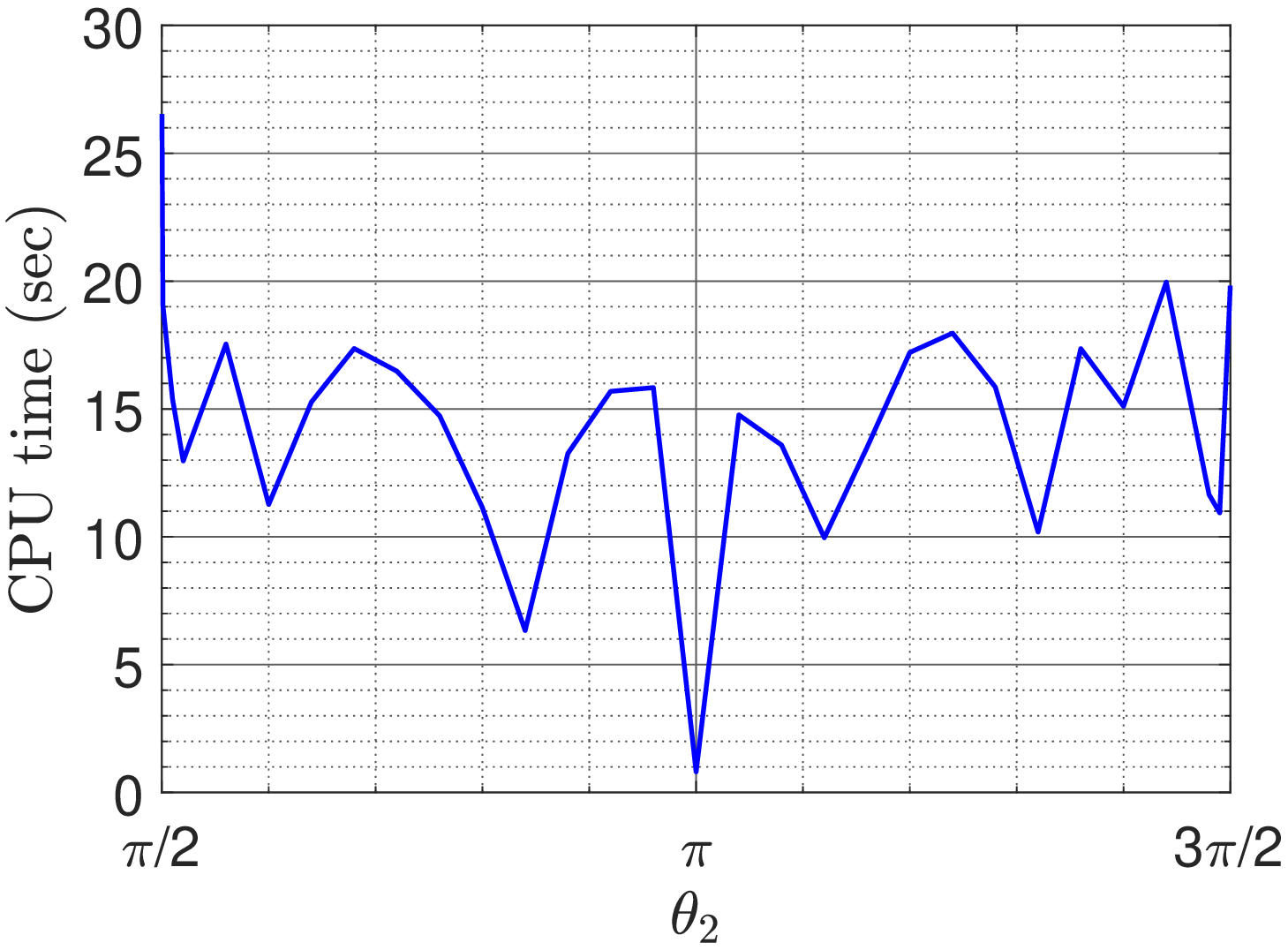}}
\hfill
\scalebox{0.35}{\includegraphics[trim=0cm 0cm 0cm 0cm,clip]{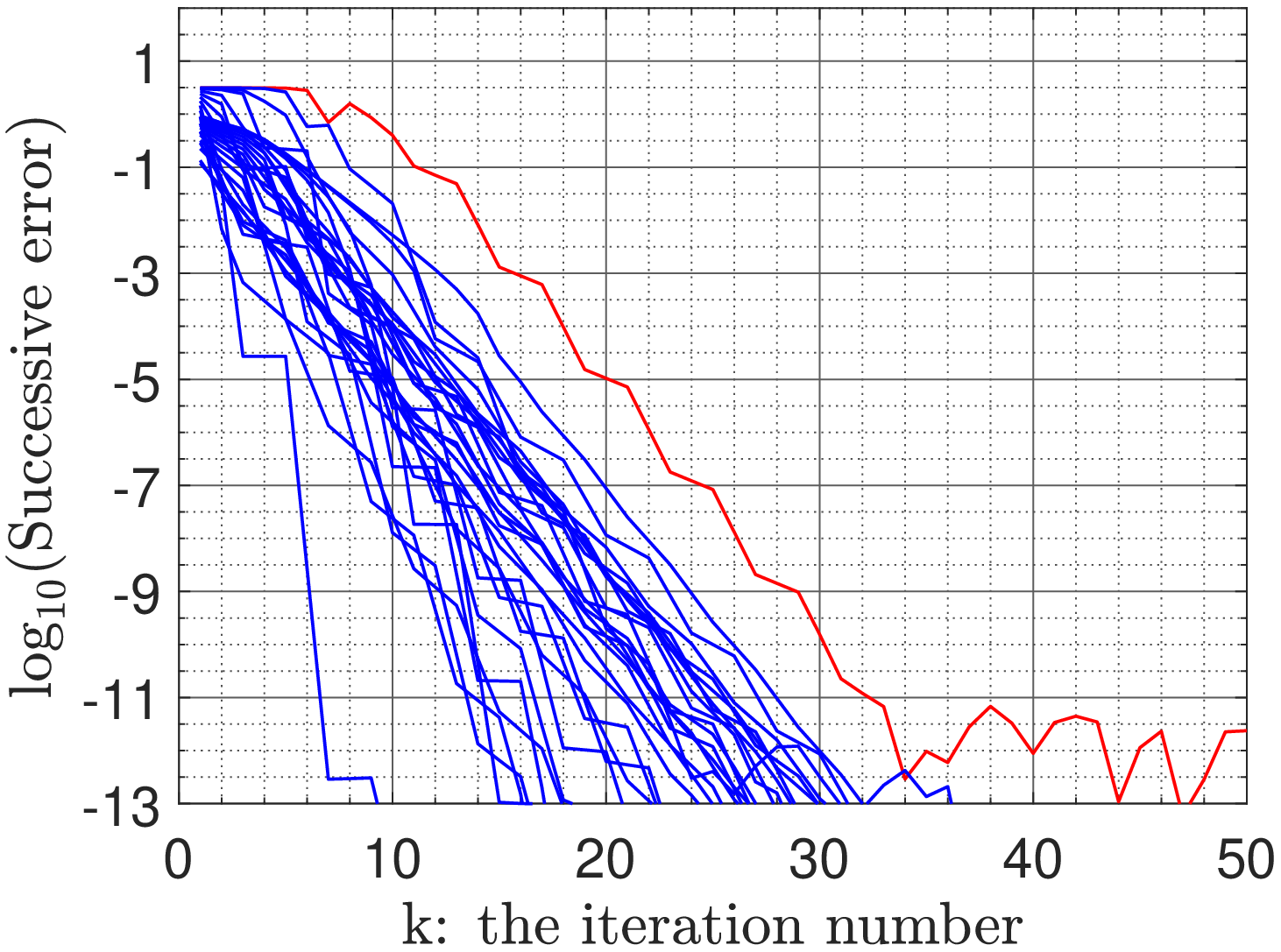}}
}
\caption{Comparison with the exact formula~\eqref{eq:r-mu} for $n=2^{13}$.}
\label{fig:err-ex}
\end{figure}

\begin{figure}[ht] %
\centerline{
\scalebox{0.35}{\includegraphics[trim=0cm 0cm 0cm 0cm,clip]{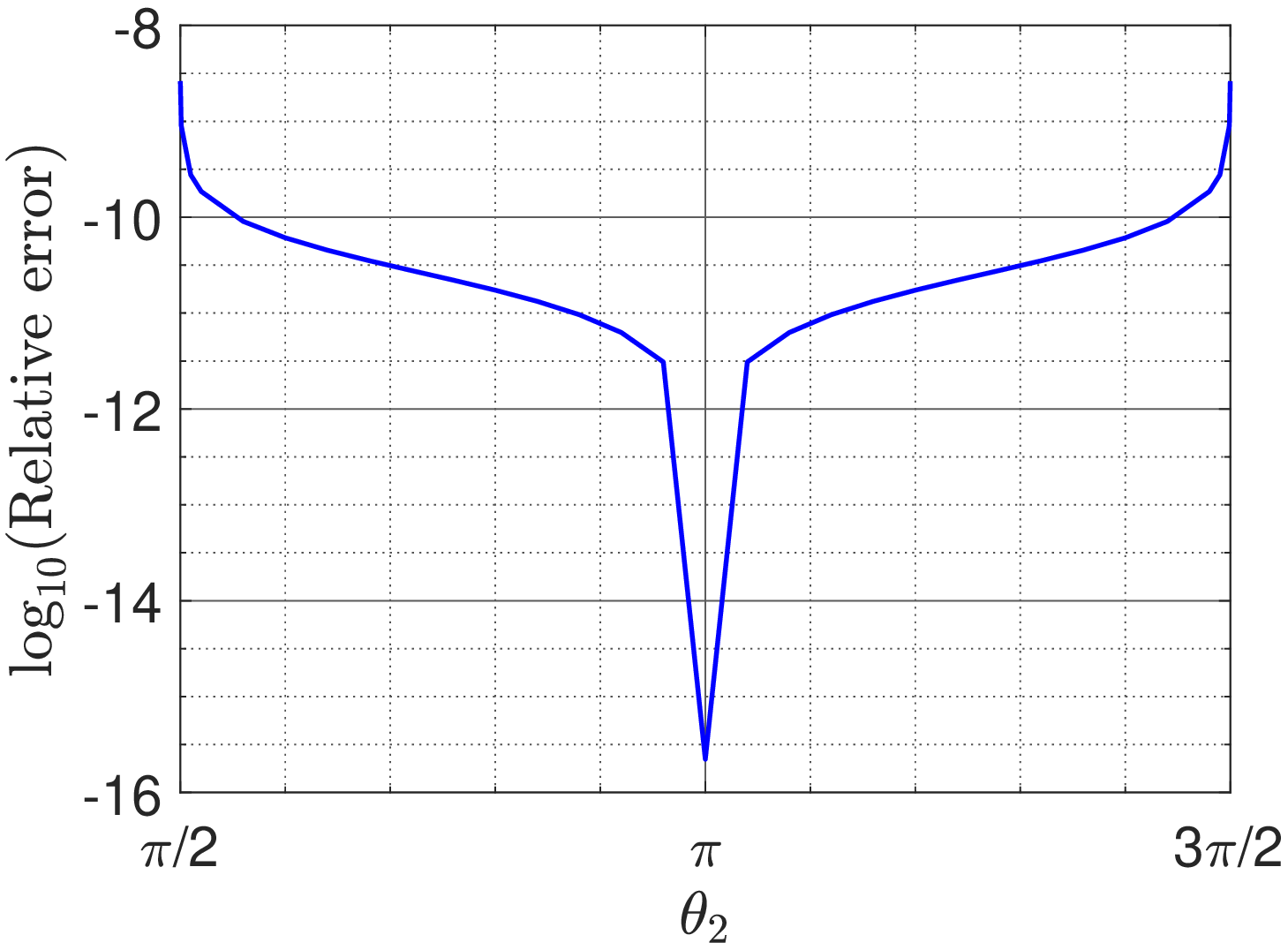}}
\hfill
\scalebox{0.35}{\includegraphics[trim=0cm 0cm 0cm 0cm,clip]{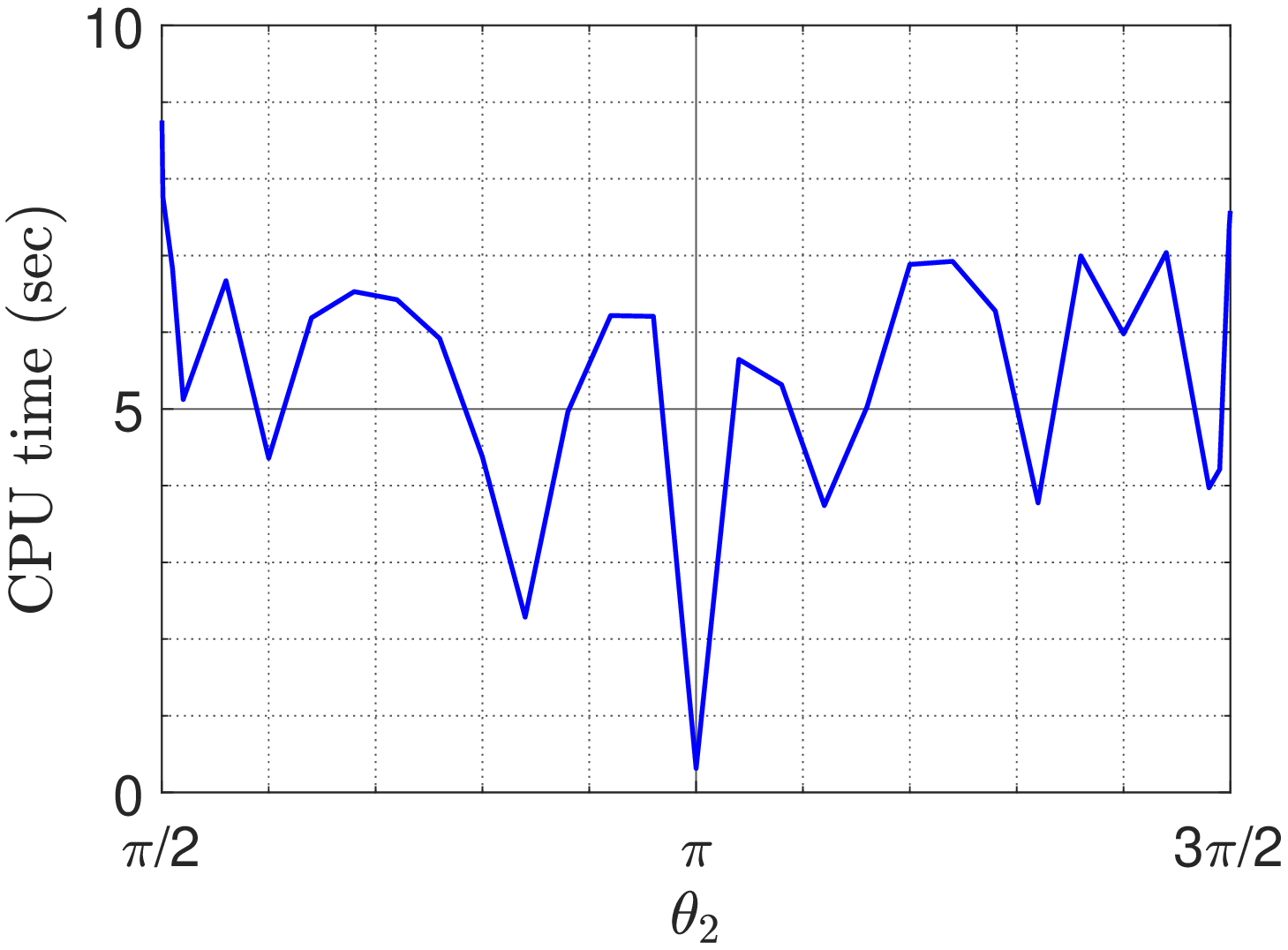}}
\hfill
\scalebox{0.35}{\includegraphics[trim=0cm 0cm 0cm 0cm,clip]{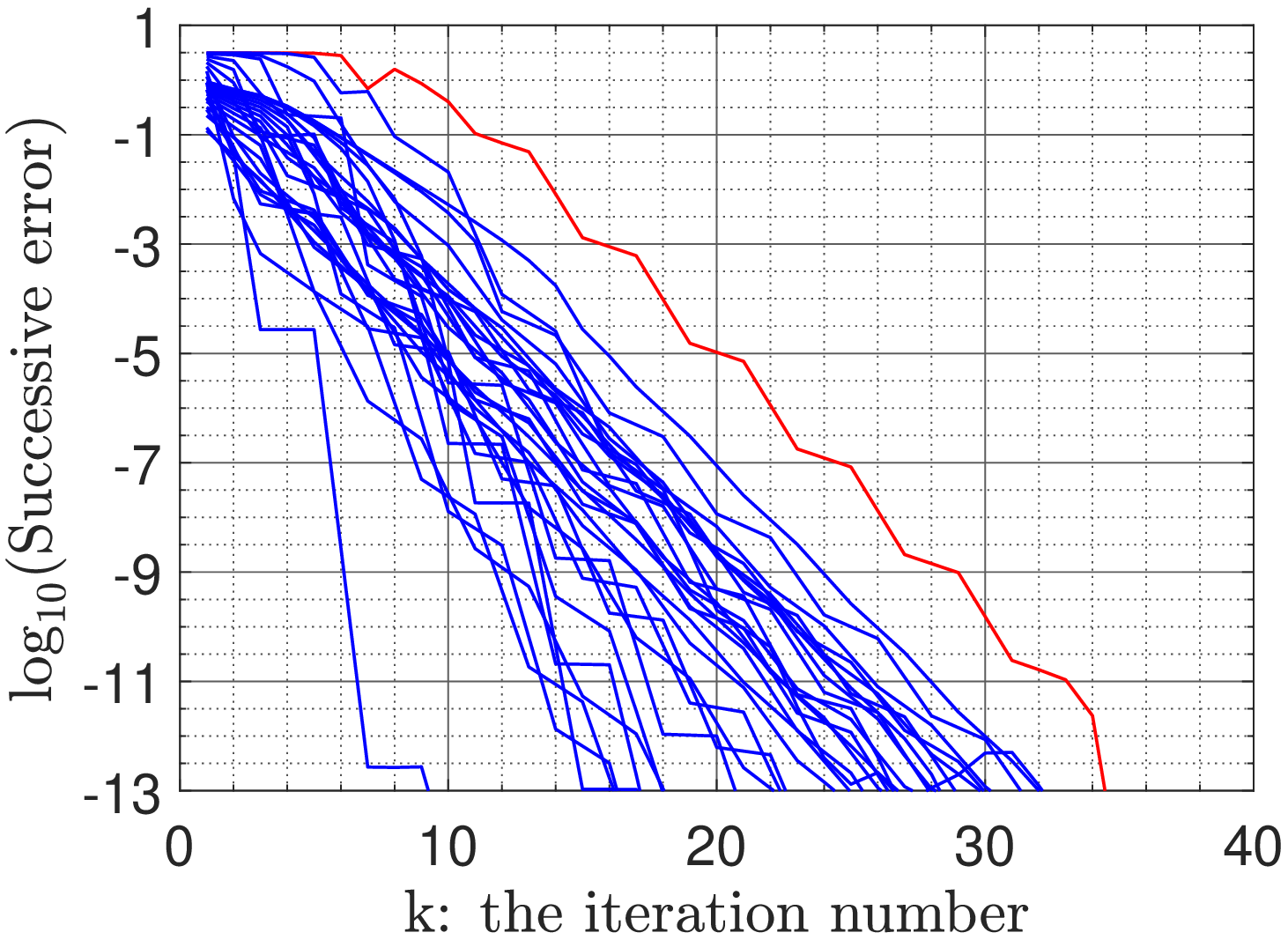}}
}
\caption{Comparison with the exact formula~\eqref{eq:r-mu} for $n=2^{10}$.}
\label{fig:err-ex10}
\end{figure}

To test the Algorithm~\ref{sc:algorithm}, we fix $\theta_1=0.5\pi$ and $\theta_3=1.5\pi$. Then, we choose values for $\theta_2$ between $0.5001\pi$ and $1.4999\pi$. The  numerical results obtained with $n=2^{13}$ are shown in Figure~\ref{fig:err-ex}. Figure~\ref{fig:err-ex} shows the relative error in the computed values of $r$ for vs $\theta_2$ (left), the total CPU time (in seconds) required for computing each value of $r$ vs $\theta_2$ (center), and the successive error $|r_{k}-r_{k-1}|$ for each value of $\theta_2$ vs the iteration number $k$ (right). We see from the figure, with less than $40$ iterations, the successive error for the iterative method method is less than $10^{-13}$ for all values of $\theta_2$ except for $\theta_2=0.5001\pi$ (red line). As expected, the relative error in the computed values of $r$ is very small when $\theta_2$ is a way from $\theta_1$ and $\theta_3$. The numerical results obtained with $n=2^{10}$ are shown in Figure~\ref{fig:err-ex10}.

\nonsec{\bf The crowding phenomenon.}
According to~\cite[pp.~20-21]{dt}, the term \emph{crowding} was coined in 1980~\cite{mz80} to describe the error/instability in numerical computing of conformal mapping. Thereafter it has become like a ``benchmark issue'' for all numerical conformal mapping software. 
As explained in~\cite[p.~77]{ps}, mapping a rectangle with aspect ratio $m$ to the unit disk seems to be impossible for $m=24$. Problems start already with $m=8$ and become more serious with increasing $m$. The critical value of $m$ depends on the computer 
floating point arithmetic with $10<m<20$ for double precision arithmetics \cite[pp.20-21]{dt}, \cite[pp.75-77]{ps}.

For rectangle $R_r$ in~\eqref{eq:R_r}, the aspect ratio is $m=r$ for $r>1$ and $m=1/r$ for $r<1$. Assume that $\theta_1=\pi/2$ and $\theta_3=3\pi/2$ are fixed as above and $\pi/2<\theta_2<3\pi/2$. 
In this subsection, we use the analytic example presented in Subsection~\ref{ex:explicit} to find the critical value of $r$ for mapping the quadrilateral $(D;1,\i,e^{\i\theta_2},-\i)$ onto the rectangle $R_r$ with double precision arithmetics.

\begin{figure}[ht] %
\centerline{
\scalebox{0.5}{\includegraphics[trim=0cm 0cm 0cm 0cm,clip]{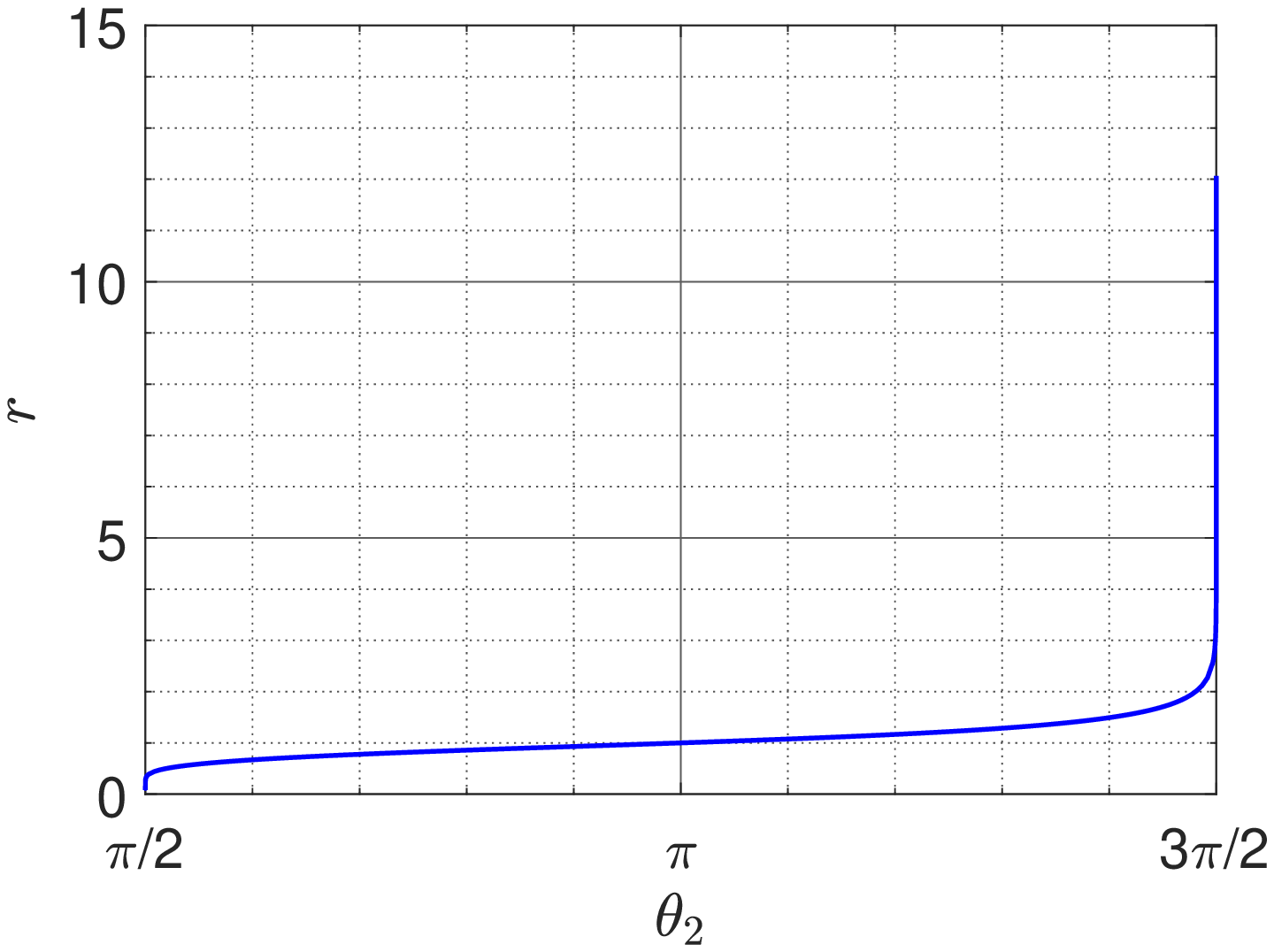}}
\hfill
\scalebox{0.5}{\includegraphics[trim=0cm 0cm 0cm 0cm,clip]{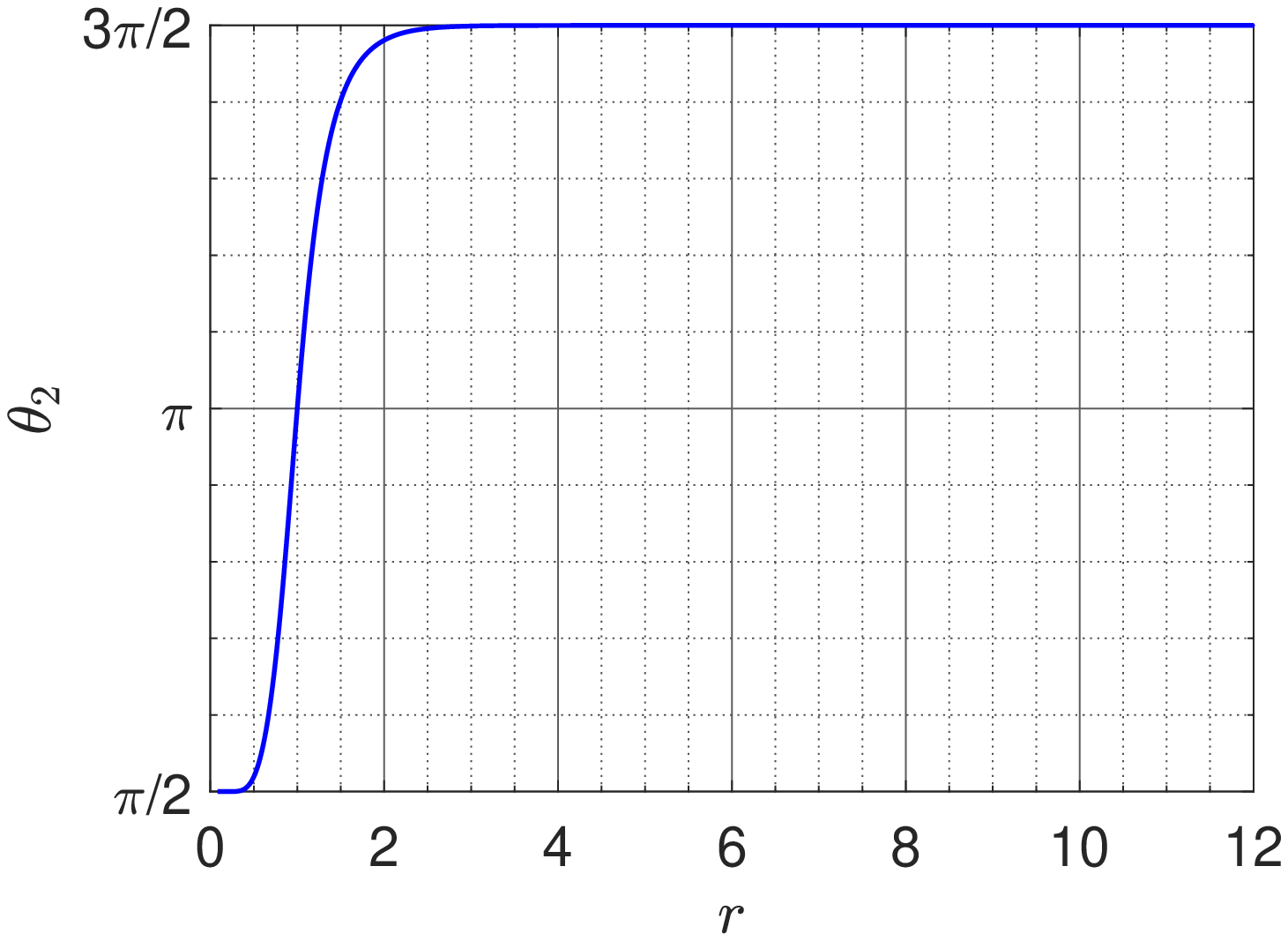}}
}
\caption{The values of $r$ in terms of $\theta_2$ (left) and the values of $\theta_2$ in terms of $r$ (right).}
\label{fig:r-thet2}
\end{figure}

In view of~\eqref{eq:r-sp1}, we have
\[
1+s=\frac{2}{1+\cot(\theta_2/2)}.
\]
Thus, by~\eqref{eq:r-mu} the value of the modulus $r$ can be written in terms of $\theta_2$,
\begin{equation}\label{eq:r-thet2}
r=\frac{2}{\pi}\mu\left(\sqrt{\frac{1+\cot(\theta_2/2)}{2}}\right).
\end{equation}
Also, the value of $\theta_2$ can be written in terms of the modulus $r$,
\begin{equation}\label{eq:thet2-r}
\theta_2=2\cot^{-1}\left(2\left(\mu^{-1}(r\pi/2)\right)^2-1\right).
\end{equation}
The values of the modulus $r$ obtained with the formula~\eqref{eq:r-thet2} for $\pi/2+10^{-15}<\theta_2<3\pi/2-10^{-15}$ are shown in Figure~\ref{fig:r-thet2} (left). Similarly, Figure~\ref{fig:r-thet2} (right) shows values $\theta_2$ obtained with the formula~\eqref{eq:thet2-r} for $1/12<r<12$. We see from Figure~\ref{fig:r-thet2} (right) that the values of $\theta_2$ become very close to $\theta_3=3\pi/2$ even for small values of $r$. In fact, for $r=12$, the value of $\theta_2$ obtained with formula~\eqref{eq:thet2-r} satisfies
\[
3\pi/2-\theta_2=1.776\times10^{-15}.
\]
Similarly, for $r=1/12$, we have
\[
\theta_2-\pi/2=1.776\times10^{-15}.
\]

\begin{figure}[ht] %
\centerline{
\scalebox{0.5}{\includegraphics[trim=0cm 0cm 0cm 0cm,clip]{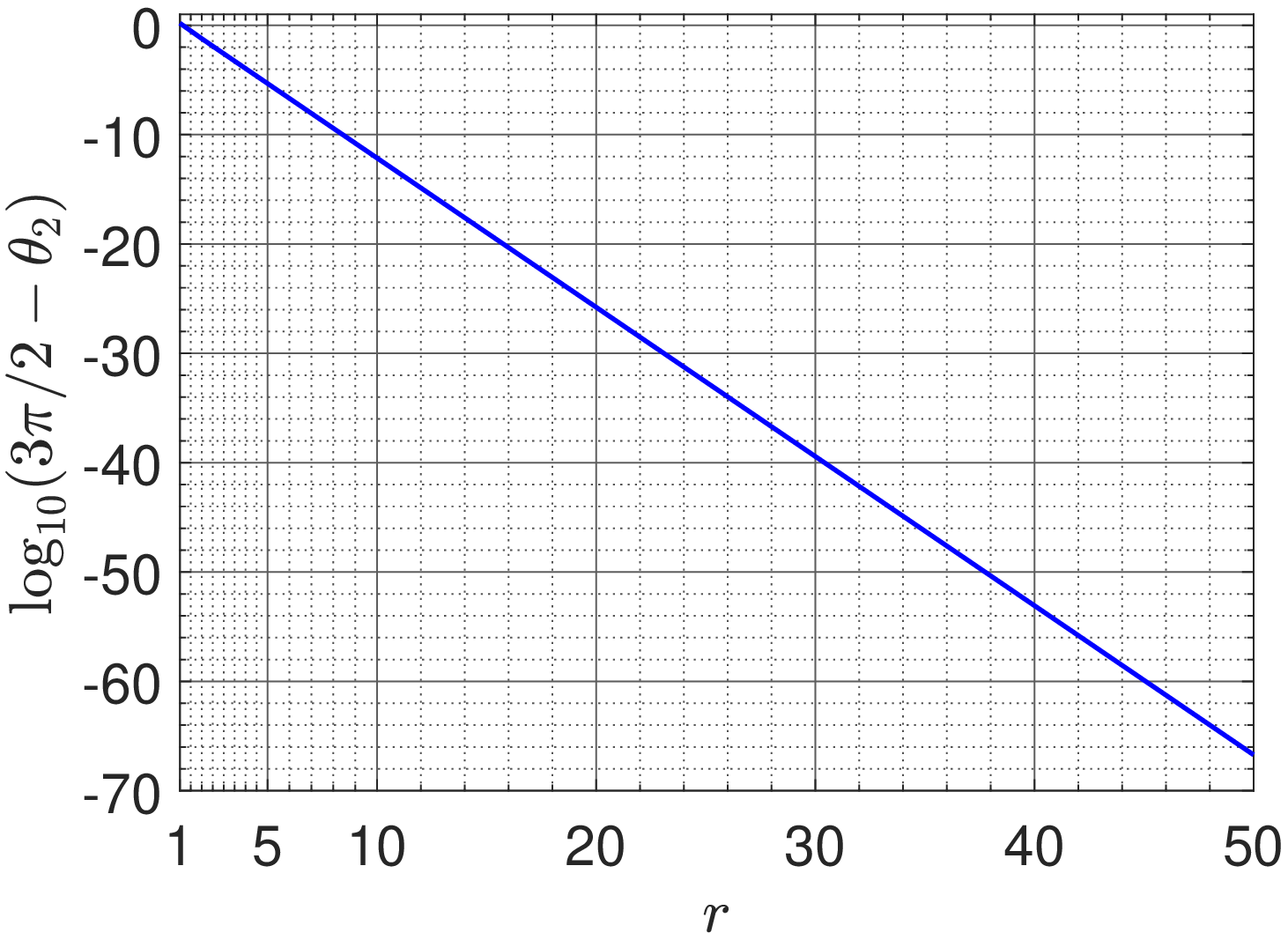}}
\hfill
\scalebox{0.5}{\includegraphics[trim=0cm 0cm 0cm 0cm,clip]{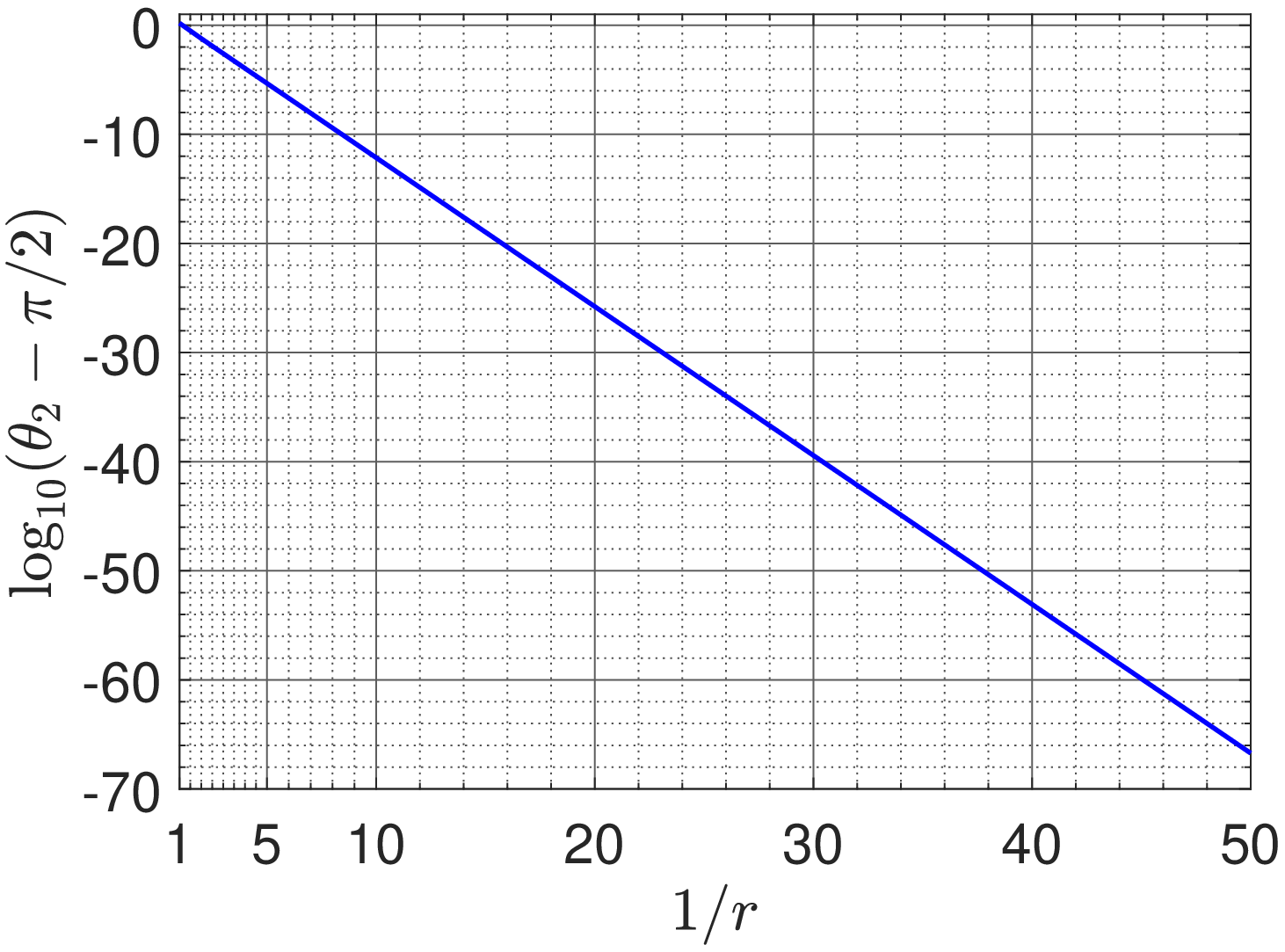}}
}
\caption{The relation between $r$ and $\theta_2$ for $r>1$ (left) and $0<r<1$ (right).}
\label{fig:r-thet2-sym}
\end{figure}

For further investigation, we use MATLAB symbolic toolbox to obtain more accurate results for the formulas~\eqref{eq:r-thet2} and~\eqref{eq:thet2-r}. For $r\ge1$, the values of $\log_{10}(3\pi/2-\theta_2)$ are shown in Figure~\ref{fig:r-thet2-sym} (left). 
It is clear from Figure~\ref{fig:r-thet2-sym} (left) that, up to double precision accuracy of the computer, $\theta_2=3\pi/2$ for the value of $r$ as small as $r=13$. 
Similarly, for $0<r\le1$, the values of $\log_{10}(\theta_2-\pi/2)$ are shown in Figure~\ref{fig:r-thet2-sym} (right). Up to double precision accuracy of the computer, it follows from Figure~\ref{fig:r-thet2-sym} (right) that $\theta_2=\pi/2$ for $r=1/13$. As a consequence, up to double precision accuracy of the computer, mapping a rectangle onto a quadrilateral $(D;1,\i,e^{\i\theta_2},-\i)$ is impossible when the aspect ratio of the rectangle is as small as $r=13$.

We see from Figure~\ref{fig:r-thet2-sym} (left) that the relation between $r$ and $\log_{10}(3\pi/2-\theta_2)$ is linear for $r>1$. We use MATLAB function \verb|polyfit| to find the coefficients of the line, and hence, we can estimate
\begin{equation}\label{eq:thet2-app-1}
\theta_2(r)\approx \frac{3\pi}{2}-32.3663566817311\times10^{-1.36452159123521r}, \quad r>1. 
\end{equation}
Similarly, Figure~\ref{fig:r-thet2-sym} (right) shows that there is a linear relationship between $1/r$ and $\log_{10}(\theta_2-\pi/2)$ for $0<r<1$. By using MATLAB function \verb|polyfit| to find the coefficients of the line, we can estimate
\begin{equation}\label{eq:thet2-app-2}
\theta_2(r)\approx \frac{\pi}{2}+32.3665310118084\times10^{-1.36452172896714/r}, \quad 0<r<1. 
\end{equation}
Equation~\eqref{eq:thet2-app-1} illustrates how fast the value of $\theta_2$ approaches $\theta_3=3\pi/2$ even for small values of $r$, $r>1$. Similarly, Equation~\eqref{eq:thet2-app-2} illustrates that the value of $\theta_2$ approaches $\theta_1=\pi/2$ even for small values of $1/r$, $0<r<1$.
